\documentclass{article}
\usepackage{fullpage}
\usepackage[utf8]{inputenc}
\usepackage{amsmath}
\usepackage{algorithm,algorithmic}
\usepackage{psfrag}
\usepackage{amsfonts}
\usepackage{amssymb}
\usepackage{amsthm}
\usepackage{color}
\usepackage{color}
\usepackage{graphicx,subfigure}
\usepackage{array}
\usepackage{tabularx}
\usepackage{longtable}
\usepackage{ulem}
\usepackage{cancel}
\usepackage{multicol}
\usepackage{url}
\newtheorem{lemma}{Lemma}

\newcommand{\beq}{\begin{equation}}
\newcommand{\eeq}{\end{equation}}

\newcommand{\dpar}[2]{\dfrac{\partial #1}{\partial #2}}

\usepackage{draftcopy}

\newcommand{\bff}{\mathbf{f}}
\newcommand{\bA}{\mathbf{A}}
\newcommand{\bF}{\mathbf{F}}
\newcommand{\bR}{\mathbf{R}}
\newcommand{\bL}{\mathbf{L}}

\DeclareMathOperator{\diver}{div}

\newcommand{\bn}{\mathbf{n}}
\newcommand{\bx}{\mathbf{x}}

\definecolor{etonblue}{rgb}{0.59, 0.78, 0.64}
\definecolor{airforceblue}{rgb}{0.36, 0.54, 0.66}

\newtheorem{remark}{Remark}

\usepackage{hyperref}
\hypersetup{
pdftitle={High-order residual distribution scheme\\ for the time-dependent Euler equations of fluid dynamics}
pdfauthor={R. Abgrall, Paola Bacigaluppi and Svetlana Tokareva},
pdfpagemode=UseOutlines,
linkbordercolor=0 0 0,
linkcolor=red,
citecolor=blue,
colorlinks=true,
bookmarks = true
}

\newcommand{\rev}[1]{#1}

\begin{document}

   \title{High-order residual distribution scheme\\ for the time-dependent Euler equations of fluid dynamics}
\author{R. Abgrall, P. Bacigaluppi, S. Tokareva\\
Institute of Mathematics, University of Zurich\\ Winterthurerstrasse 190, 8057 Z\"urich, Switzerland\\
$\{$remi.abgrall,paola.bacigaluppi,svetlana.tokareva$\}$@math.uzh.ch}
\date{April 15th, 2018}
\maketitle

\begin{abstract}
In the present work, a high order finite element type residual distribution scheme is designed in the framework of multidimensional compressible Euler equations of gas dynamics. The strengths of the proposed approximation rely on the generic spatial discretization of the model equations using a continuous finite element type approximation technique, while avoiding the solution of a large linear system with a sparse mass matrix which would come along with any standard ODE solver in a classical finite element approach to advance the solution in time. In this work, we propose a new Residual Distribution (RD) scheme, which provides an arbitrary explicit high order approximation of the smooth solutions of the Euler equations both in space and time.
\rev{The design of the scheme via the coupling of the RD formulation \cite{mario,abg} with a Deferred Correction (DeC) type method \cite{shu-dec,Minion2}, allows to have the matrix associated to the update in time, which needs to be inverted, to be diagonal. The use of Bernstein polynomials as shape functions, guarantees that this diagonal matrix is invertible and ensures strict positivity of the resulting diagonal matrix coefficients. }
This work is the extension of \cite{enumath,Abgrall2017} to multidimensional systems.
We have assessed our method on several challenging benchmark problems for one- and two-dimensional Euler equations and the scheme has proven to be robust and to achieve the theoretically predicted high order of accuracy on smooth solutions. 
\end{abstract}

\section{Introduction}
Consider a generic multidimensional time-dependent hyperbolic system of equations
\begin{equation}\label{system1}
\dpar{U}{t}+\diver\bF(U)=0 
\end{equation}
defined on a space-time domain $\Omega \times T$, with necessary initial and boundary conditions.
We are interested in a numerical approximation of \eqref{system1} by means of a finite element (FE) type technique. 
In \cite{enumath,Abgrall2017}, we have shown how one can solve a scalar version of \eqref{system1} with a method that approximates the spatial term using a Residual Distribution (RD) approach, without having to solve a large linear system with a sparse mass matrix. This means  that we are able to avoid any mass matrix "inversion" and have also an explicit scheme. This is achieved by first approximating the time operator in a consistent way with the spatial term. A priori, this would lead either to an implicit method in case of a nonlinear approximation, as done in order to avoid spurious oscillations in the case of discontinuous solutions, or at minima the inversion of a sparse but non diagonal matrix.  This apparent difficulty can be solved by applying a Deferred Correction (DeC) type time-stepping method and the use of proper basis functions. It is demonstrated in \cite{enumath,Abgrall2017} that Bernstein polynomials are a suitable choice, but this is not the only possible one. \rev{The idea to use as shape functions the Bernstein polynomials, instead of the more typical Lagrange polynomials, has been discussed in \cite{Abgrall2010, Abgrall2017} applied to the context of high order residual distribution schemes and very recently, in \cite{Lohmann2017}, this idea has been applied for a different class of methods, namely, the flux-corrected transport method.}

The purpose of this paper is to show how these ideas can be further extended for solving the Euler equations of fluid dynamics for the simulation of flows involving strong discontinuities.
The RD formulation used here is based on the finite element approximation of the solution as a globally continuous piecewise polynomial.
\rev{The design principle of the new RD scheme guarantees a compact approximation stencil even for high order accuracy, which would hold for Discontinuous Galerkin (DG), but not for example for Finite Volume (FV) methods and allows to consider a smaller number of nodes than DG (\cite{enordhigh, Cangiani2013, AbgrallViville2017}).}

The format of this paper is the following. In Section~\ref{sec:RD-steady}, we recall the idea of the residual distribution schemes for steady problems and in Section~\ref{sec:DeC} we describe the time-stepping algorithm and adapt the method developed in \cite{enumath,Abgrall2017} to multidimensional systems. We illustrate the robustness and accuracy of the proposed method by means of rigorous numerical tests and discuss the obtained results in Section~\ref{sec:results}. Finally, we give the conclusive remarks and outline further perspectives.

\section{Basic ideas of residual distribution schemes}
\label{sec:RD-steady}

\subsection{Governing equations and approximate solution}

Let us consider a generic system of PDEs 
\begin{equation}
\begin{cases}
\dpar{U}{t} + \text{div\,} \mathbf{F}(U) =0\quad \text{on}\; \;\Omega \times [0,T]\\[0.1em]
U(x,0) = U^0(x),\\
\end{cases}
\label{system}
\end{equation}
which could, for example, represent the Euler equations of gas dynamics,
with $U=[\rho,\; \rho u,\; \rho v,\; E]^T$, with the fluxes $\bF = (\bff_1, \bff_2)$ defined as
\begin{align*}
& \bff_1 = (\rho u,\; \rho u^2 + p,\; \rho u v,\; u(E+p))^T, \\
& \bff_2 = (\rho v,\; \rho u v,\; \rho v^2 + p,\; v(E+p))^T 
\end{align*}

We discretize the system \eqref{system} using the residual distribution approach. In this section, we shall give a brief overview of the RD method for steady problems and discuss the spatial discretization. The reader may refer to \cite{SWjcp,abg,Ricchiuto2007} for further details on the construction of generic residual distribution schemes. 

We consider the spatial domain $\Omega$ and its triangulation $\Omega_h$, and denote by $K$ a generic element of the mesh and by $h$ the characteristic mesh size. We also introduce the time discretization with time steps $\Delta t_n = t_{n+1} - t_n$.
\begin{figure}[H]
\centering
\includegraphics[scale=0.5]{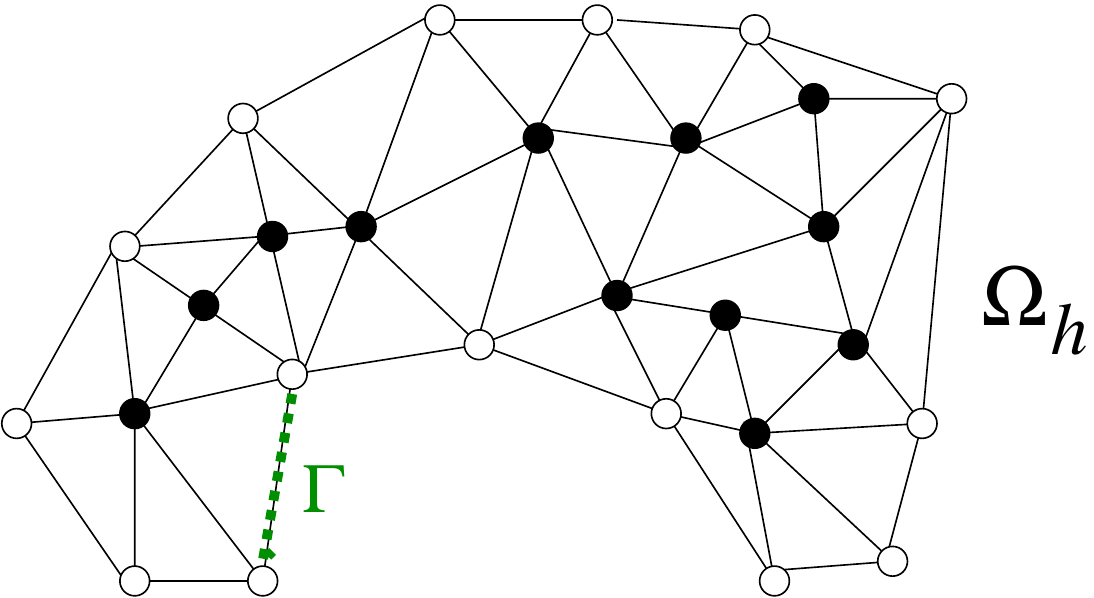}
\caption{Discretized domain $\Omega_h$ and its boundary $\Gamma$}
\label{Step1_RD}
\end{figure}

Following the ideas of the Galerkin finite element method (FEM), the solution approximation space $V_h$ is given by globally continuous piecewise polynomials of degree $k$:
\begin{equation}
V_h=\{U \in L^2(\Omega_h) \cap C^{0}(\Omega_h), U_{|K} \in \mathcal{P}^{k}, \forall K \in \Omega_h\},
\end{equation}
so that the numerical solution $U_h^{n} \approx U(\mathbf{x},t_{n})$ can be written as a linear combination of shape functions $\varphi_{\sigma} \in V_h$:
\begin{equation}
\label{approx_uh}
U_h^{n}=\sum_{\sigma \in {\Omega_h}} U_\sigma^{n} \varphi_\sigma = \sum_{K \in {\Omega_h}} \sum_{\sigma \in K} U_\sigma^{n} \varphi_\sigma,
\end{equation}
with coefficients $U_\sigma^n$ to be calculated by a numerical method.

\subsection{Residual distribution scheme for steady problems}

Consider first a steady scalar version of system \eqref{system}: 
\begin{equation}
\nabla_x\cdot \mathbf{F}(U)=0.
\label{steadysys}
\end{equation}

\begin{figure}[H]
\centering
\subfigure[Step 1: Compute fluctuation]{\includegraphics[scale=0.3]{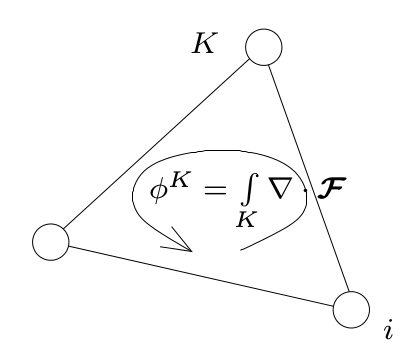}\label{Steps_RDa}}
\hspace{0.18cm}\subfigure[Step 2: Split distribution]{\includegraphics[scale=0.3]{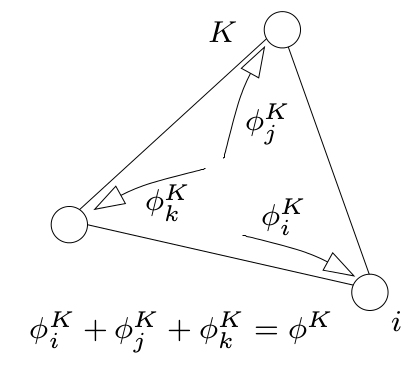}\label{Steps_RDb}}
\hspace{0.33cm}\subfigure[Step 3: Gather residuals, evolve]{\includegraphics[scale=0.4]{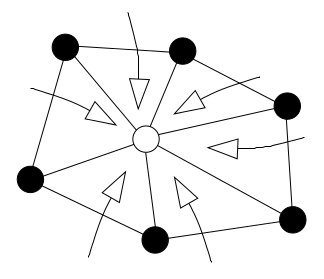}\label{Steps_RDc}}
\caption{Illustration of the three steps of the residual distribution approach \rev{for linear triangular elements}.}
\label{Steps_RD}
\end{figure}

The main steps of the residual distribution approach could be summarized as follows, \rev{see also Fig.~\ref{Steps_RD} where the approach is illustrated for linear FEM on triangular elements}:
\begin{enumerate}
\item
We define  $\forall K \in \Omega_h$ a fluctuation term (total residual) $\phi^K=\int_K \nabla_x\cdot \mathbf{F}(U)\,d\mathbf{x}$ (see Fig.~\ref{Steps_RDa})

\item
We define a nodal residual $\phi_\sigma^K$ as the contribution to the fluctuation term $\phi^K$ from a degree of freedom (DoF) $\sigma$ within the element $K$, so that the following conservation property holds (see Fig.~\ref{Steps_RDb}):
\begin{equation}
\phi^K(U_h)=\sum_{\sigma \in K} \phi_\sigma^K, \quad \forall K \in \Omega_h, \quad\forall  \Omega_h
\label{RD_distrib}
\end{equation}

The distribution strategy, i.e. how much of the fluctuation term has to be taken into account on each DoF $\sigma \in K$, is defined by means of the so-called distribution coefficients $\beta_\sigma$:
\begin{equation}
\phi_\sigma^K=\beta_\sigma^K \; \phi^K,
\label{def_phiK}
\end{equation}
where, due to \eqref{RD_distrib},
\begin{equation*}
\sum_{\sigma \in K} \beta_\sigma^K = 1.
\end{equation*}

\item The resulting scheme is obtained by collecting all the residual contributions $\phi_\sigma^K$ from elements $K$ surrounding a node $\sigma \in \Omega_h$ (see Fig.~\ref{Steps_RDc}), that is
\begin{equation}
\sum_{K | \sigma \in K} \phi_\sigma^K = 0, \quad \forall \sigma \in \Omega_h,
\label{RD_nodal}
\end{equation}
which allows to calculate the coefficients $U_\sigma$ in the approximation \eqref{approx_uh}.
\end{enumerate}

\subsection{On the choice of the spatial discretization}

In \cite{abgrall2017jcp,Abgrall2017,Abgrall99} it has been shown that any known finite element or finite volume scheme (such as SUPG, DG, FV-WENO, etc.) can be written in a generic residual distribution form \eqref{RD_nodal}.
In case $\sigma \in \Gamma$, equation \eqref{RD_nodal} \rev{can be split} for any degree of freedom (DoF) $\sigma$ into the internal and boundary contributions:
\begin{equation}
\label{RD_nodaldecomp}
\sum\limits_{K \subset \Omega_h,\sigma\in K} \phi_{\sigma,\bx}^{K}(U_h)+\sum\limits_{\gamma\subset \Gamma, \sigma\in \gamma}\phi_{\sigma,\bx}^{\gamma}(U_h)=0,
\end{equation}
where $\gamma$ is an edge on the boundary $\Gamma$ of the computational domain $\Omega_h$. 
The values $\phi_{\sigma,\bx}^{K}$ and $\phi_{\sigma,\bx}^{\gamma}$ are the residuals corresponding only to the spatial discretization, which is emphasized by the subscript $\mathbf{x}$.
Assuming that $u=g$ on $\Gamma$, both residuals satisfy the following conservation relations
\begin{equation}
\begin{split}
& \sum_{\sigma \in K}\phi_{\sigma,\bx}^K (U_h) = \int_{\partial K} \bF(U_h) \cdot \bn, \quad \forall K\\[0.3em]
& \sum_{\sigma\in \Gamma} \phi_{\sigma,\bx}^{\Gamma}(U_h)=\int_{\Gamma} (\mathcal{F}_{\bn}(U_h,g)-\bF(U_h)\cdot\bn), \quad \forall \Gamma.
\end{split}
\label{RD_nodaldecomp_contrib}
\end{equation}

Below we outline some of the schemes written in terms of residuals which satisfy the conservation relations \eqref{RD_nodaldecomp_contrib}, see also \cite{Abgrall2017}:
\begin{itemize}
    \item the SUPG scheme \cite{hughes}:
    \begin{multline}
    \label{SUPG}
    \phi_{\sigma,\bx}^K(U_h)=\int_{\partial K}\varphi_\sigma \bF(U_h)\cdot \bn\,d\gamma - \int_K \nabla \varphi_\sigma\cdot \bF(U_h)\,d\bx \\ + h_K
    \int_K \bigg (\nabla_U\,\bF(U_h)\cdot \nabla \varphi_\sigma \bigg )\tau \bigg (\nabla_U\,\bF(U_h)\cdot \nabla U_h \bigg )\,d\bx
    \end{multline}
    with $\tau>0$.
    \item the Galerkin FEM scheme with jump stabilization \cite{burman}:
    \begin{equation}
    \label{burman}
    \phi_{\sigma,\bx}^K(U_h)=\int_{\partial K}\varphi_\sigma \bF(U_h)\cdot \bn\,d\gamma -\int_K \nabla \varphi_\sigma\cdot \bF(U_h)\,d\bx +\sum_{e \in K}
    \theta h_e^2 \int_e \,[\nabla U_h]\cdot [\nabla \varphi_\sigma]\,d\gamma
    \end{equation}
    with $\theta\geq0$ (see \cite{burman} for details).
   Note that in this case if the mesh is conformal, any edge $e$ (or face in 3D) is the intersection of the element $K$ and another element denoted by $K^+$.
    For any function $\psi$, we define $[\nabla \psi ]=\nabla \psi_{|K}-\nabla \psi_{| K^+}$.  
    \item for the boundary approximation, it is possible to follow the same technique as in \cite{DeSantis2015}, so that for $\sigma \in \gamma \subset \Gamma$ we have 
$$ \phi_{\sigma}^{\gamma,\bx}(U_h) =\int_{\gamma} \Big( \mathcal{F}(U_h,g)-\bF(U_h)\cdot \bn\Big) d\gamma. $$     
\end{itemize}

It is also possible to consider schemes that do not have a straightforward variational formulation, as for example the limited residual distribution scheme (RDS) \cite{enordhigh,DeSantis2015,icm,CanadaCFD}:
\begin{equation}
\label{schema RDS SUPG}
    \phi_{\sigma,\bx}^K(U_h)=\beta_\sigma^K \int_{\partial K}\bF(U_h)\cdot \bn\,d\gamma + h_K
    \int_K \bigg (\nabla_U\,\bF(U_h)\cdot \nabla \varphi_\sigma \bigg )\tau \bigg (\nabla_U\,\bF(U_h)\cdot \nabla U_h \bigg )\,d\bx
\end{equation}
or
\begin{equation}
\label{schema_RDS_jump}
\phi_{\sigma,\bx}^K(U_h)=\beta_\sigma^K \int_{\partial K}\bF(U_h)\cdot \bn\,d\gamma + \sum_{\text{edges of }K}
\theta h_e^2 \int_e \,[\nabla U_h]\cdot [\nabla \varphi_\sigma]\,d\gamma.
\end{equation}
\rev{where $\beta_\sigma^K$ are parameters that guarantee conservation and $\sum_{\sigma \in K} \beta_{\sigma}^K=1.$}

\begin{remark}
\rev{
One may notice that in \eqref{schema RDS SUPG} and \eqref{schema_RDS_jump} the streamline diffusion term and jump term are introduced due to the possible existence of spurious modes in the solution, but their role is somehow different compared to \eqref{SUPG} and \eqref{burman} where they are introduced to stabilize the Galerkin schemes (see \cite{ENORD,enordhigh,Abgrall2017} for more details). Indeed, \eqref{schema RDS SUPG} without the streamline term and \eqref{schema_RDS_jump} without the jump term satisfy a discrete maximum principle, and adding this terms violates formally the maximum principle, while experimentally this violation results to be extremely small if not non-existent.}\\
It is important to remark that, at least formally, the exact solution cancels the residuals in the case of SUPG and RDS-SUPG, while in case of Burman's jump stabilization, we are able to rewrite the scheme as 
$$\phi_\sigma^K(U_h)=\int_K\psi_\sigma \diver\bF(U)\,d\bx + R_\sigma(U_h)$$
    with  
    $$R_\sigma=\sum_{\text{edges of }K}
    \theta h_e^2 \int_e\;[\nabla U_h]\cdot [\nabla \varphi_\sigma]\,d\gamma.$$
 where  $\sum\limits_{\sigma\in K} R_\sigma=0$. Here, $R_\sigma$ is not zero, except for the exact solution unless this solution has continuous normal gradients, see \cite{burman} for more details. 
\end{remark}

\section{An Explicit High Order timestepping approach}
\label{sec:DeC}

\subsection{Iterative timestepping method}

In the previous sections, we have shown how system \eqref{steadysys} can be discretized in terms of residual distributions approach. The main target of this paper is to extend this approximation to unsteady problems. Moreover, we aim to have a high order and explicit approximation method in time. In the rest of this paragraph, we rephrase \cite{Abgrall2017}, since the discussion on the scalar case extends in a straightforward manner to the system case.

Here we describe the timestepping algorithm that we use in combination with the RD discretization in space to achieve high order accuracy in time. We consider $M$ subintervals within each time step $[t_n, t_{n+1}]$, so that $t_n = t_{n,0} < t_{n,1} < \dots < t_{n,m} < \dots < t_{n,M} = t_{n+1}$.
Next, for each subinterval $[t_{n,m},t_{n,m+1}]$, we introduce the corrections $r=0,\dots,R$ and denote 
the solution at the $r$-th correction and the $m$-th substep $t_{n,m}$ as $U_h^{n,m,r}$ and the solution at $t_n$ by $U_h^n$. In addition, we define the solution vector
\[ U^{(r)} = (U_h^{n,1,r},...,U_h^{n,m,r},...,U_h^{n,M,r}). \]

We propose a timestepping method that can be interpreted as a deferred correction method and proceed as follows within the time interval $[t_n,t_{n+1}]$:
\begin{enumerate}
\item for $r=0$, set $U_h^{(0)} = (U_h^{n,1,0},\dots,U_h^{n,m,0},\dots,U_h^{n,M,0}) = (U_h^n,\dots,U_h^n,\dots,U_h^n)$;
\item for each correction $r>0$, knowing $U^{(r)}$, evaluate $U^{(r+1)}$ as the solution of 
\begin{equation}
\mathcal{L}^1(U^{(r+1)})=\mathcal{L}^1(U^{(r)})-\mathcal{L}^2(U^{(r)})
\label{HO_timestepping}
\end{equation}
\item set the solution $U_h^{n+1}=U_h^{n,M,R}.$
\end{enumerate}

Formulation \eqref{HO_timestepping} relies on a  Lemma which has been proven in  \cite{Abgrall2017}.
\begin{lemma}\label{ZeLem} If two operators $\mathcal{L}^1_\Delta$ and $\mathcal{L}^2_\Delta$ depending on a parameter $\Delta$ are such that:
\begin{enumerate}
\item There exists a unique $U^\star_\Delta$ such that $\mathcal{L}^2_\Delta(U^\star_\Delta)=0$;
\item There exists $\alpha_1>0$ independent of $\Delta$, such that for any $U$ and $V$ the operator $\mathcal{L}^1_\Delta$ is coercive, i.~e. 
\begin{equation}
\label{coercive}
\alpha_1 ||U-V||\leq ||\mathcal{L}^1_\Delta (U)-\mathcal{L}^1_\Delta (V)||;
\end{equation}
\item There exists $\alpha_2>0$ independent of $\Delta$, such that, for any $U$ and $V$
\begin{equation}
\label{error}
\bigg |\bigg | \big (\mathcal{L}^1_\Delta(U)-\mathcal{L}^2_\Delta(U)\big )-\big (\mathcal{L}^1_\Delta(V)-\mathcal{L}^2_\Delta(V)\big )\bigg |\bigg |\leq \alpha_2 \Delta ||U-V||.
\end{equation}
This last condition is nothing more than saying that the operator $\mathcal{L}^1_\Delta -\mathcal{L}^2_\Delta$ is uniformly Lipschitz continuous with Lipschitz constant
$\alpha_2 \Delta$.
\end{enumerate}
Then if $\nu =\frac{\alpha_2}{\alpha_1}\Delta <1$ the deferred correction method \eqref{HO_timestepping} is convergent, and after $R$ iterations the error is smaller than $\nu^R || U^{(0)} - U^\star_\Delta ||$.
\end{lemma}

The differential operators $\mathcal{L}^1$ and $\mathcal{L}^2$ will be defined in the following sections; the detailed error analysis can be found in \cite{Abgrall2017}.

\subsection{On the low order differential operator $\mathcal{L}^1$}\label{Sec_L1operator}

Our discretization in time relies on the fact that the system \eqref{system} can be formally integrated on $[0,t]$ as
$$U(\bx,t)=U(\bx,0)+\int_0^t \diver \bF(U(x,s))\,ds,$$
and the solution can be further approximated using a quadrature formula as
\begin{equation}
 U(\bx,t)\approx U(\bx,0)+t\;\sum_{l=0}^{r} \omega_l \diver\bF(U(\bx, s_l)),
\end{equation}
with the same conventions as in the ODE case in \cite{Abgrall2017}.

For any $\sigma\in K$, define $\mathcal{L}^1_\sigma$ as:
 \begin{equation}
 \label{L1_complete}
 \mathcal{L}^1_\sigma(U^{(r)})= \mathcal{L}^1_\sigma(U^{n,1,r}, \ldots , U^{n,M,r})=\begin{pmatrix}
|C_\sigma|\big ( U_\sigma^{n,M,r} - U_\sigma^{n,0} \big)\, + \sum\limits_{K|\sigma\in K} \displaystyle\int_{t_{n,0}}^{t_{n,M}}  \mathcal{I}_{0}\big(\phi_{\sigma,\mathbf{x}}^K(U^{(r)}),s\big) \,ds\\
\vdots\\
|C_\sigma| \big ( U_\sigma^{n,1,r} - U_\sigma^{n,0}\big )+ \sum\limits_{K|\sigma\in K} \displaystyle\int_{t_{n,0}}^{t_{n,1}}  \mathcal{I}_{0}\big(\phi_{\sigma,\mathbf{x}}^K(U^{(r)}),s\big)\,ds
\end{pmatrix}.
\end{equation}
where $\mathcal{I}_{0}$ represents any first order piecewise-constant interpolant and where we have adopted a notation
\[\mathcal{I}_{0}\big(\phi_{\sigma,\mathbf{x}}^K(U^{(r)}),s\big) = \mathcal{I}_{0}\big(\phi_{\sigma,\mathbf{x}}^{K}(U^{n,0,r}),\ldots,\phi_{\sigma,\mathbf{x}}^{K}(U^{n,M,r}),s\big).\] 

In order to simplify \eqref{L1_complete} and make the operator explicit in time, we take the interpolant $\mathcal{I}_{0}$ as a simple approximation at $U^{n,0}$, so that \eqref{L1_complete} becomes 
\begin{equation}
\mathcal{L}^1_\sigma(U^{(r)})= \mathcal{L}^1_\sigma(U^{n,1,r},\ldots , U^{n,M,r})=\begin{pmatrix}
 |C_\sigma|(U^{n,M,r}-U^{n,0})+\xi_M\;\Delta t\; \sum\limits_{K|\sigma\in K} \phi_{\sigma,\mathbf{x}}^K(U^{n,0})\\
\vdots\\
 |C_\sigma|(U^{n,1,r}-U^{n,0})+\xi_1\;\Delta t \; \sum\limits_{K|\sigma\in K} \phi_{\sigma,\mathbf{x}}^K(U^{n,0})\\\
\end{pmatrix},
\label{L1}
\end{equation}
\rev{where the weights $\xi_m$ with $m=1,..,M$ are chosen in $[t_n,t_{n+1}]$ to satisfy $t_{n,m}=t_{n}+\xi_m\Delta t$ and $0=\xi_0< ...<\xi_m<\xi_{m+1}<...<\xi_M=1$.\\}
In this system the coefficients $|C_\sigma|$ play the role of the dual cell measures and are defined as
\begin{equation}
|C_\sigma|:= \int_K \varphi_\sigma d\bx 
\label{Ci}
\end{equation}
A direct consequence of Lemma \ref{ZeLem} is that in order for \eqref{L1} to be solvable, we have to satisfy the constraint 
\begin{equation}
\label{Ci:constraint}
|C_\sigma|>0.
\end{equation}
This means that we are not free to choose any family of polynomials as shape functions but only those guaranteeing this property. For instance, for the family of the Lagrangian polynomials $\mathcal{P}^k$ the condition \eqref{Ci:constraint} doesn't hold for $k > 1$. \rev{Therefore, in this work we investigate the use of Bernstein polynomials for high-order residual distribution approximations.}\\
The low order differential operator $\mathcal{L}^1_\sigma$ constructed this way is explicit in time.\\

 The drawback of using Bernstein polynomials is that not all degrees of freedom $U_{\sigma}^n$ in the expansion \eqref{approx_uh} will represent the solution values at certain nodes, however, the advantage of this family of shape functions is their positivity on $K$ that will enforce \eqref{Ci:constraint}.

\rev{We next provide the expressions for the families of Bernstein polynomials used in this paper. On triangular elements, given the barycentric coordinates $x_1$, $x_2$, $x_3$, the Bernstein shape functions are defined as follows.}
\begin{multicols}{2}
\begin{itemize}
	\item \rev{Order 1 ('B1'): 
	\begin{equation*}
	\varphi_1 = x_1, \ \varphi_2 = x_2, \ \varphi_3 = x_3. \quad \quad \quad \quad
	\end{equation*}}
	
	\item \rev{Order 2 ('B2'):
	\begin{align*}
	& \varphi_1 = x_1^2, \ \varphi_2 = x_2^2, \ \varphi_3= x_3^2, \\
	& \varphi_4 = 2 x_1 x_2, \ \varphi_5 = 2 x_2 x_3, \ \varphi_6 = 2 x_1 x_3.
	\end{align*}}

	\item \rev{Order 3 ('B3'):
	\begin{align*}
	& \varphi_1 = x_1^3, \ \varphi_2 = x_2^3, \ \varphi_3= x_3^3, \\
	& \varphi_4 = 3 x_1^2 x_2, \ \varphi_5 = 3 x_1 x_2^2, \ \varphi_6 = 3 x_2^2 x_3, \\
	& \varphi_7 = 3 x_2 x_3^2, \ \varphi_8 = 3 x_1 x_3^2, \ \varphi_9 = 3 x_1^2 x_3, \\
	& \varphi_{10} = 6 x_1 x_2 x_3.
	\end{align*}
	}
\end{itemize}
	
	\begin{figure}[H]
	\centering
\includegraphics[scale=0.3]{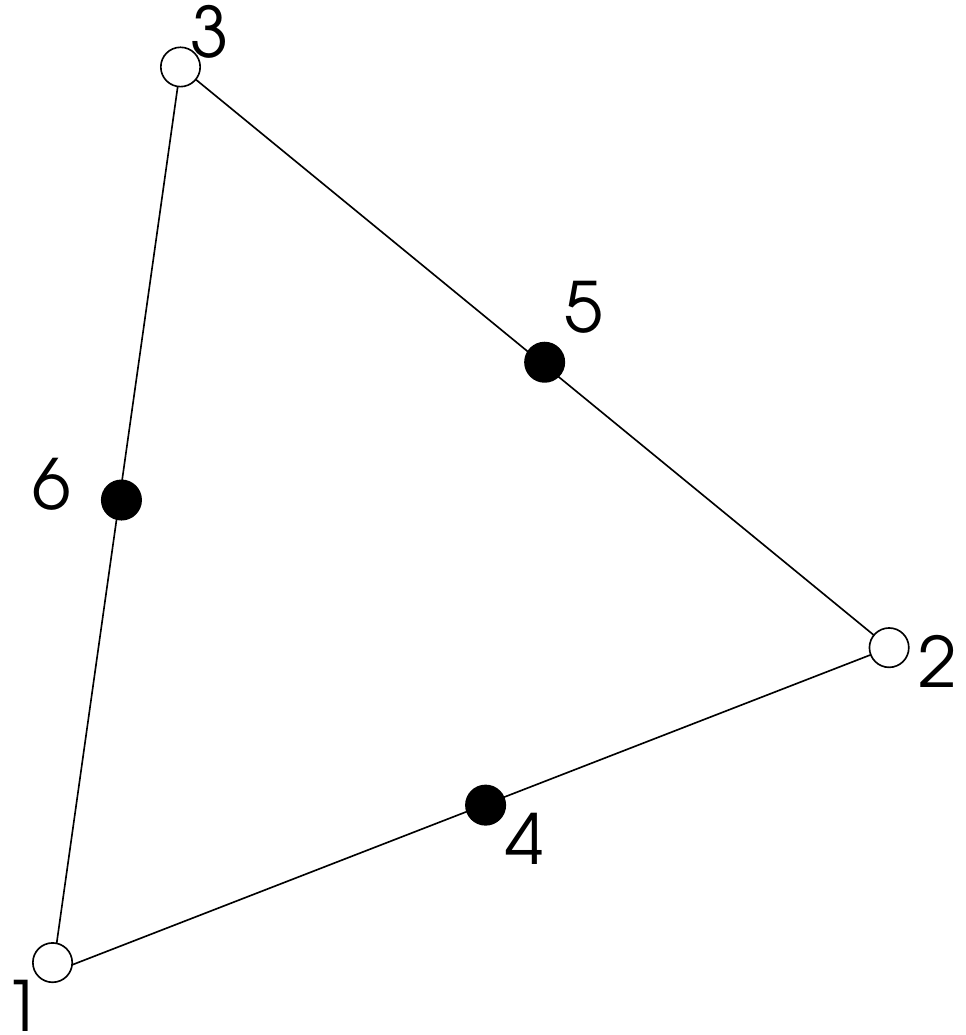}\\
\includegraphics[scale=0.3]{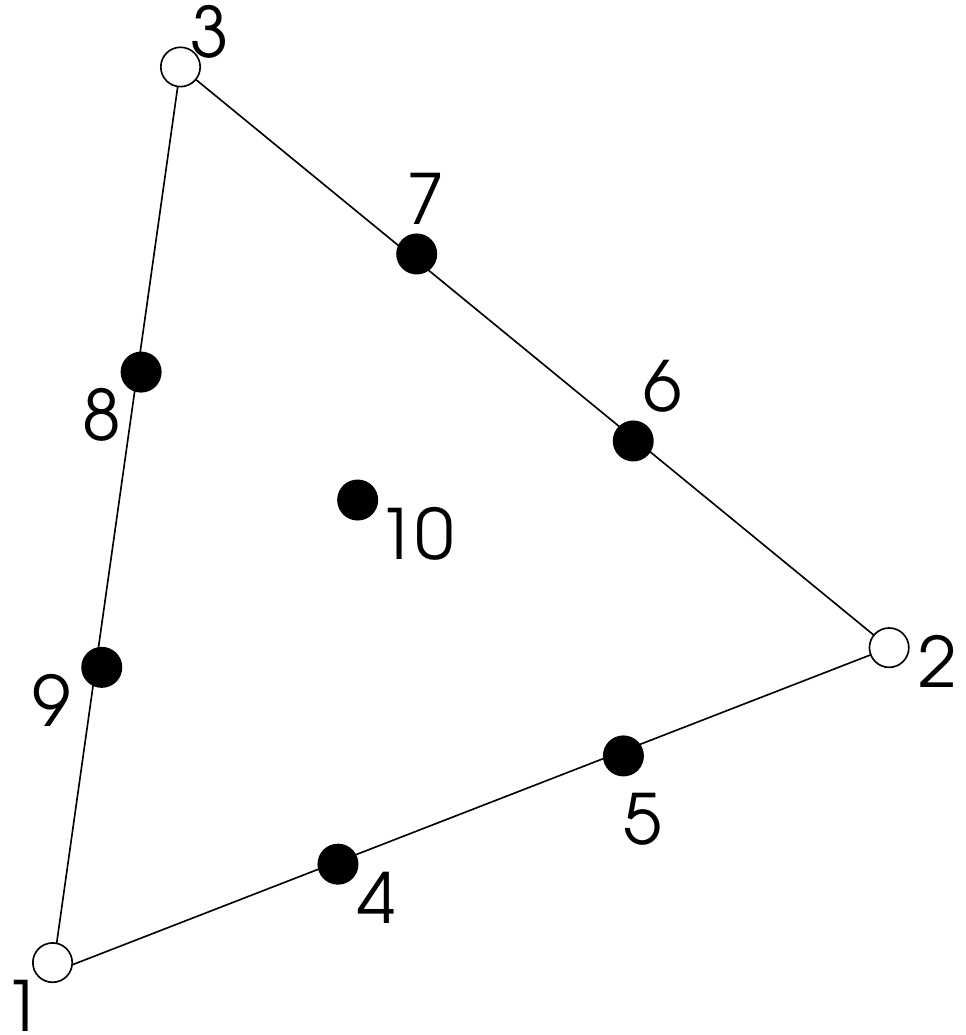}
\caption{\rev{Nomenclature of the DoFs within a $B_2$ (upper triangle) and a $B_3$ (lower triangle) element. }}
\end{figure}
\end{multicols}

\subsection{On the high order differential operator $\mathcal{L}^2$}

The high order differential operator $\mathcal{L}^2_\sigma$ reads
\begin{equation}
\begin{split}
 \mathcal{L}^2_\sigma(U^{(r)})&= \mathcal{L}^2_\sigma(U^{n,1,r}, \ldots , U^{n,M,r})\\&=
\begin{pmatrix}
\sum\limits_{K|\sigma\in K} \Big{(} \int_K \psi_\sigma \big( U_\sigma^{n,M,r} - U_\sigma^{n,0} \big)\,d\mathbf{x} + \displaystyle\int_{t_{n,0}}^{t_{n,M}}  \mathcal{I}_{M} \big( \phi_{\sigma,\mathbf{x}}^{K}(U^{(r)}),s\big) \,ds \Big{)}\\
\vdots\\
\sum\limits_{K|\sigma\in K} \Big{(} \int_K \psi_\sigma \big ( U_\sigma^{n,1,r} - U_\sigma^{n,0} \big)\,d\mathbf{x} + \displaystyle\int_{t_{n,0}}^{t_{n,1}}  \mathcal{I}_{M} \big( \phi_{\sigma,\mathbf{x}}^{K}(U^{(r)}),s\big) \,ds  \Big{)}
\end{pmatrix}.
\end{split}
\label{L2_simplif}
\end{equation}

In practice, we compute the coefficients of the interpolating polynomial $\mathcal{I}_{M}$ of degree $M$ and perform exact integration to obtain the approximation for every row of \eqref{L2_simplif} in the form
\begin{equation}
\label{Interp_M}
\int_{t_{n,0}}^{t_{n,m}}\mathcal{I}_{M}\big(\phi_{\sigma,\mathbf{x}}^{K}(U^{n,0,r}),\ldots,\phi_{\sigma,\mathbf{x}}^{K}(U^{n,M,r}),s\big) \,ds=\sum_{l=0}^{M}\theta_{m,l} \;\phi_{\sigma,\mathbf{x}}^{K}(U^{n,l,r}).
\end{equation}
This high order differential operator ensures a high order approximation of the space-time term $\partial_t U + \text{div}\,\bF(U)$, but is implicit in time, and therefore the iterative formulation \eqref{HO_timestepping} is used in the timestepping to obtain an explicit scheme which ensures high order of accuracy both in space and time.

\subsection{On the choice of the sub-time and correction steps}

\rev{
As outlined in \cite{Abgrall2017}, after $R$ corrections we have $\mathcal{L}^1(U^{(R+1)})=\mathcal{O}(h^{R+1})$, as for each correction holds
\begin{equation}
\mathcal{L}^1(U^{(r+1)})=\mathcal{L}^1(U^{(r)})-\mathcal{L}^2(U^{(r)})=\mathcal{O}(h^{r+1}).
\label{HO_timestepping_order}
\end{equation}
The approximation $ \mathcal{L}^1(U)=0$ corresponds to a two level scheme for each of the sub-time steps $m$, and, thus, the solution $U_h^{n,M,R}$ is obtained from a two-level scheme that is perturbed by an $\mathcal{O}(h^{r+1})$ term. From a result in \cite{Morton}, we see that, given a norm, the stability condition of the method corresponds to $ \mathcal{L}^1$. Further, $ \mathcal{L}^1$ is recast in terms of a forward Euler method and to obtain a method of order $M$ in space, the time step must be divided by $M$ with respect to the time step needed for the first order in space scheme.
Therefore, as a general rule, the idea is to take as many sub-timesteps $M$ as corrections $R$, in order to provide the desired order.
}
\subsection{Extension to systems}\label{Section_Extend_sys}

Out of the schemes described in the previous section, we have considered an approximation given by the limited RDS with an additional jump stabilization in the form of \eqref{schema_RDS_jump}, since the results obtained in \cite{Abgrall2017} have shown the supremacy of the jump formulation with respect to the SUPG scheme in terms of dispersive errors.

The boundary term $\int_{\partial K}\varphi_\sigma \bF(U_h)\cdot \bn$ is evaluated with the same quadrature formula as the face term
$ \int_e [\nabla U_h]\cdot [\nabla \varphi_\sigma]$.
The volume term is evaluated by quadrature as well, and the accuracy requirements on these quadrature formula are similar to those of the discontinuous Galerkin (DG) methods.

The approximation of $\bF$ which we denote by $\bF(U_h)$ can be done in two possible ways. Either from the data $U_h$ one evaluates the values of the flux at the DoFs $\bF_\sigma$ and defines $\bF(U_h)$ as:
\begin{equation}
\label{flux:approximation}
\bF(U_h)\approx \sum_{\sigma\in K}\bF_\sigma \varphi_\sigma,
\end{equation}
which leads to a quadrature-free implementation since the integrals of the shape functions and/or gradients can be evaluated explicitly. Alternatively, one can define $\bF(U_h)$ as the flux evaluated for the local value of $U_h$ at the quadrature point, since both approaches are formally equivalent from the accuracy point of view. The  $\varphi_\sigma$ have degree $k$.
 
Let us now explain how the nonlinear residual \eqref{schema_RDS_jump} is calculated in case of systems, omitting the jump stabilization term for simplicity. We also omit the correction index $r$ and describe the calculation for the $m$-th substep in time.
 
We start by introducing a local Lax-Friedrichs type nodal residual on the steady part of \eqref{system}:
\begin{equation}
\phi_{\sigma,\mathbf{x}}^{K,LxF}(U_h)=\int_{\partial K}\varphi_\sigma \mathbf{F}(U_h)\cdot \mathbf{n}\,d\gamma -\int_K \nabla \varphi_\sigma\cdot \mathbf{F}(U_h)\,d \mathbf{x} +\alpha_K(U_\sigma-\overline{U}_h^K)
\label{phi_LxF_xI}
\end{equation} 
and define the nodal residual in space and time of \eqref{system} as
\begin{equation}
\phi_\sigma^{K,LxF}(U_h) = \int_K \psi_\sigma \big( U_\sigma^{n,m} - U_\sigma^{n,0} \big)\,d\mathbf{x} +   \displaystyle\int_{t_{n,0}}^{t_{n,\rev{m}}}  \mathcal{I}_{M} \big( \phi_{\sigma,\mathbf{x}}^{K,LxF}(U),s\big) \,ds 
\label{phi_LxF}
\end{equation}     
where $\overline{U}_h^K$ is the arithmetic average of all degrees of freedom defining $U_h$ in $K$. \rev{The coefficient $\alpha_K$ is defined via the spectral radius of the flux Jacobian matrix $\bA(U) = \nabla_U\bF(U) \cdot \bn$ as follows 
$$\alpha_K = \max\limits_{\sigma\in K} \bigg ( \rho_S\Big ( \nabla_U\big{(}\bA(U_\sigma) \Big )\bigg ).$$}
\rev{The use of this classic formulation of the Lax-Friedrichs results, nevertheless, in a very dissipative scheme when dealing with higher than second order schemes (see cf. Fig. \ref{Fig:SO_200comp}).
This observation has led to a reformulation of the term $\alpha_K(U_\sigma-\overline{U}_h^K)$ of equation \eqref{phi_LxF_xI}.
This reformulation requires the fullfillment of a Lax-Wendroff like theorem \cite{abg2001d,AbgrallViville2017} that sets the constraint on the conservation of \eqref{RD_nodaldecomp_contrib} at an element interface level and not globally on the element. }

 {\color{blue} To achieve this, we consider the flux approximation \eqref{flux:approximation}. If the basis functions are of degree $k$, the approximation of \eqref{flux:approximation} is denoted by $\bF^{(k)}$.
Following the idea of \cite{AbgrallViville2017}, we can rewrite  the  residuals $\int_K \text{ div }\bF^{(k)} d\mathbf{x}$ as 
\begin{equation}
\int_K \text{div }\bF^{(k)} d\mathbf{x}=\sum_{K_i\subset K} \omega_{K_i} \int_{K_i}\text{ div }\bF^{(1)}\;d\mathbf{x},
\label{divkinto1}
\end{equation}
which is a weighted sum of the first order residuals $\int_{K_i} \text{ div }\bF_{K_i}^{(1)} \; d\bx$, where $\bF_{K_i}^{(1)}$ is the piecewise linear interpolation of the flux $\bF$ and $\bF_\sigma$ represents the values at the vertices $\sigma$ of $K_i$. The weights $\omega_{K_i}$ are positive (refer to \ref{appendix A} for more details).\\
Equation \eqref{divkinto1} allows to reformulate a new version of  \eqref{phi_LxF_xI} as
\begin{equation}
\sum_{K_i\in K,\sigma \in K_i}  \omega_{K_i} \,\big{[} \int_{\partial K_i}\varphi_\sigma \mathbf{F}^{(1)}\cdot \mathbf{n}\,d\gamma_i -\int_{K_i} \nabla \varphi_\sigma\cdot \mathbf{F}^{(1)}\,d \mathbf{x} + \alpha_{K_i}(U_{\sigma,K_i}-\overline{U}_{h,K_i})\big{]},
\label{phi_LxF_sub}
\end{equation}
which corresponds to recast the Lax-Friedrichs term as a weighted sum over each node $i$ belonging to a sub-cell $K_i$ within a cell $K$. In the proposed formulation, we do not write the Galerkin term for degree $k$
$$\int_K \varphi_\sigma \text{ div }\bF^{(k)}\; d\mathbf{x}=-\int_K\nabla\varphi_\sigma \cdot \bF^{(k)}\; d\mathbf x +\int_{\partial K} \varphi_\sigma \bF^{(k)}\cdot\bn \; d\gamma,$$
as a weighted sum of the Galerkin term for degree $1$ (though this is also possible, with positive weights), but we consider the Lax Friedrichs scheme for the sub-elements $K_i$ and we weight them in such a way that the conservation at the element-level $K$ is recovered.

To get a non oscillatory scheme, the dissipation terms $\alpha_{K_i}$ are defined by
$$\alpha_{K_i} = \max\limits_{\sigma\in K_i} \bigg ( \rho_S\Big (\bA(U_\sigma) \Big )\bigg ).$$
The weights $\omega_{K_i}$ are set to $1$ in the one dimensional case, so that the sum corresponds basically to a telescopic sum over the sub-cells. In two-dimensions, we have set for the Bernstein approximation of order 2
\begin{equation}
\omega_{K_i}=
\begin{cases}
\dfrac{2}{3},\quad \text{for}\,\, i=1,2,3\\ 
2,\quad \text{for}\,\, i=4\\ 
\end{cases}
\end{equation}
and for the Bernstein approximation of order 3 
\begin{equation}
\omega_{K_i}=
\begin{cases}
\dfrac{1}{2},\quad \text{for}\,\, i=1,2,..,6\\
1,\quad \text{for}\,\, i=7,8,9\\
\end{cases}
\end{equation}
following the sub-cell nomenclature as in Fig. \ref{Fig:sub_elements} (see for more details \cite{AbgrallViville2017}).
\begin{figure}[H]
\begin{center}
\includegraphics[scale=0.4]{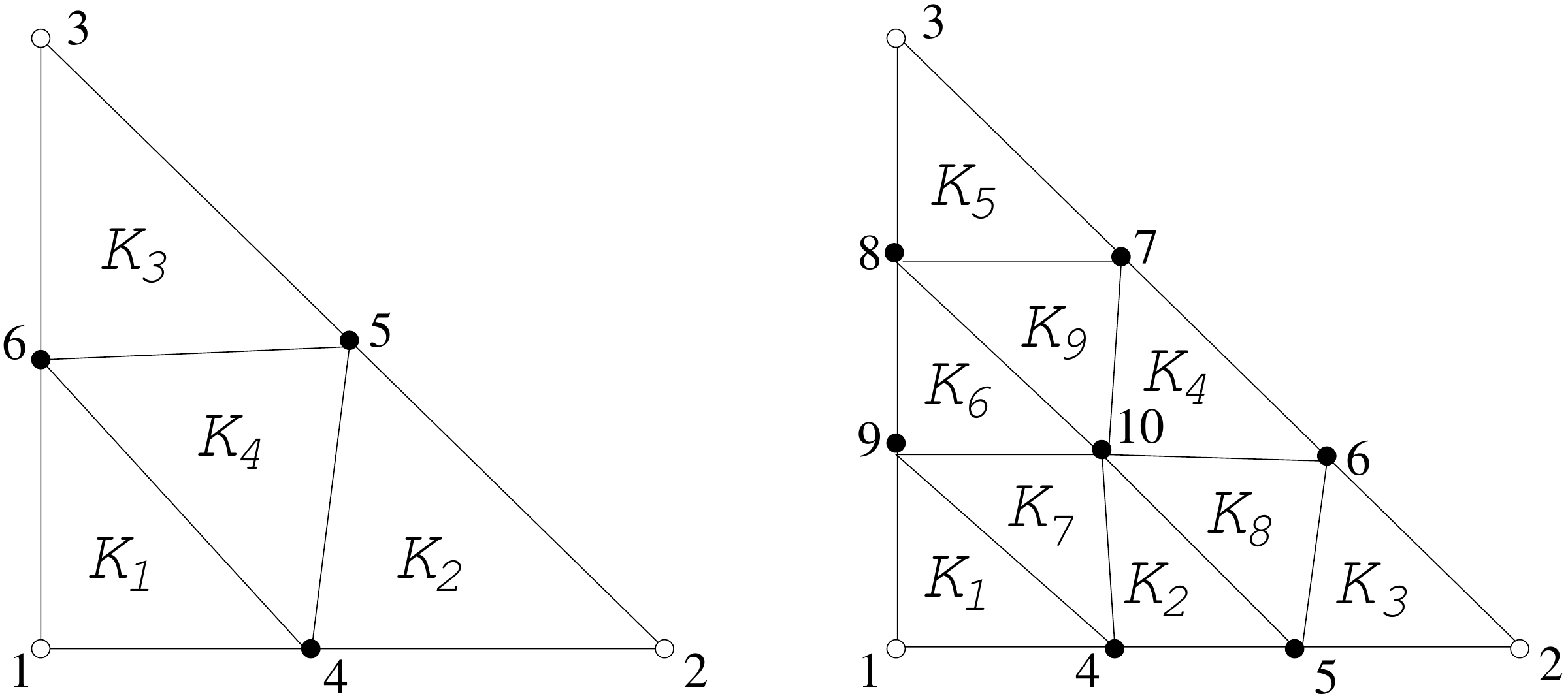}
\caption{Nomenclature of the sub-cells $K_i$ within a generic linear triangular element $K$ for 'B2' (left) and 'B3' (right). }
\label{Fig:sub_elements}
\end{center}
\end{figure}
}

In order to achieve high order accuracy and guarantee the monotonicity of the solution near strong discontinuities, we proceed as follows. For a scalar problem we would compute the distribution coefficients $\beta_\sigma^K$ as
\begin{equation}
\beta_\sigma^K(U_h)=\dfrac{\max\big( \frac{\phi_\sigma^{K,LxF}}{\phi^K},0 \big)}{\sum\limits_{j\in K} \max\big( \frac{\phi_j^{K,LxF}}{\phi^K},0 \big)}, \quad \phi^K = \sum_{\sigma\in K} \phi_\sigma^{K,LxF}
\label{beta_singularity}
\end{equation}
\rev{In case of systems, to allow less dissipation, \eqref{beta_singularity} is applied to each variable by considering their characteristic decomposition as described e.~g. in \cite{DeSantis2015}}.
To this end, one considers the eigen-decomposition of the Jacobian matrix $\bA(U) = \nabla_U\bF(U) \cdot \bn$ of the flux $\mathbf{F}$ with respect to the state $\overline{U}_h^K$, where as $\bn$ we take the average fluid velocity vector or we choose an arbitrary direction (for example the $x$-coordinate) in case the average velocity vanishes. The matrix composed of the right eigenvectors of $\bA(U)$ is denoted by $\bR$, so that $\bL = \bR^{-1}$ is the matrix of left eigenvectors.

More precisely, and as described in \cite{DeSantis2015}, the distribution coefficients for the system of equations are calculated in local characteristic variables by projecting the first order residuals onto a space of left eigenvalues, as
\begin{equation}
\label{phi_char}
\hat{\phi}_\sigma^{K,LxF} = \bL\,\phi_\sigma^{K,LxF}, \quad \hat{\phi}^K = \bL\,\phi^K.
\end{equation}

The high order nodal limited residuals are then obtained as follows. We first calculate the distribution coefficients according to
\begin{equation}
\beta_\sigma^K = \dfrac{\max\big( \frac{\hat{\phi}_\sigma^{K,LxF}}{\hat{\phi}^K},0 \big)}{\sum\limits_{j\in K} \max\big( \frac{\hat{\phi}_j^{K,LxF}}{\hat{\phi}^K},0 \big)}, \quad \hat{\phi}^K = \sum_{\sigma\in K} \hat{\phi}_\sigma^{K,LxF}.
\label{beta_char}
\end{equation}
Next, we apply the blending scheme
\begin{equation}
\hat{\phi}_\sigma^{K,\star} = (1-\Theta)\,\beta_\sigma^K \hat{\phi}^K + \Theta\,\hat{\phi}_\sigma^{K,LxF},
\label{limiting}
\end{equation} 
where the blending coefficient $\Theta$ is defined by
\begin{equation}
\label{theta}
\Theta = \dfrac{\big| \hat{\phi}^K \big|}{\sum\limits_{\sigma'\in K} \big| \hat{\phi}_{\sigma'}^{K,LxF} \big|}.
\end{equation}
Clearly, $0 \leq \Theta \leq 1$, and $\Theta = O(h)$ for a smooth solution, thus ensuring accuracy and $\Theta = O(1)$ at the discontinuity, thus ensuring monotonicity \cite{abg}.
Finally, the high order nodal residuals are projected back to the physical space:
\begin{equation}
\label{phi_phys}
\phi_\sigma^{K,\star} = \bR\,\hat{\phi}_\sigma^{K,\star}.
\end{equation}
This guarantees that the scheme is high order in time and space and (formally) non-oscillatory, see \cite{mario,abg} for more details. 

\begin{remark} \rev{ When possible singularities may arise as, for example, due to pressures close to zero, instead of applying the characteristic limiting \eqref{phi_char}-\eqref{phi_phys}, we adopt locally in the affected cells the limiting as in \eqref{beta_singularity}.These two limiting strategies are both Lipschitz continuous, so that the switch does indeed not affect the property \eqref{error}, causing convergence problems.
Further, the only situation when the pressure becomes close to zero, or, eventually negative, typically occurs in problems with strong interacting discontinuities, and across shocks one would have a first order monotone method due to the limiting formulation.}
\end{remark}

\rev{
After the application of the limiter, we add the jump stabilization term 
\begin{equation}
\phi_{\sigma,\mathbf{x}}^{K,jump}(U_h)= \sum_{\text{edges of }K}
    \theta_1 h_e^2 \int_e [\nabla U_h]\cdot [\nabla \varphi_\sigma]\,d\gamma+\sum_{\text{edges of }K}
    \theta_2 h_e^4 \int_e [\nabla^2 U_h\bn]\cdot [\nabla^2 \varphi_\sigma\bn]\,d\gamma
    \label{phi_burman}
\end{equation} 
where $\bn$ is a normal to $e$.
}
\rev{
In general, since we are adding the edge stabilization terms of \eqref{phi_burman} to the residual distribution scheme after the high order limiting, the question may arise, wheather the inclusion of an unlimited high-order stabilization term might destroy monotonicity-preserving properties \cite{Jameson1995}. Numerical experiments show that this method is essentially non-oscillatory. Formally the monotonicity property is violated, in practice, nevertheless, this is not the case: we can see very small undershoots/overshoots. Moreover, in our experiments, we have observed that an appropriate choice of the coefficient $\theta_1$ and $\theta_2$ in \eqref{phi_burman} does not lead to any spurious oscillations at shocks and has only the beneficial effect to stabilize the solution for high order, ensuring, thus the aimed accuracy.
}

\section{Numerical results}
\label{sec:results}

To assess the accuracy and robustness of the proposed high order residual distribution scheme,  in the following section we perform the convergence analysis for the wave equation and isentropic flow and study several benchmark problems in one and two spatial dimensions for the Euler equations of gas dynamics. In the following, we shall refer to the second order scheme obtained by using linear shape functions on each element as 'B1'. Higher order approximations are obtained by choosing quadratic ('B2') or cubic ('B3') Bernstein polynomials as shape functions.
\rev{ For the B1 approximation we consider $M=2$ and $R=2$, for B2 $M=3$ and $R=3$ and, finally, for B3 $M=4$ and $R=4$ in algorithm \eqref{HO_timestepping}.
All test cases are advanced in time using the Courant-Friedrichs-Lewy condition $\Delta t= \text{CFL} \cdot\Delta x$ which is then updated by computing $\Delta t= \text{CFL}\cdot \min_{\sigma} \big{(} \frac{\Delta x_{\sigma}}{|u_{\sigma}+c_{\sigma}|}\big{)} $, where $\Delta x_{\sigma}$ represents the volume of the cell corresponding to the considered degree of freedom $\sigma$ and $|u_{\sigma}+c_{\sigma}|$ the spectral radius of the solution in $\sigma$.
We have set for all the considered tests a fixed $CFL=0.1$.
The parameters of \eqref{phi_burman} $\theta_1$ and $\theta_2$ depend on the order of accuracy and on the typology of considered system, i.e. they change from 1D to 2D and from the wave equation to the Euler system. In the following considered benchmark problems, we set empirically the values of $\theta_1$ and $\theta_2$ that show a robust stabilization capability. 
}

\subsection{Numerical results for 1D test cases}

\subsubsection{Convergence study: wave equation}

We start by considering the one dimensional wave equation $q_{tt} - a^2 q_{xx} = 0$
which we rewrite as a first-order system of PDEs with respect to the variables $u = q_t$ and $v = q_x$:
\begin{alignat*}{3}
& u_t - a^2 &&v_x &&= 0, \\
& v_t - &&u_x &&= 0.
\end{alignat*}

We perform the convergence analysis on the smooth problem with initial condition
\[ q(x,0) = \exp\big(-\beta (x-1/2)^2\big)\sin(\alpha x), \quad -1 \leq x \leq 2, \]
where we set $\alpha = \beta = 100$ and $a=1$. The initial condition for the new variables $u$ and $v$ is derived from $q(x,0)$ as
\begin{align*}
& u(x,0) = 0, \\
& v(x,0) = \exp\big(-\beta (x-1/2)^2\big)\big(\alpha\cos(\alpha x) - 2\beta(x-1/2)\sin(\alpha x)\big).
\end{align*}
The final time of the computation is $T=0.5$.

The computations have been performed for B1, B2 and B3 shape functions, leading to schemes of second, third and fourth order of accuracy, respectively. 
\rev{In order to stabilize the approximation, we set w.r.t. \eqref{phi_burman} for B1 $\theta_1=0.2$ and $\theta_2=0$, for B2 $\theta_1=0.1$ and $\theta_2=0$ and for B3 $\theta_1=2.$ and $\theta_2=4.$}

\begin{figure}[h!]
\begin{center}
d{\includegraphics[width=0.45\textwidth]{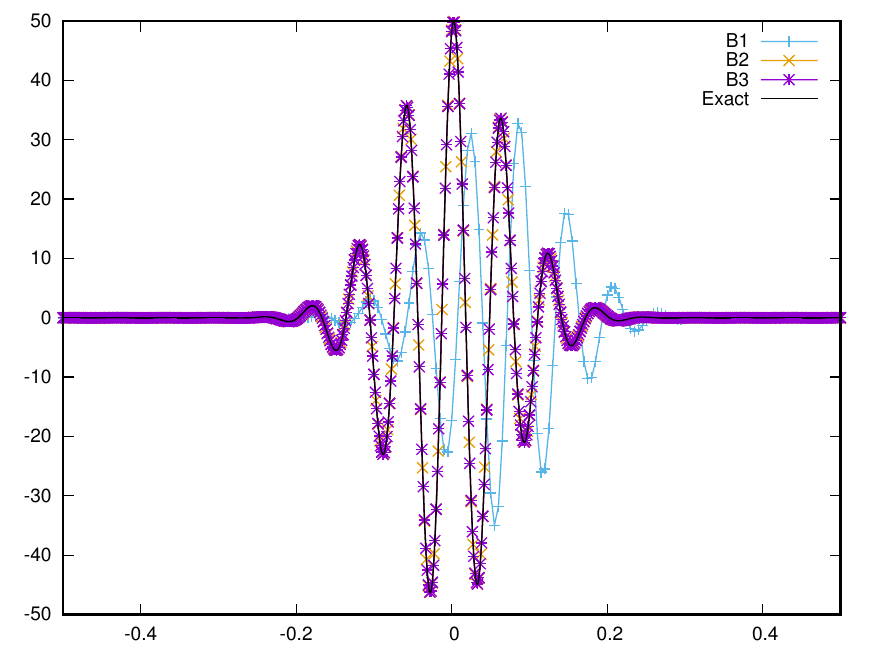}}
\end{center}
\caption{Numerical solution for the 1D wave system using $600$ cells at $T=0.5$ in the range $-0.5 \leq x \leq 0.5$}
\label{wave1D:sol}
\end{figure}

\begin{figure}[H]
\begin{center}
\subfigure[u]{\includegraphics[width=0.47\textwidth]{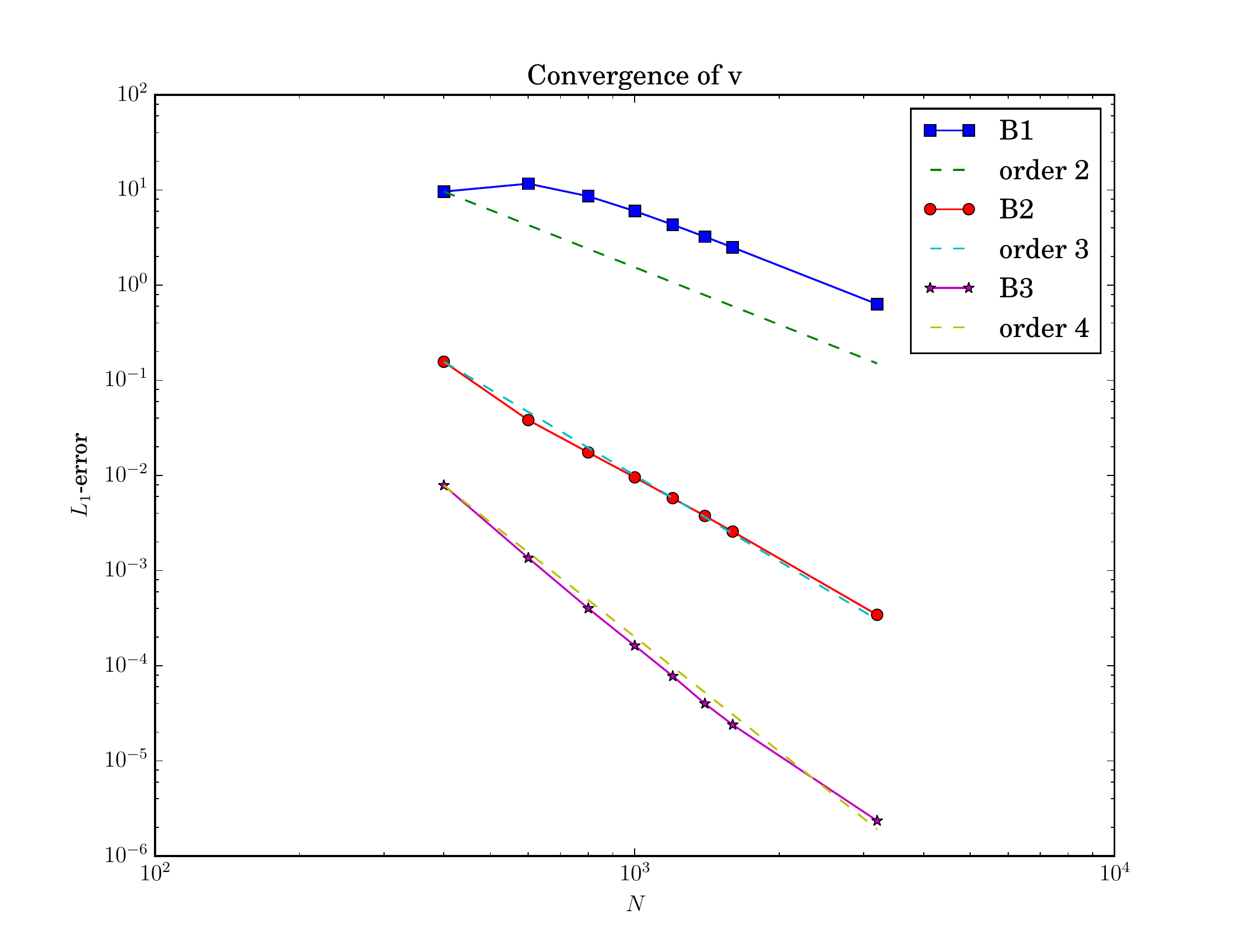}}
\subfigure[v]{\includegraphics[width=0.47\textwidth]{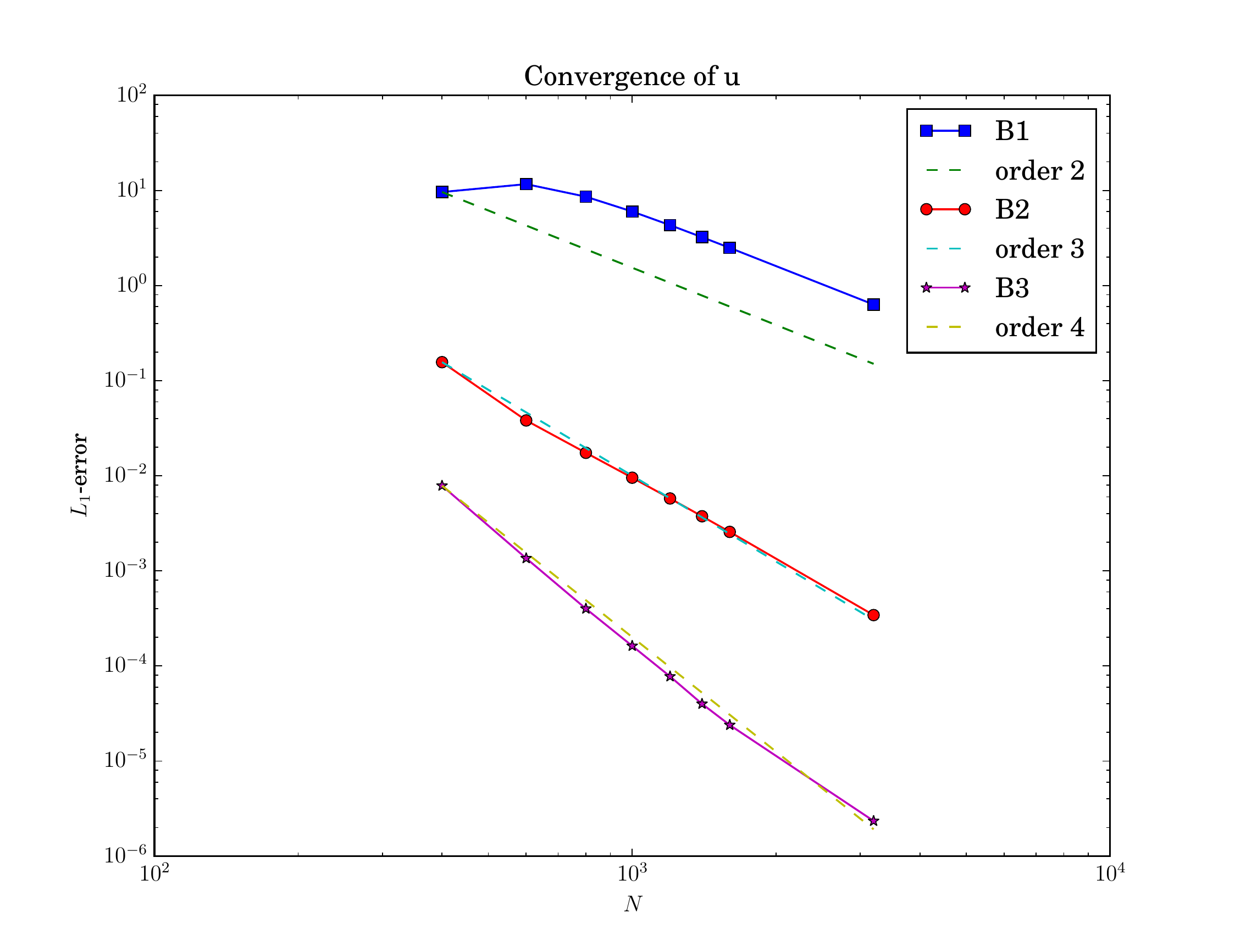}}
\caption{\rev{Convergence plot for the wave equation in 1D at $T=0.5$. B3 counts $R=8$ corrections, w.r.t. Table \ref{Table:wave1D_convergence}}}
\label{Fig:wave1D_convergence}
\end{center}
\end{figure}

\rev{The numerical solution obtained with $600$ grid cells is shown in Fig.~\ref{wave1D:sol} and the convergence plot is given in Fig.~\ref{Fig:wave1D_convergence} along its Table \ref{Table:wave1D_convergence}. 
In Table \ref{Table:wave1D_convergence} we have reported two different convergence studies for B3 elements: one with $R=4$ corrections which gives a decreasing order of accuracy along a mesh refinement, and one with the double amount of corrections, which results in the predicted convergence rates of fourth order.}
We see that the solution given by B1 elements fails to capture the correct location of the waves in discretized domains with low number of cells, while the B2 and B3 elements at the same mesh are already able to provide a very accurate solution, however, the situation improves for B1 elements as the mesh is refined, which can be seen from the convergence plot. The scheme reaches the theoretically predicted convergence rates for all approximation orders that we have tested here.

\begin{table}
\begin{center}
\rev{{\small
\begin{tabular}{c|| l l |l l| l l| l l }
 & B1       & &B2 & & B3&  & B3 \\
 &$M=1$ &  &  $M=2$                          &   &   $M=4$,                          &   &           $M=4$,                       \\
 &$R=2$ &  &   $ R=3$                        &    &    $R=4$                         &   &    $R=8$                             \\
 & & &  & &  & &  &\\
$\log_{10}(h)$ & $L_1$-error  & slope & $L_1$-error & slope & $L_1$-error & slope & $L_1$-error & slope\\\\\hline\\
$2.6021$ & $9.623$ & $-$ & $1.565 \cdot 10^{-1}$ & $-$ & $ 5.105 \cdot 10^{-2}$ & $-$ & $7.857 \cdot 10^{-3}$& $-$\\
$2.7782$ & $11.655$ & $-0.47$ & $3.817 \cdot 10^{-2}$ & $3.48$ & $1.092 \cdot 10^{-2}$ & $3.803$ & $1.356\cdot 10^{-3}$ & $4.33$\\
$2.9031$ & $8.601$ & $1.06$ & $ 1.744 \cdot 10^{-2}$ & $2.72$ & $4.537 \cdot 10^{-3}$ & $3.05$ & $4.004\cdot 10^{-4}$ & $4.24$\\
$3.0000$ & $6.006$ & $1.61$ & $9.543 \cdot 10^{-3}$ & $2.70$ &$2.520\cdot 10^{-3}$ & $2.64$ & $1.625 \cdot 10^{-4}$  & $4.04$ \\
$3.0792$ & $4.319$ & $1.81$ & $ 5.772 \cdot 10^{-3}$ & $2.76$ & $1.600 \cdot 10^{-3}$ & $2.49$ & $7.654 \cdot 10^{-5}$ &$4.05$\\
$3.1461$ & $3.228$  & $1.89$ & $3.752 \cdot 10^{-3}$ & $2.79$ & $1.109\cdot 10^{-3}$ &$2.38$ & $3.996 \cdot 10^{-5}$ & $4.31$\\
$3.2041$ & $2.492$ &$1.94$ & $3.955 \cdot 10^{-4}$ & $1.57$& $8.237 \cdot 10^{-4}$ & $2.23$ & $2.395 \cdot 10^{-5}$ & $3.83$\\
$3.5051$ & $0.632$ & $1.98$ & $1.134 \cdot 10^{-4}$ & $1.80$& $1.873 \cdot 10^{-4}$& $2.14$ & $2.346 \cdot 10^{-6}$ & $ 3.35$
\end{tabular}
\caption{Convergence study for the wave system in 1D at $T=0.5$}
\label{Table:wave1D_convergence}
}}
\end{center}
\end{table}

\subsubsection{Convergence study: smooth isentropic flow}\label{Isoflow_1D}

The next considered test case in 1D \rev{for the Euler system} is performed to assess the accuracy of our scheme on a smooth isentropic flow problem introduced in \cite{ChengShu2014}. The initial data for this test problem is the following:
\begin{equation*}
\rho_0(x) = 1 + 0.9999995\sin(\pi x), \quad u_0(x) = 0, \quad p_0(x) = \rho^{\gamma}(x,0),
\end{equation*}
with $x \in [-1,1]$, $\gamma=3$ and periodic boundary conditions. 

The exact density and velocity in this case can be obtained by the method of characteristics and is explicitly given by
\begin{equation*}
\rho(x,t) = \dfrac12\big( \rho_0(x_1) + \rho_0(x_2)\big), \quad u(x,t) = \sqrt{3}\big(\rho(x,t)-\rho_0(x_1) \big),
\end{equation*}
where for each coordinate $x$ and time $t$ the values $x_1$ and $x_2$ are solutions of the non-linear equations
\begin{align*}
& x + \sqrt{3}\rho_0(x_1) t - x_1 = 0, \\
& x - \sqrt{3}\rho_0(x_2) t - x_2 = 0.
\end{align*}

\rev{The smooth isoentropic flow test has been run with the parameters in \eqref{phi_burman} as follows: B1 $\theta_1=1$ and $\theta_2=0$; B2 $\theta_1=1$ and $\theta_2=0$; B3 $\theta_1=3$ and $\theta_2=10$.}

\begin{figure}[H]
\begin{center}
\subfigure[density]{\includegraphics[width=0.45\textwidth]{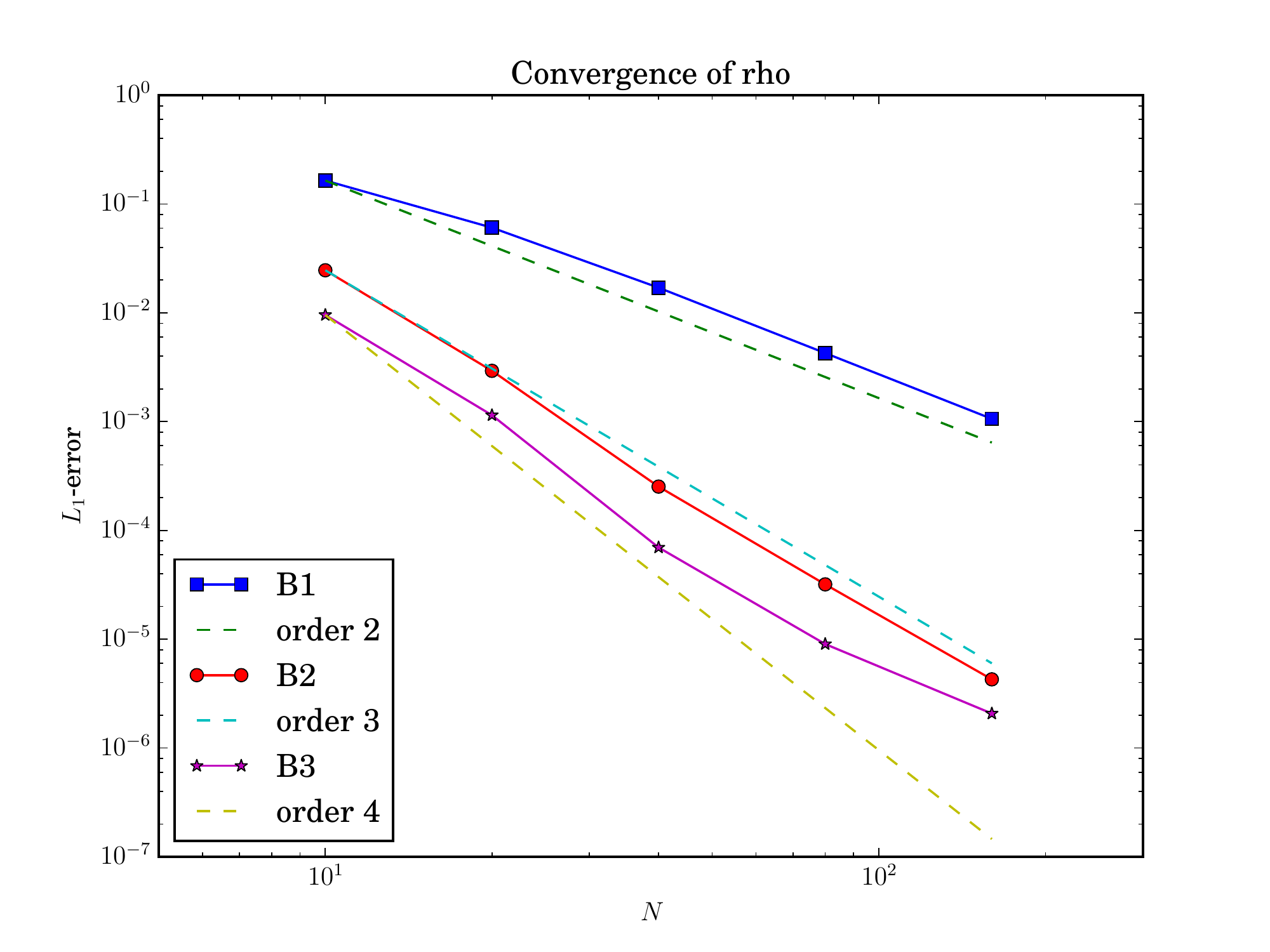}}
\subfigure[velocity]{\includegraphics[width=0.45\textwidth]{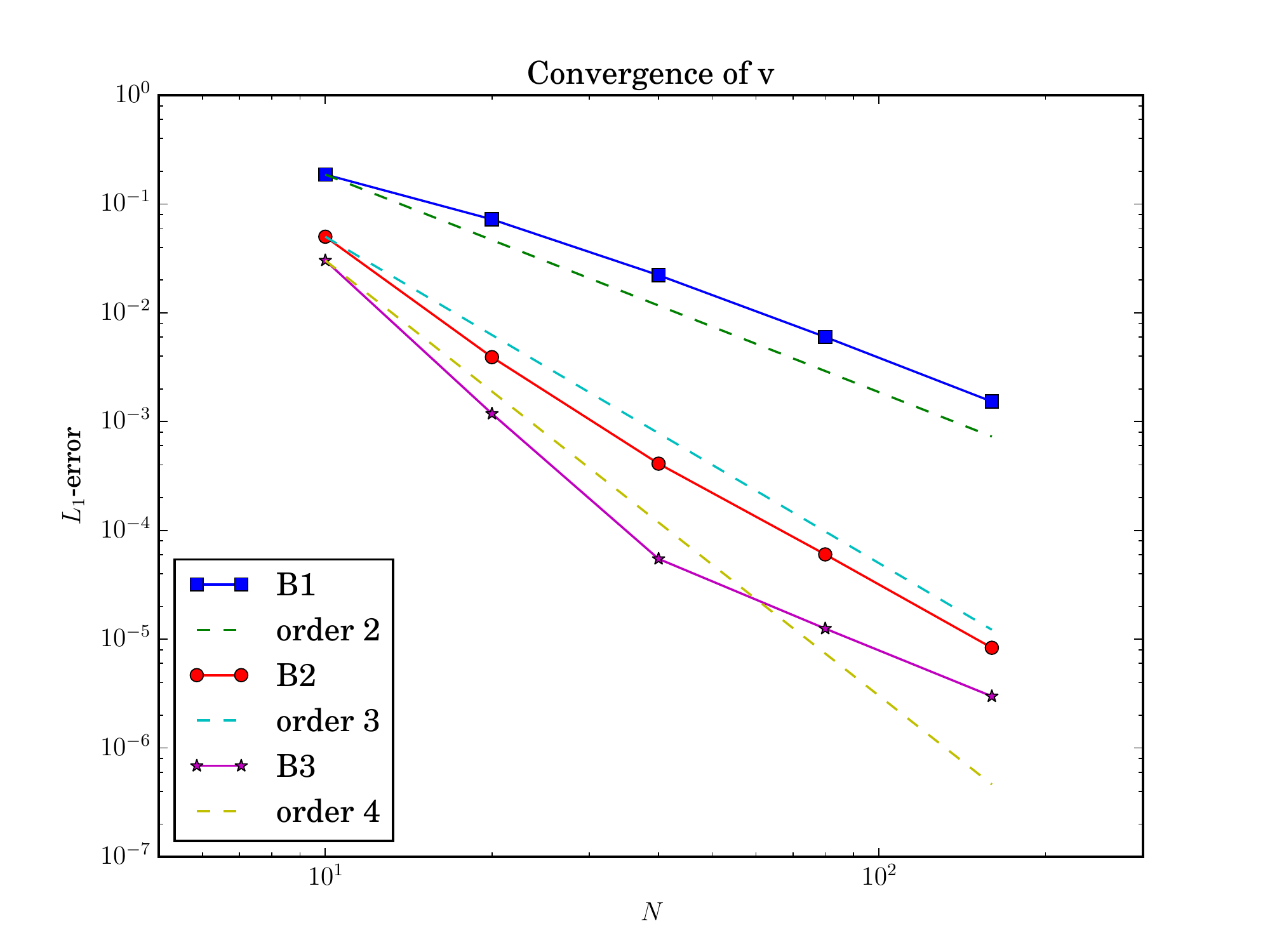}} \\
\subfigure[pressure]{\includegraphics[width=0.45\textwidth]{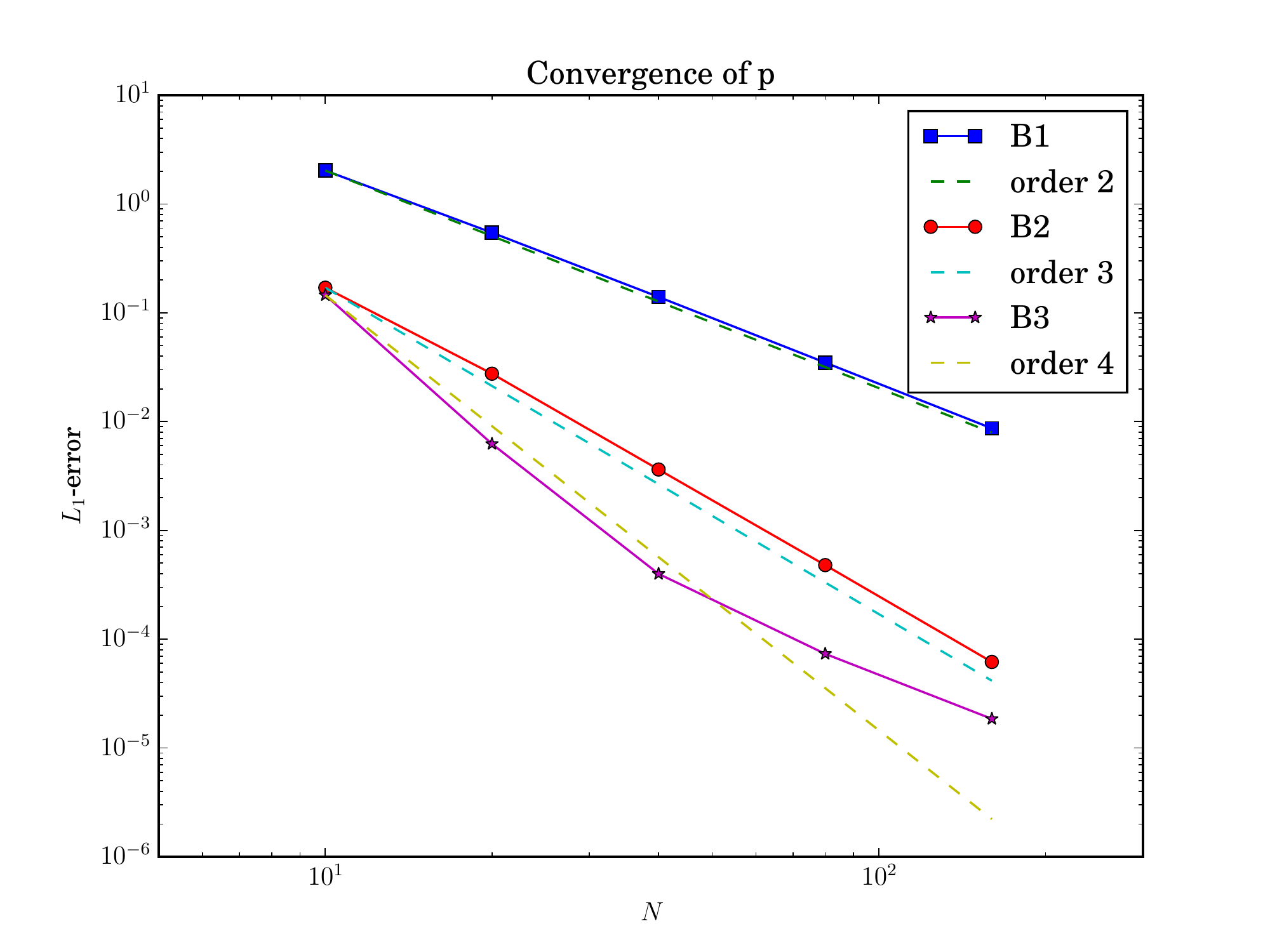}}
\end{center}
\caption{Convergence plot for a smooth isentropic flow in 1D at $T=0.1$.}
\label{Fig:convergence}
\end{figure}

\begin{figure}[H]
\begin{center}
\subfigure[density]{\includegraphics[width=0.45\textwidth]{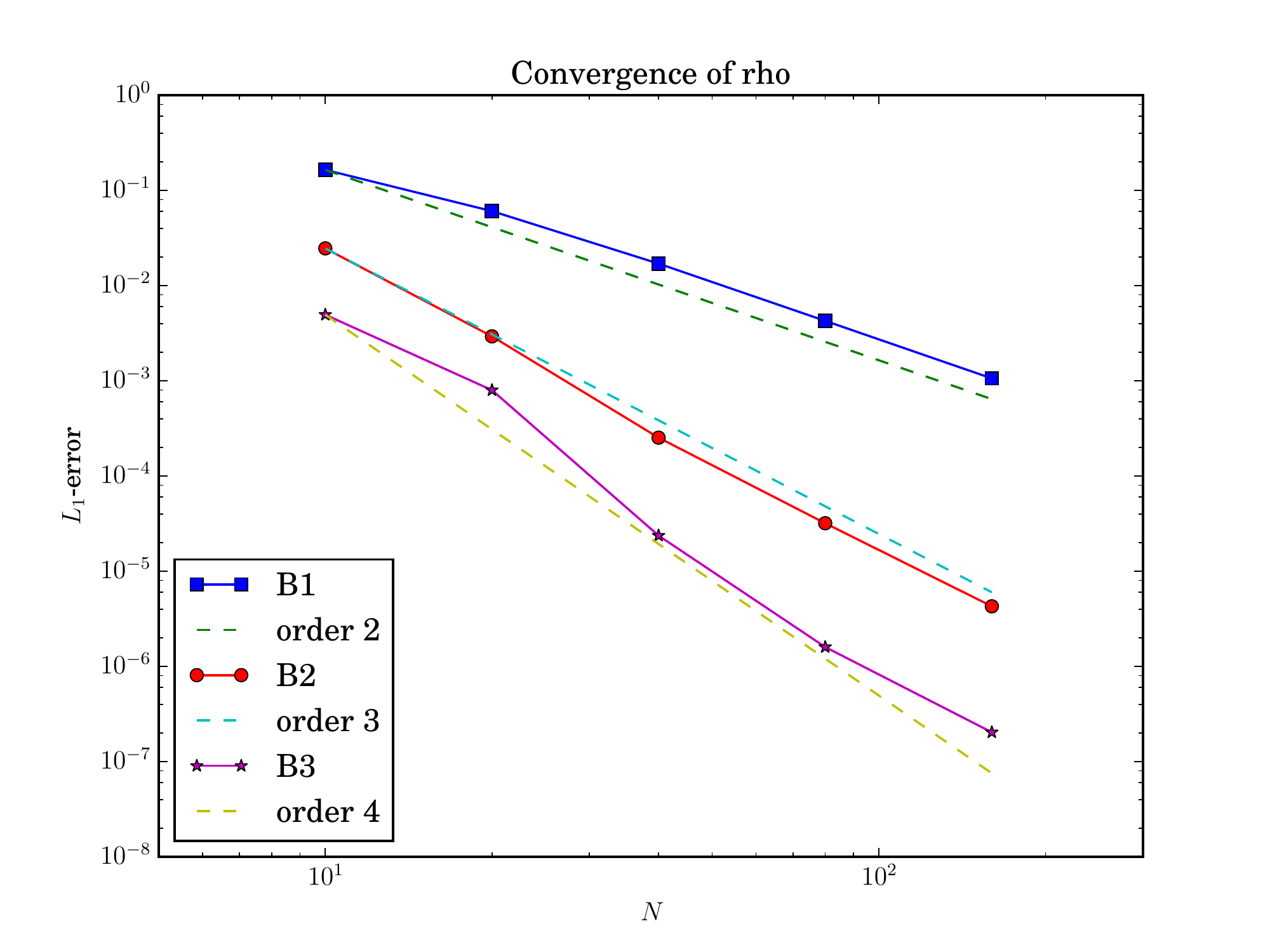}}
\subfigure[velocity]{\includegraphics[width=0.45\textwidth]{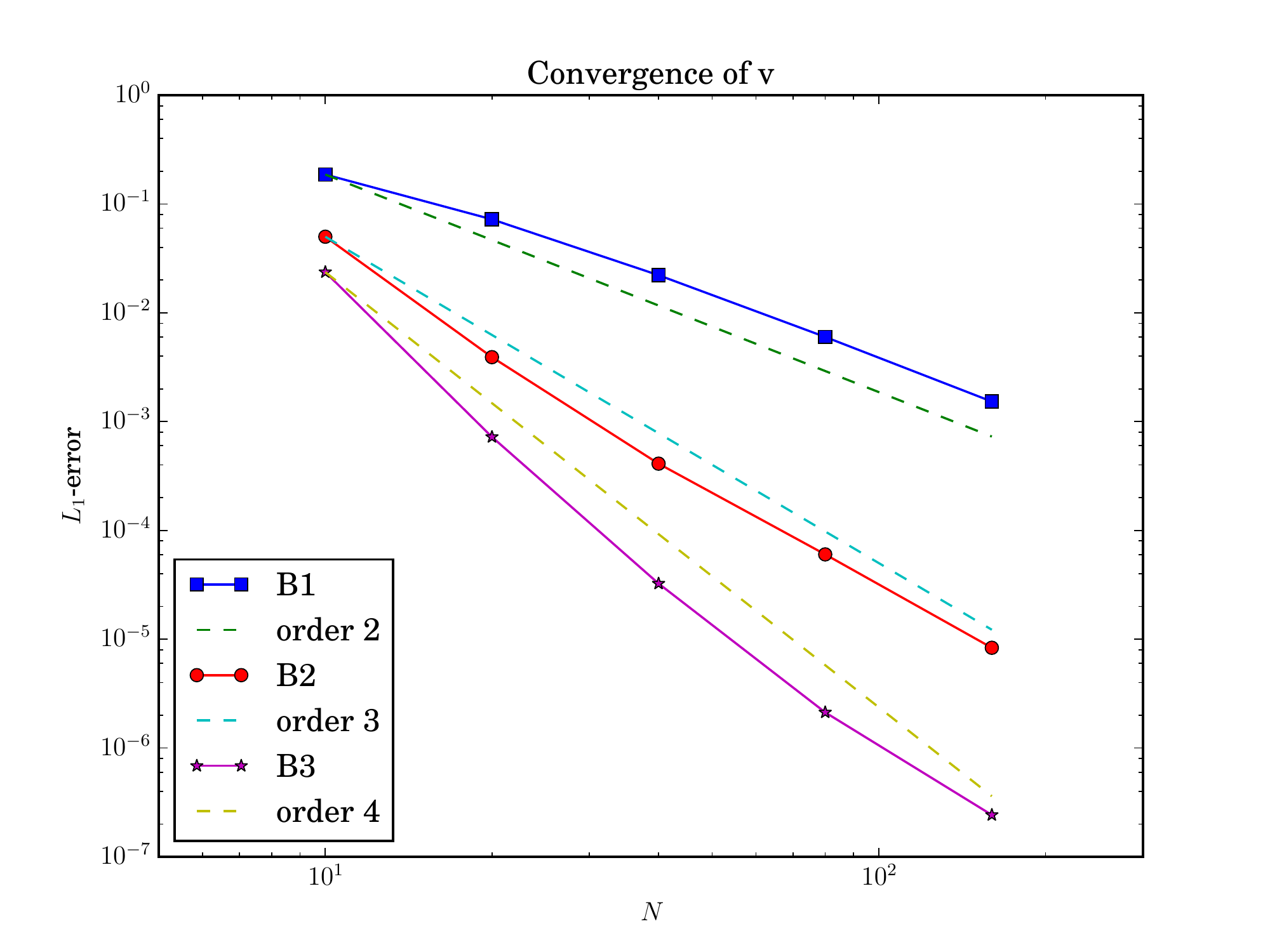}} \\
\subfigure[pressure]{\includegraphics[width=0.45\textwidth]{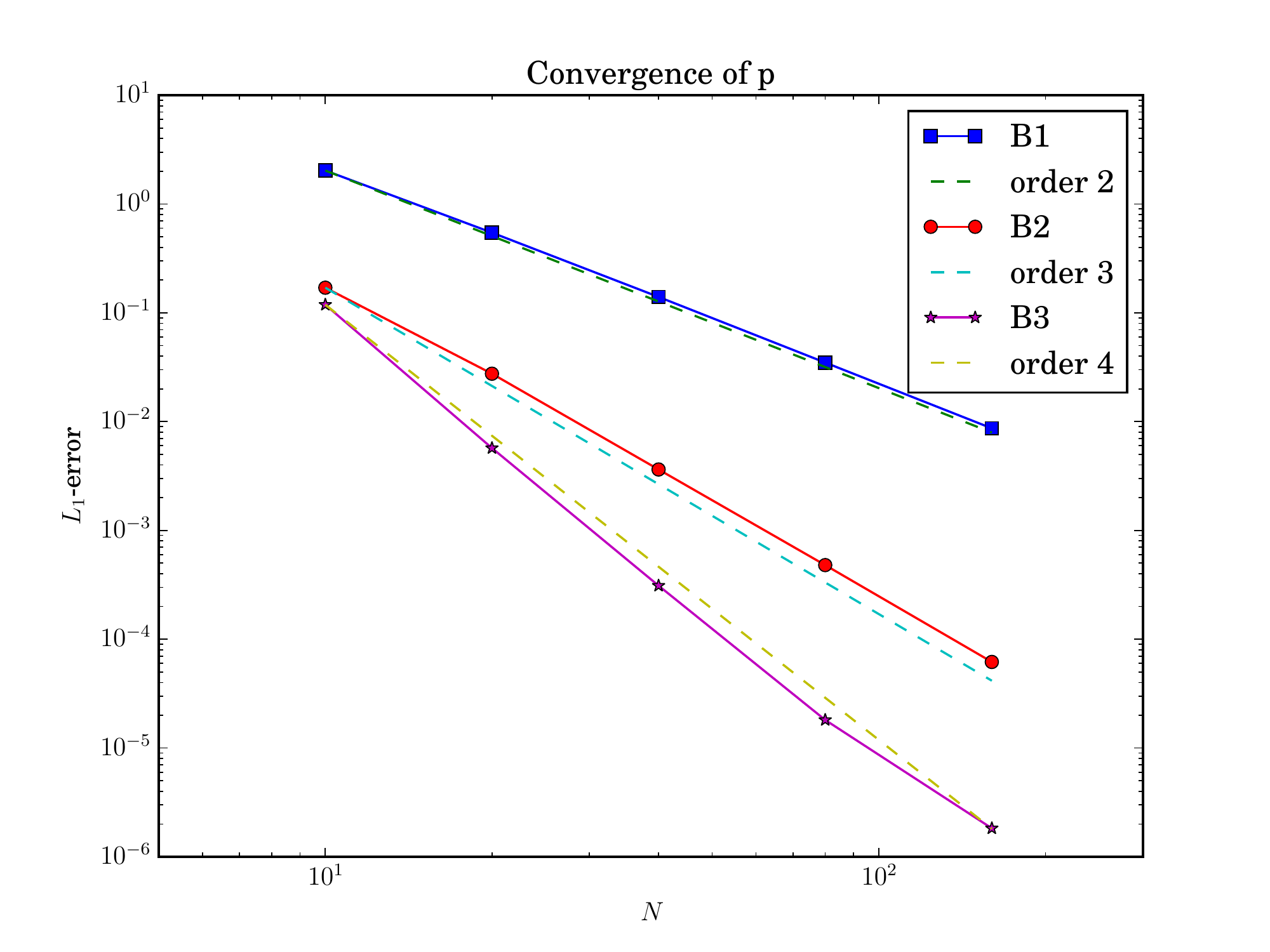}}
\end{center}
\caption{Convergence plot for a smooth isentropic flow in 1D at $T=0.1$ with doubled number of corrections for 'B3'.}
\label{Fig:convergence_10corr}
\end{figure}

\rev{The convergence of the second ('B1'), third ('B2') and fourth ('B3') order RD schemes is demonstrated in Fig.~\ref{Fig:convergence}.
We observe an overall good convergence rate for all the variables. 
It is nevertheless interesting to note, that increasing the amount of the performed corrections, i.e. setting for the fourth order scheme $R=8$ corrections, greatly improves the convergence rate or 'B3' as shown in Fig.~\ref{Fig:convergence_10corr}. }

\subsubsection{Sod's shock tube problem}

\begin{figure}[h!]
\begin{center}
\subfigure[density]{\includegraphics[width=0.45\textwidth]{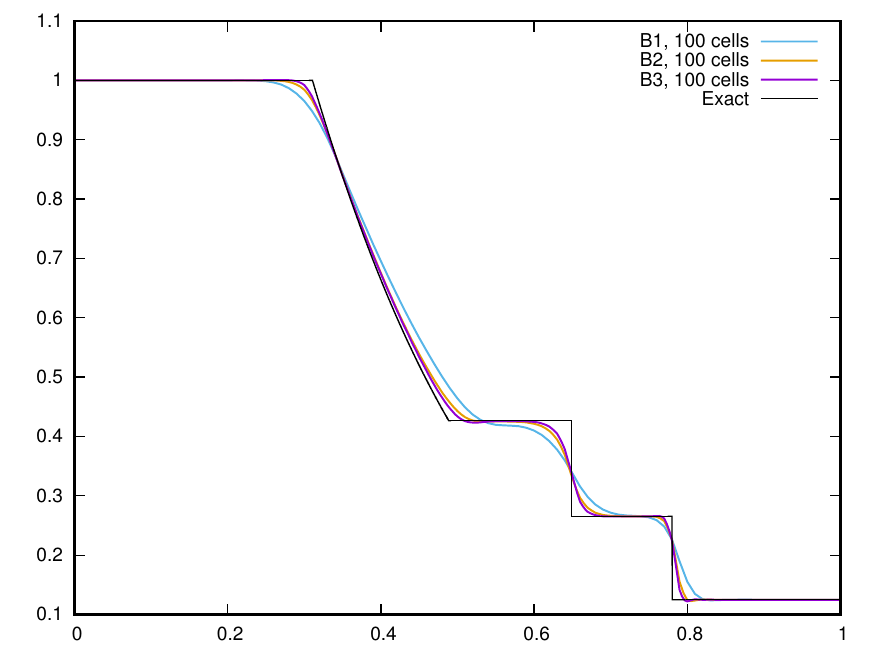}}
\subfigure[velocity]{\includegraphics[width=0.45\textwidth]{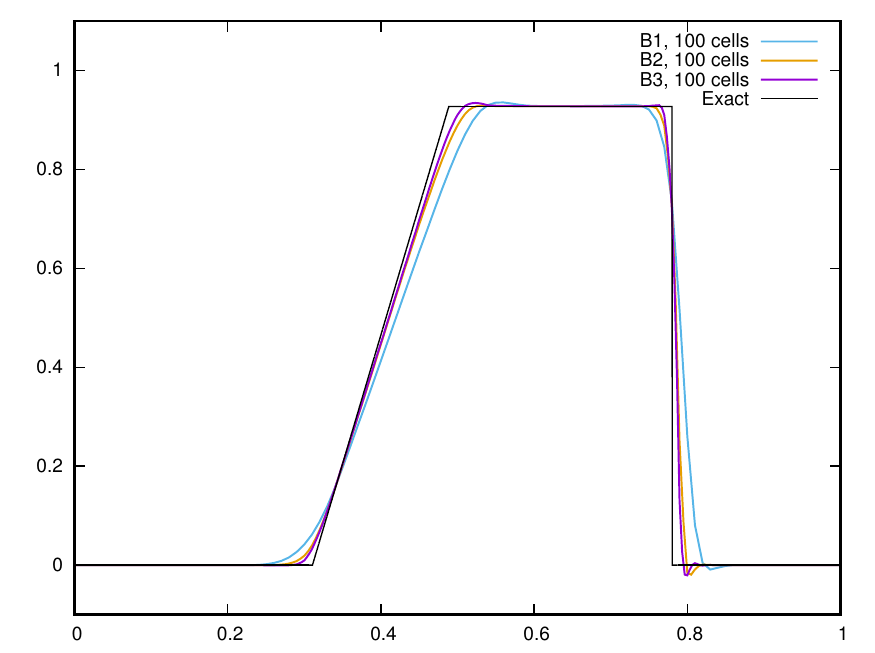}} 
\subfigure[pressure]{\includegraphics[width=0.45\textwidth]{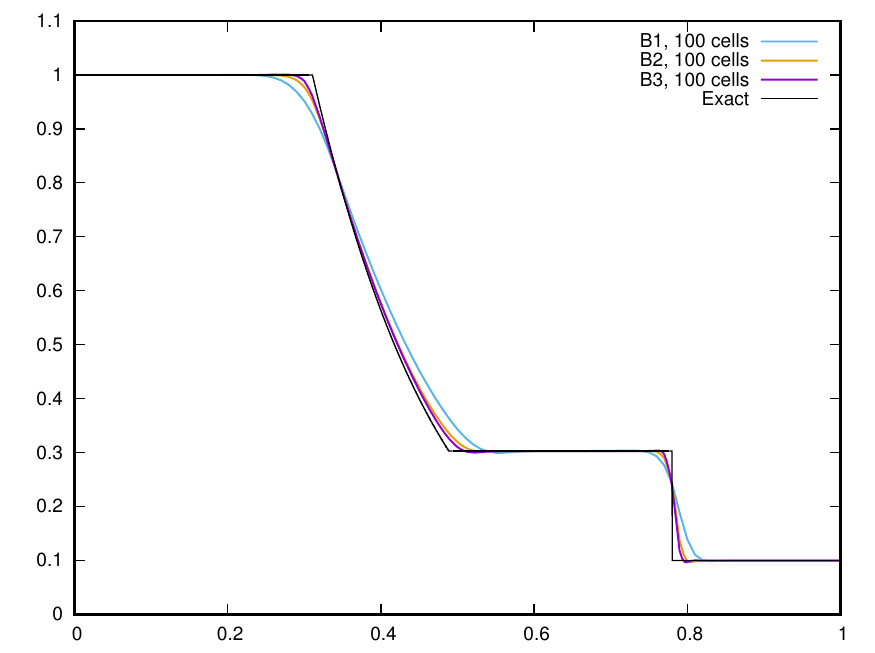}} 
\caption{Comparison between $B1,\;B2,\,B3$ for the Sod 1D test at $T=0.16$}
\label{Fig:Sod} 
\end{center}
\end{figure}

The Sod shock tube is a classical test problem for the assessment of the numerical methods for solving the Euler equations. Its solution consists of a left rarefaction, a contact and a right shock wave. The initial data for this problem is given as follows:
\begin{equation*}
(\rho_0, u_0, p_0) = 
\begin{cases}
 (1, 0, 1), \quad &x < 0, \\
 (0.125, 0, 0.1), \quad &x > 0.
\end{cases}
\end{equation*}
\rev{The jump stabilization \eqref{phi_burman} parameters have been set as in Section ~\ref{Isoflow_1D}.}
The results of the simulations comparing the second, third and fourth order RD scheme are illustrated in Fig.~\ref{Fig:Sod}. The numerical solution converges to the exact one and higher order schemes show much more accurate approximation then the second order scheme.

\subsubsection{Woodward-Colella problem}

\begin{figure}[h!]
\begin{center}
\subfigure[$400$ cells]{\includegraphics[width=0.45\textwidth]{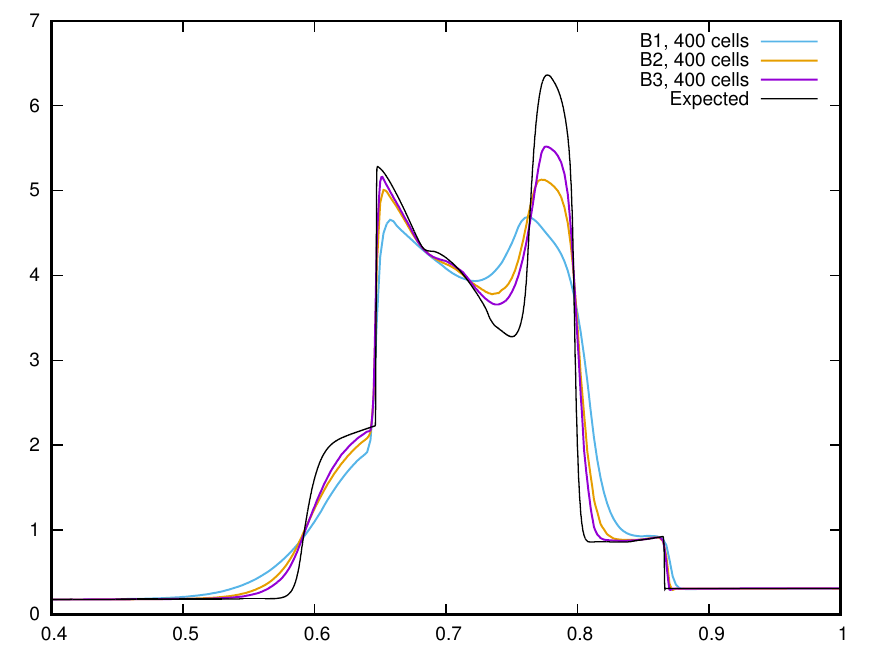}}
\subfigure[$800$ cells]{\includegraphics[width=0.45\textwidth]{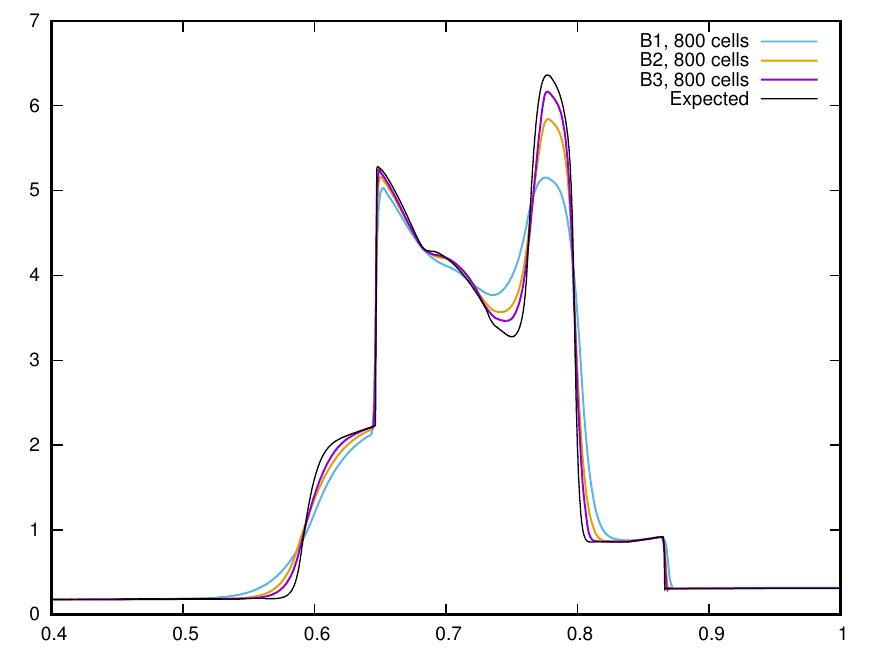}} 
\caption{Comparison between $B1,\;B2,\,B3$ for the Woodward-Colella test case at $T=0.038$}
\label{Fig:WC} 
\end{center}
\end{figure}

The interaction of blast waves is a standard low energy benchmark problem involving strong shocks reflecting from the walls of the tube with further mutual interactions. The initial data is the following:
\begin{equation*}
(\rho_0,u_0,p_0) =
\begin{cases}
[1, 0, 10^3], \; &0 \leq x \leq 0.1, \\
[1, 0, 10^{-2}], \; &0.1 < x < 0.9, \\
[1, 0, 10^2], \; &0.9 \leq x \leq 1.
\end{cases}
\end{equation*}
\rev{The jump stabilization \eqref{phi_burman} parameters have been set as in Section ~\ref{Isoflow_1D}.}
The results of the simulations using second ('B1'), third ('B2') and fourth ('B3') order RD schemes are illustrated in Fig.~\ref{Fig:WC}. The plots show a very good behavior of the numerical scheme even for this extremely demanding test case. The solution is well approximated already on a 400 cell mesh with B3, and further mesh refinement shows the expected convergence to the exact solution.

\subsubsection{Shu-Osher problem}

This test case, introduced in \cite{shuOsher1989}, is intended to demonstrate the advantages of high order schemes for problems involving some structure in smooth regions. In this test, we solve the Euler equations with initial conditions containing a moving Mach 3 shock wave which later interacts with periodic perturbations in density. The initial data for this problem is defined as follows:
\begin{equation*}
W=[\rho,u,p]=
\begin{cases}
[3.857143, 2.629369, 10.333333], \; &-5 \leq x \leq -4, \\
[1 + 0.2\sin(5x), 0, 1], \; &-4 < x \leq 5.
\end{cases}
\end{equation*}
\rev{Also in this case, the jump stabilization \eqref{phi_burman} parameters have been set as in Section ~\ref{Isoflow_1D}.}
The more accurate approximation obtained by the fourth order scheme in comparison to the second and third order is clearly visible in this benchmark problem, and increasing the number of mesh elements within the domain strongly increases the quality of the solution.
\begin{figure}[h!]
\begin{center}
\subfigure[$400$ cells]{\includegraphics[width=0.45\textwidth]{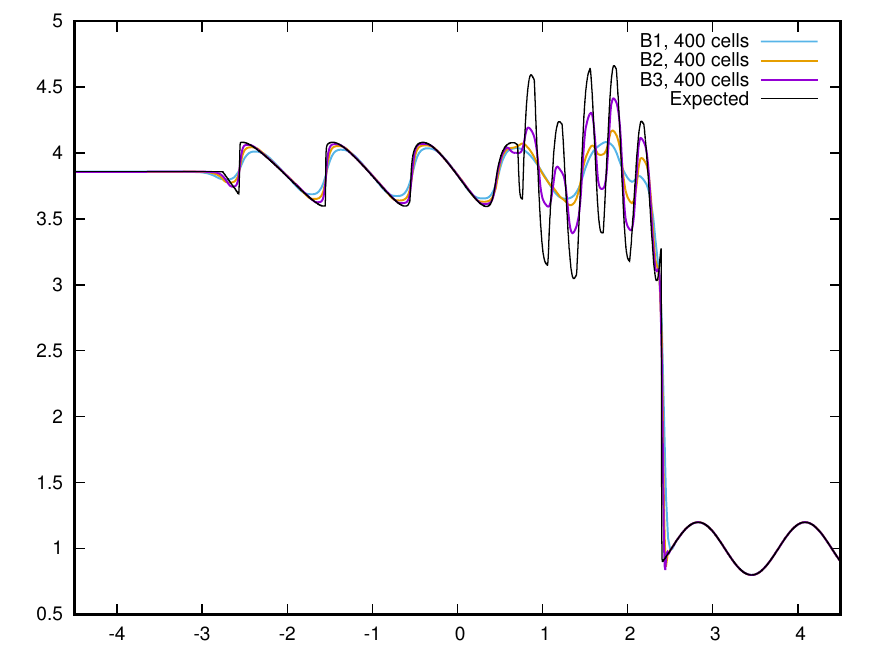}}
\subfigure[$800$ cells]{\includegraphics[width=0.45\textwidth]{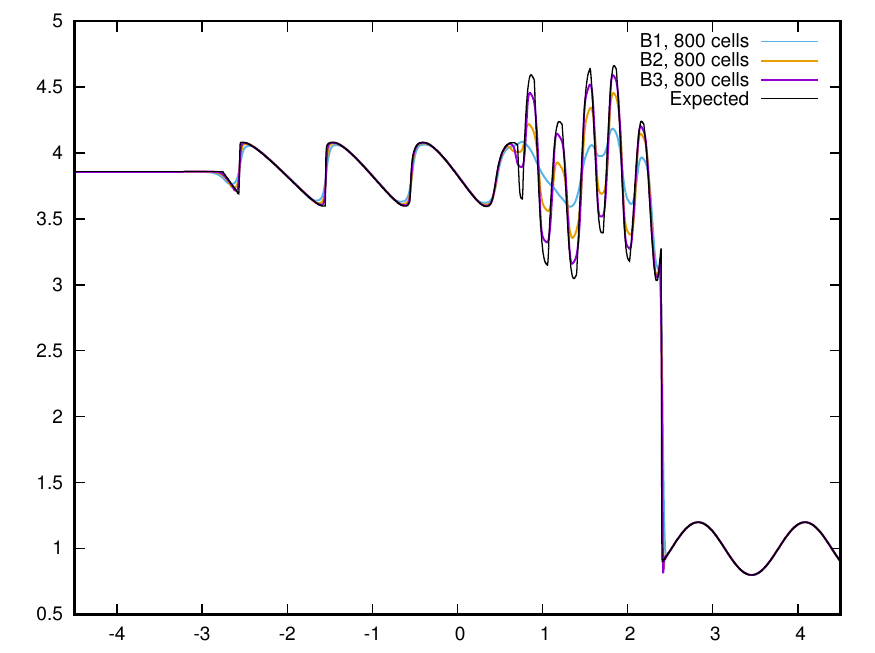}}
\end{center}   
\caption{Comparison between $B1,\;B2,\,B3$ for the Shu-Osher test case at $T=1.8$}
\label{Fig:SO}
\end{figure}
\clearpage
\rev{This benchmark problem allows, further, to highlight the dissipative character of the original local Lax-Friedrichs approximation \eqref{phi_LxF_xI} when dealing with higher than second order of accuracy and the improvement one may achieve with its reformulation in terms of a telescopic sum as in \eqref{phi_LxF_sub}. On a grid with 200 nodes we can observe in Fig.~\ref{Fig:SO_before} how the approximation has extremely clipped extrema when the proposed method with \eqref{phi_LxF_xI} is applied, while with \eqref{phi_LxF_sub} in Fig.~\ref{Fig:SO_after} results in sharper extrema on such a coarse mesh.
\begin{figure}[h!]
\begin{center}
\subfigure[Formulation \eqref{phi_LxF_xI}]{\includegraphics[width=0.45\textwidth]{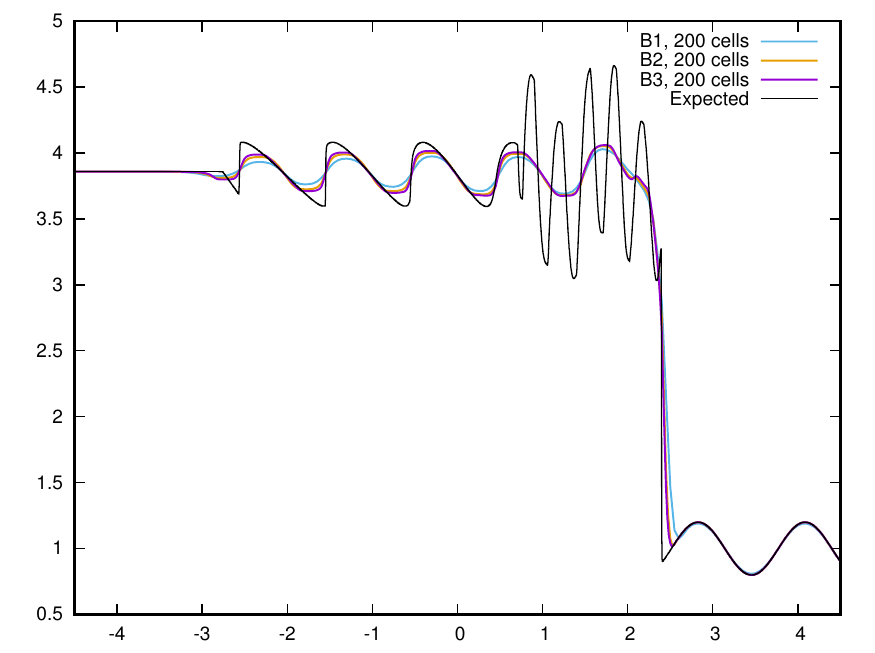}}\label{Fig:SO_before}
\subfigure[formulation \eqref{phi_LxF_sub}]{\includegraphics[width=0.45\textwidth]{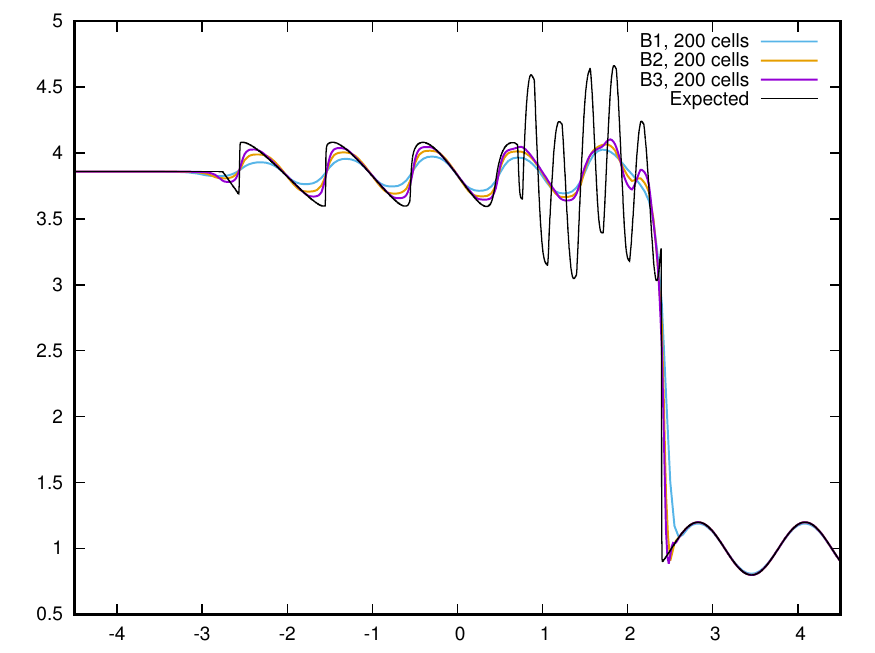}}\label{Fig:SO_after}\\  
\subfigure[Formulation \eqref{phi_LxF_xI} - zoom ]{\includegraphics[width=0.45\textwidth]{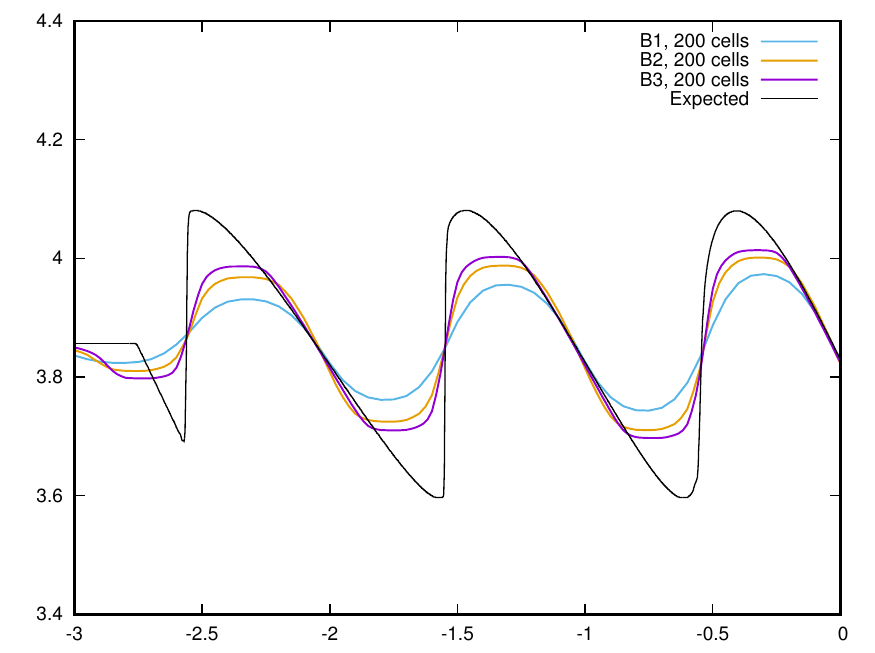}}\label{Fig:SO_before}
\subfigure[formulation \eqref{phi_LxF_sub} - zoom ]{\includegraphics[width=0.45\textwidth]{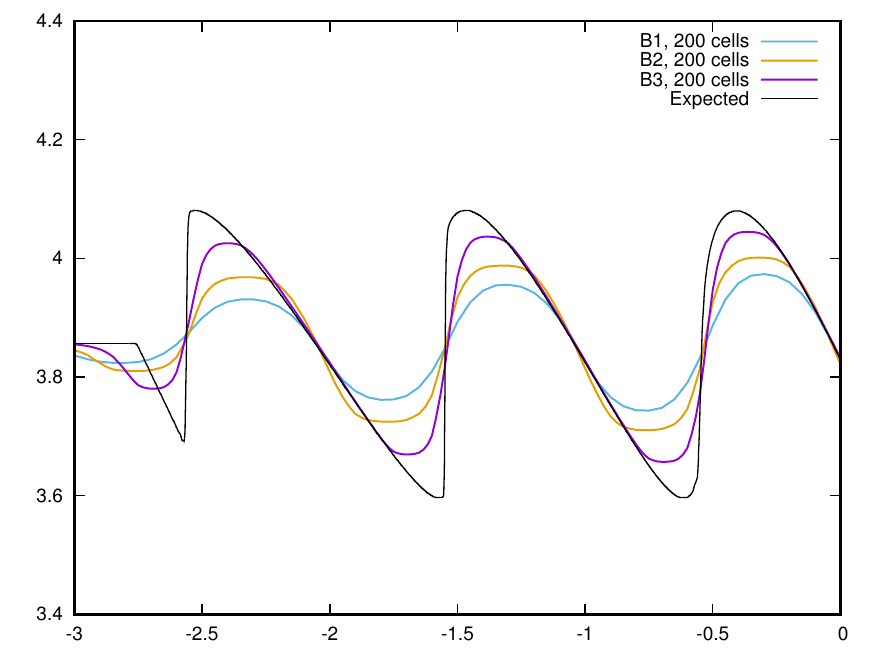}}\label{Fig:SO_after}  
\end{center} 
\caption{\rev{Comparison of the different Lax-Friedrichs approximations between $B1,\;B2,\,B3$ for the Shu-Osher test case at $T=1.8$. Compare with Fig. \ref{Fig:SO}, where we use a finer resolution for formulation \eqref{phi_LxF_sub}.}}
\label{Fig:SO_200comp}
\end{figure}
}
\clearpage

\subsection{Numerical results for 2D test cases}

\subsubsection{Stationary vortex}\label{Section:Vortex}

The first considered test case in 2D is the stationary isentropic vortex evolution problem, see e.~g. \cite{Yee1999}.
Initially, an isentropic perturbation $(\delta S=0)$ is applied to the system, such that
\begin{equation}
\label{isentropic:delta}
\begin{cases}
\delta u=-\bar{y}\frac{\beta}{2\pi}e^{(1-r^2)/2}\\[0.3em]
\delta v=\bar{x}\frac{\beta}{2\pi}e^{(1-r^2)/2}\\[0.3em]
\delta T =-\frac{(\gamma-1)\beta^2}{8\gamma\pi^2}e^{1-r^2}.
\end{cases}
\end{equation}
The initial conditions are thus set to
\begin{equation}\label{isentropic:IC}
\begin{cases}
\rho=T^{1/(\gamma-1)}=(T_{\infty}+\delta T)^{1/(\gamma-1)}=\Big{[}1-\frac{(\gamma-1)\beta^2}{8\gamma\pi^2}e^{1-r^2}\Big{]}^{1/(\gamma-1)}\\[0.3em]
\rho u=\rho(u_{\infty}+\delta u)=\rho \Big{[}1-\bar{y} \frac{\beta}{2\pi}e^{(1-r^2)/2}\Big{]}\\[0.3em]
\rho v=\rho(v_{\infty}+\delta v)=\rho \Big{[}1+\bar{x} \frac{\beta}{2\pi}e^{(1-r^2)/2}\Big{]}\\[0.3em]
\rho E=\frac{\rho^\gamma}{\gamma-1}+\frac{1}{2}\rho (u^2+v^2),
\end{cases}
\end{equation}
where $\beta=5$, $r=\bar{x}^2+\bar{y}^2$ and $(\bar{x},\bar{y})=\big{(}(x-x_0),(y-y_0)\big{)}$. The computational domain is a circle with radius of $10$ and center at $(0,0)$. The center of the vortex is set in $(x_0,y_0)=(0,0)$. The parameters of the unperturbed flow are set to $u_{\infty}=v_{\infty}=0$ for the velocities, $p_{\infty}=1$ for the pressure,  $T_{\infty}=1$ for the temperature and $\rho_{\infty}=1$ for the density.

\begin{figure}[H]
\begin{center}
\subfigure[$\rho$]{\includegraphics[width=0.45\textwidth]{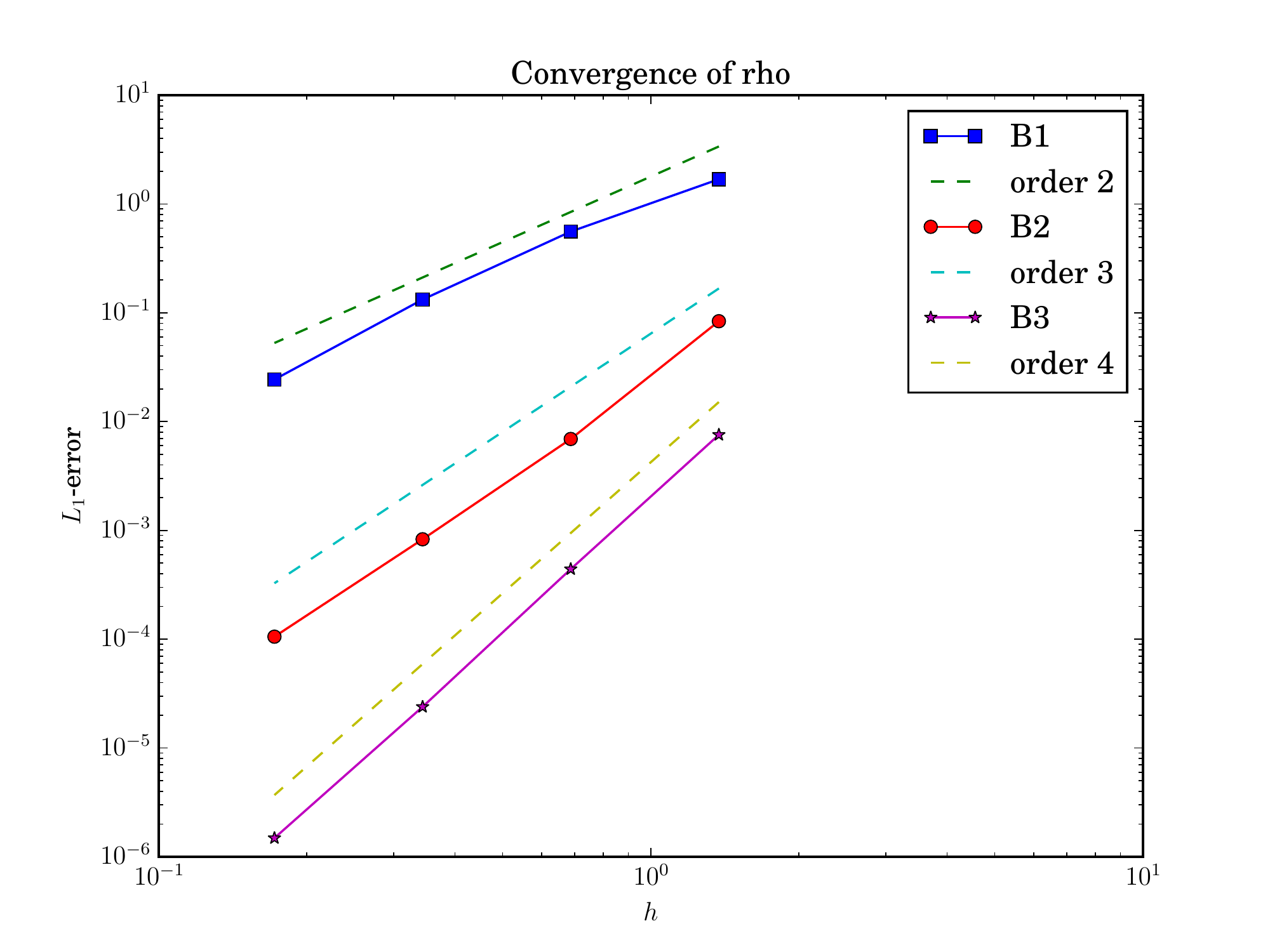}}
\subfigure[$u$]{\includegraphics[width=0.45\textwidth]{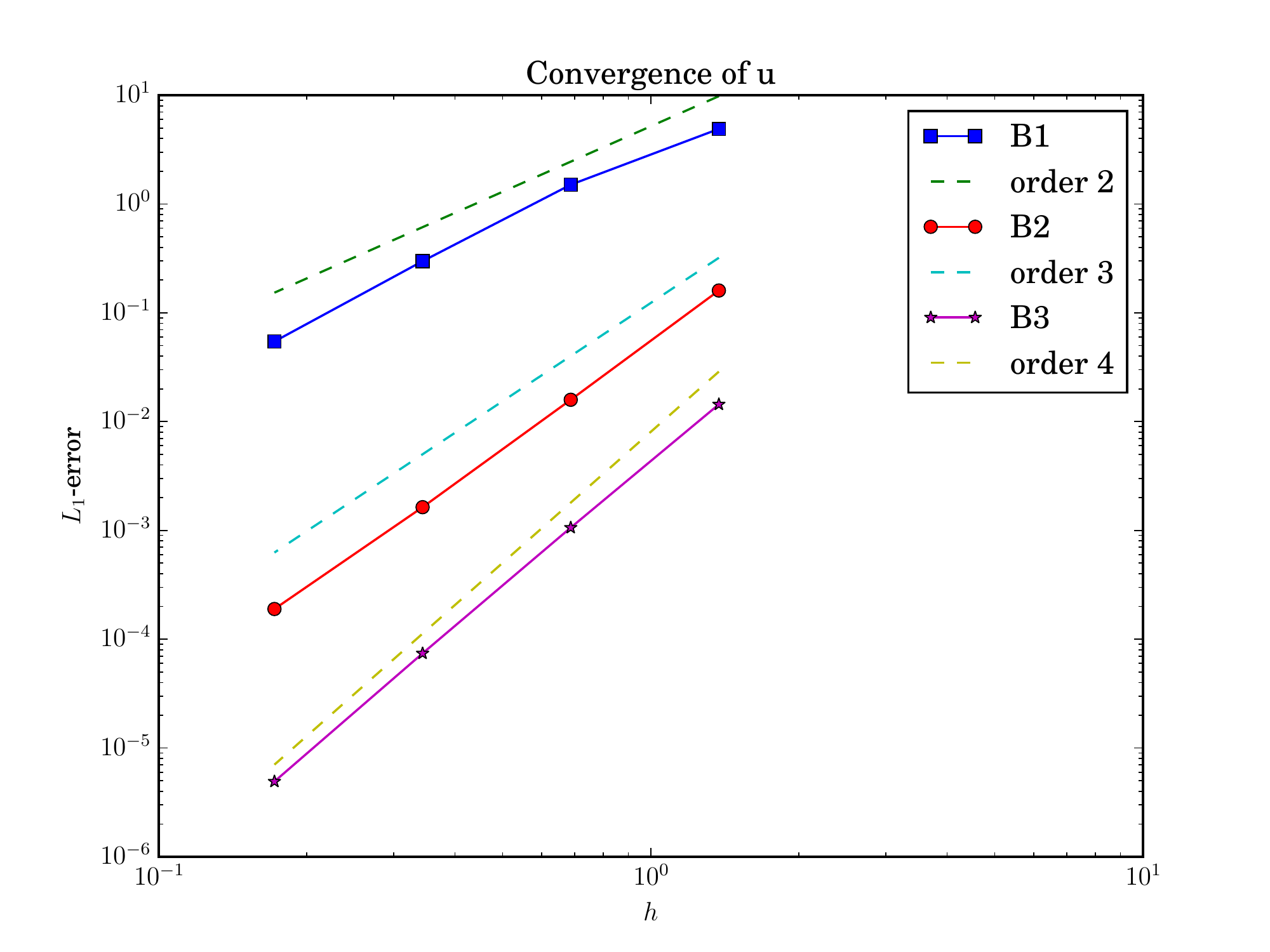}} \\
\subfigure[$p$]{\includegraphics[width=0.45\textwidth]{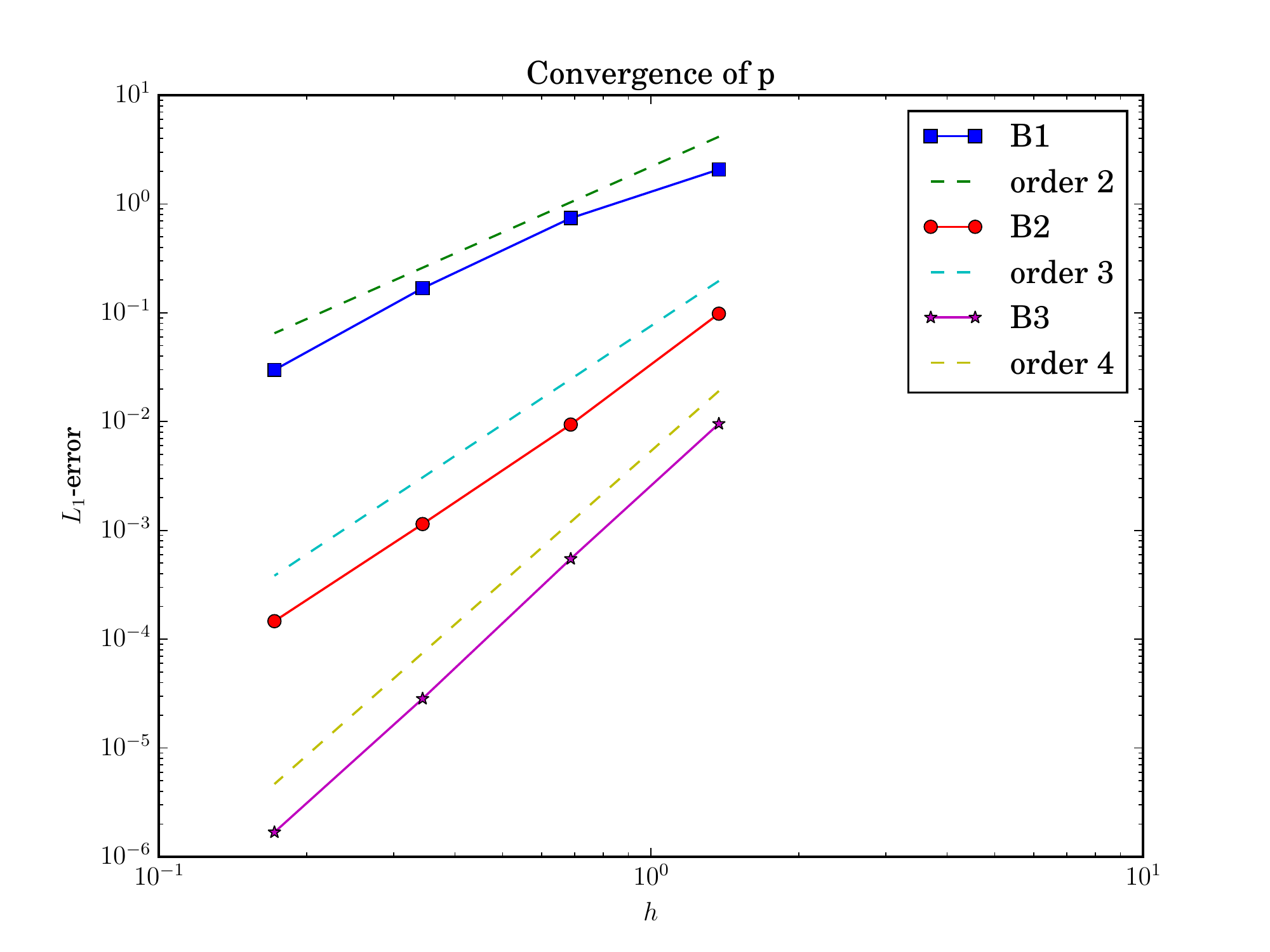}}
\subfigure[$v$]{\includegraphics[width=0.45\textwidth]{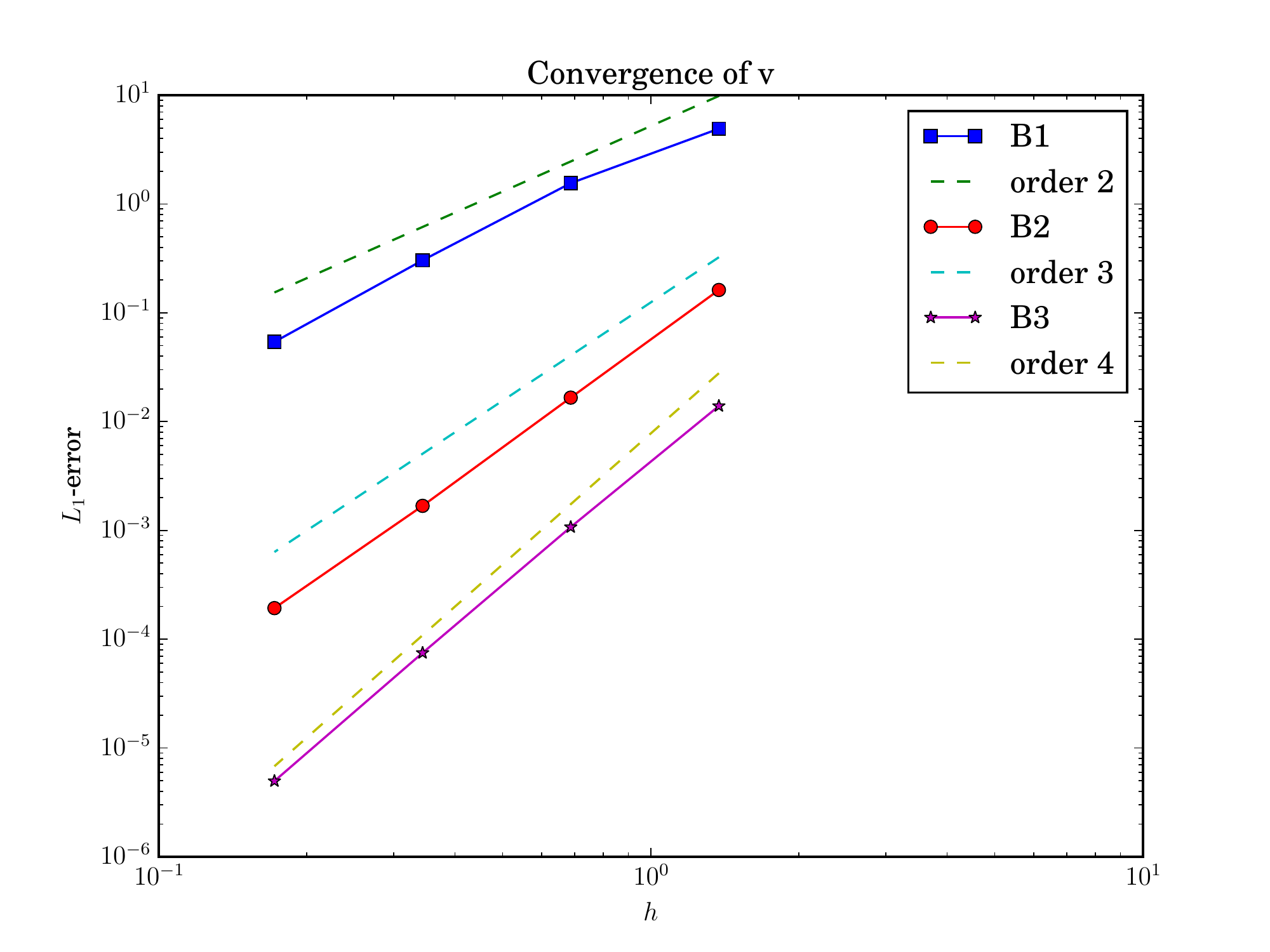}} 
\caption{2D stationary vortex convergence study.}\label{vortex_convergence_2D}
\end{center}
\end{figure}

\begin{figure}[H]
\begin{center}
\subfigure[$B1$]{\includegraphics[width=0.4\textwidth]{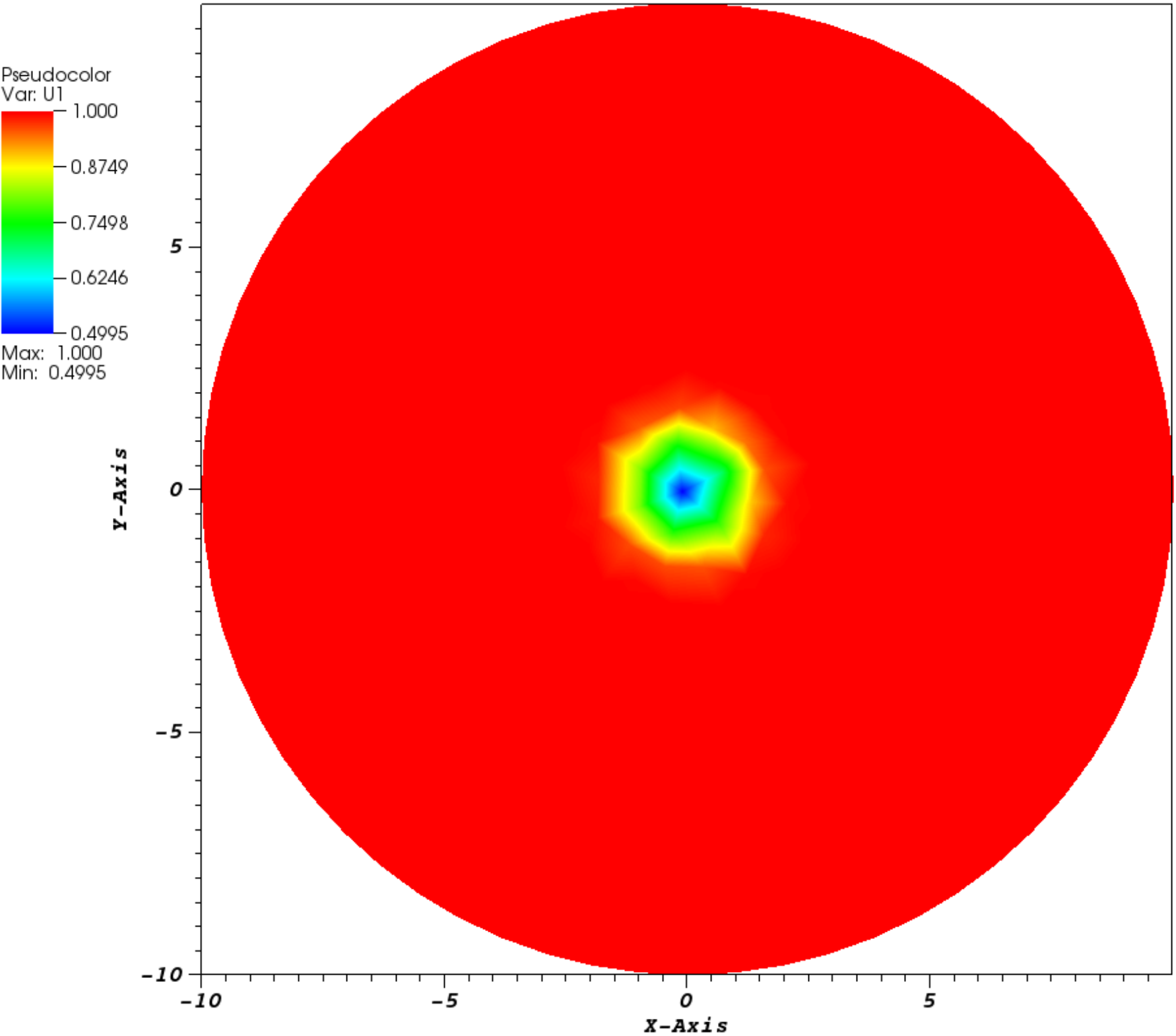}}
\subfigure[$B2$]{\includegraphics[width=0.4\textwidth]{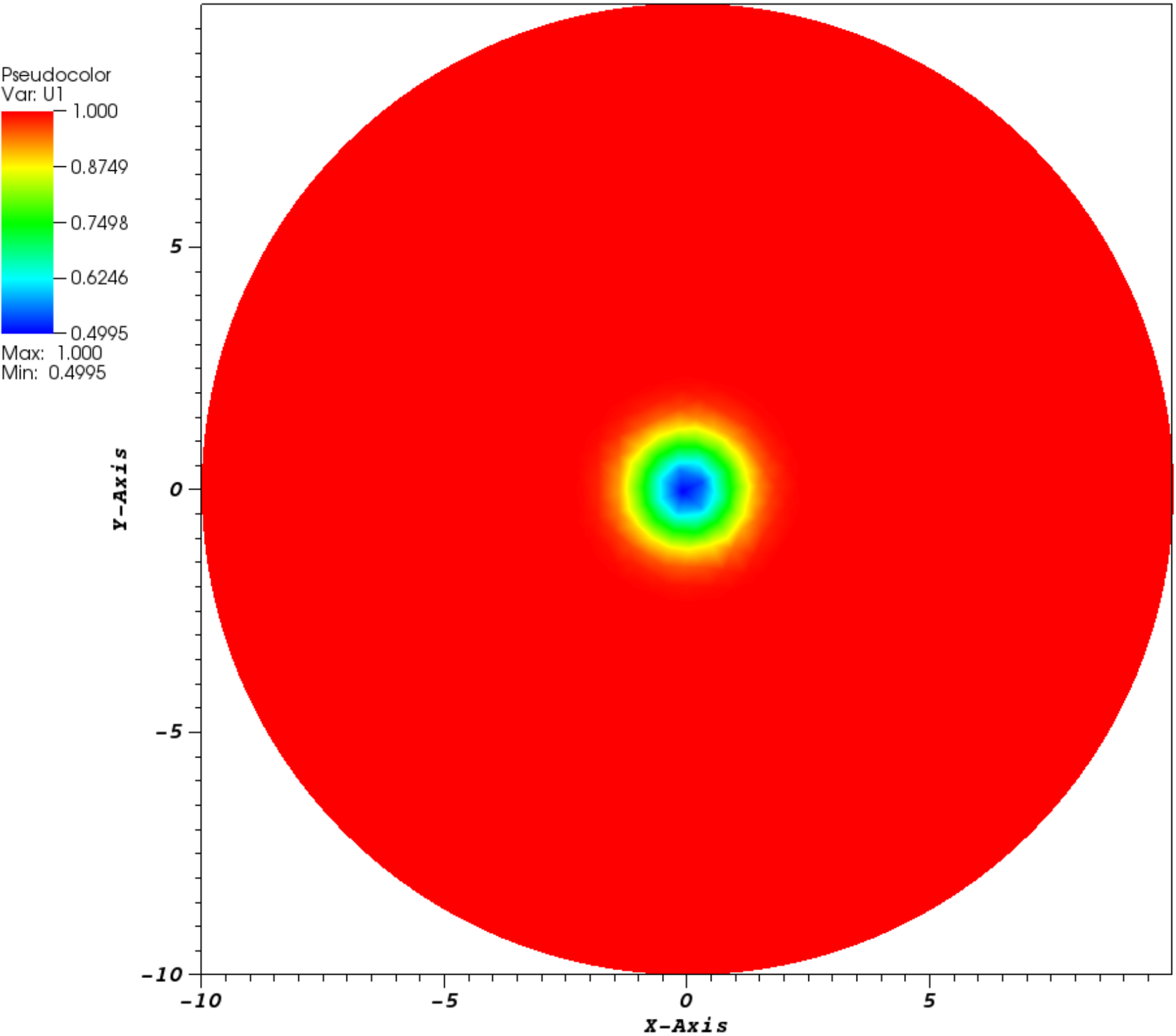}}
\subfigure[$B3$]{\includegraphics[width=0.4\textwidth]{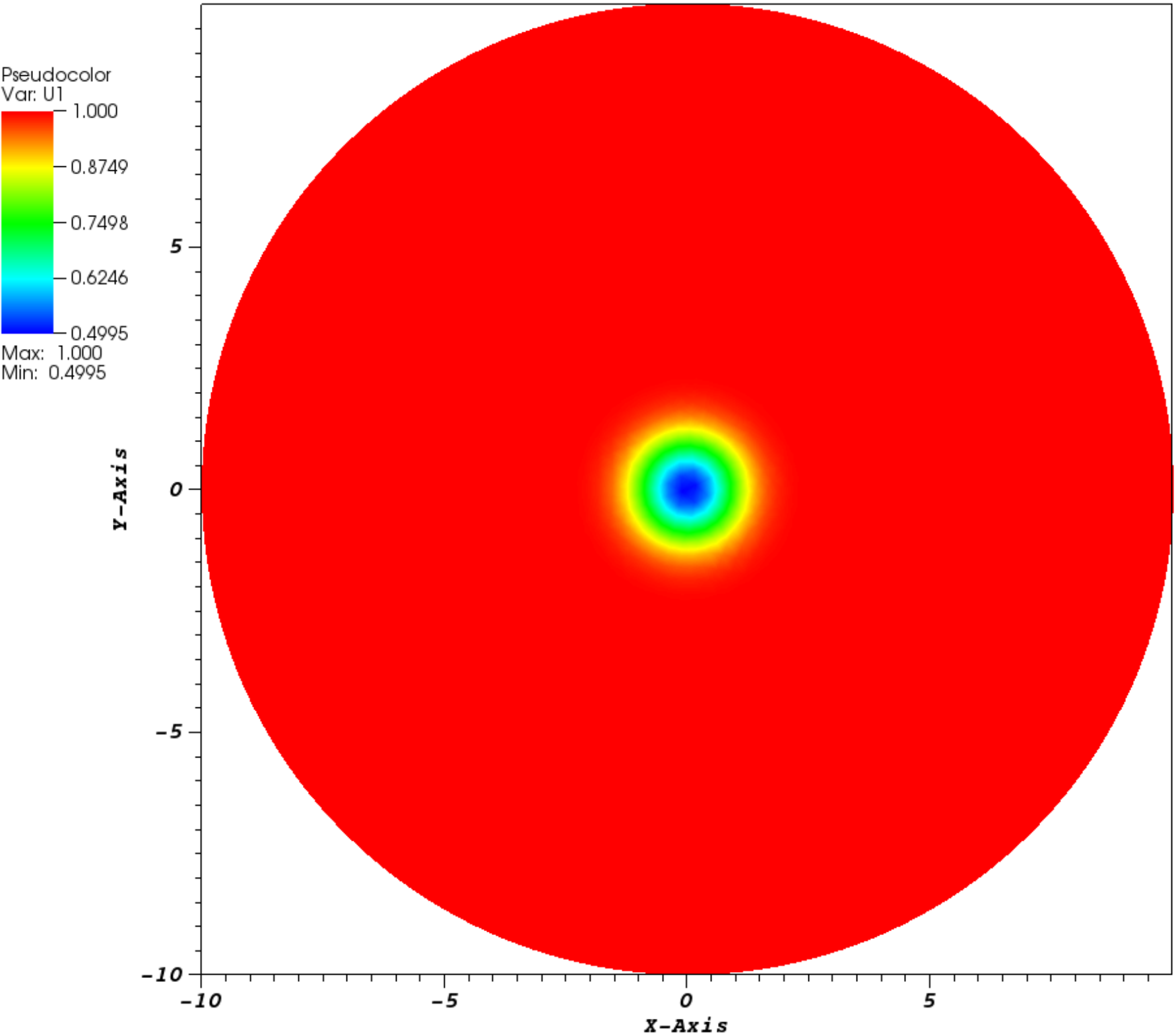}} 
\subfigure[exact]{\includegraphics[width=0.4\textwidth]{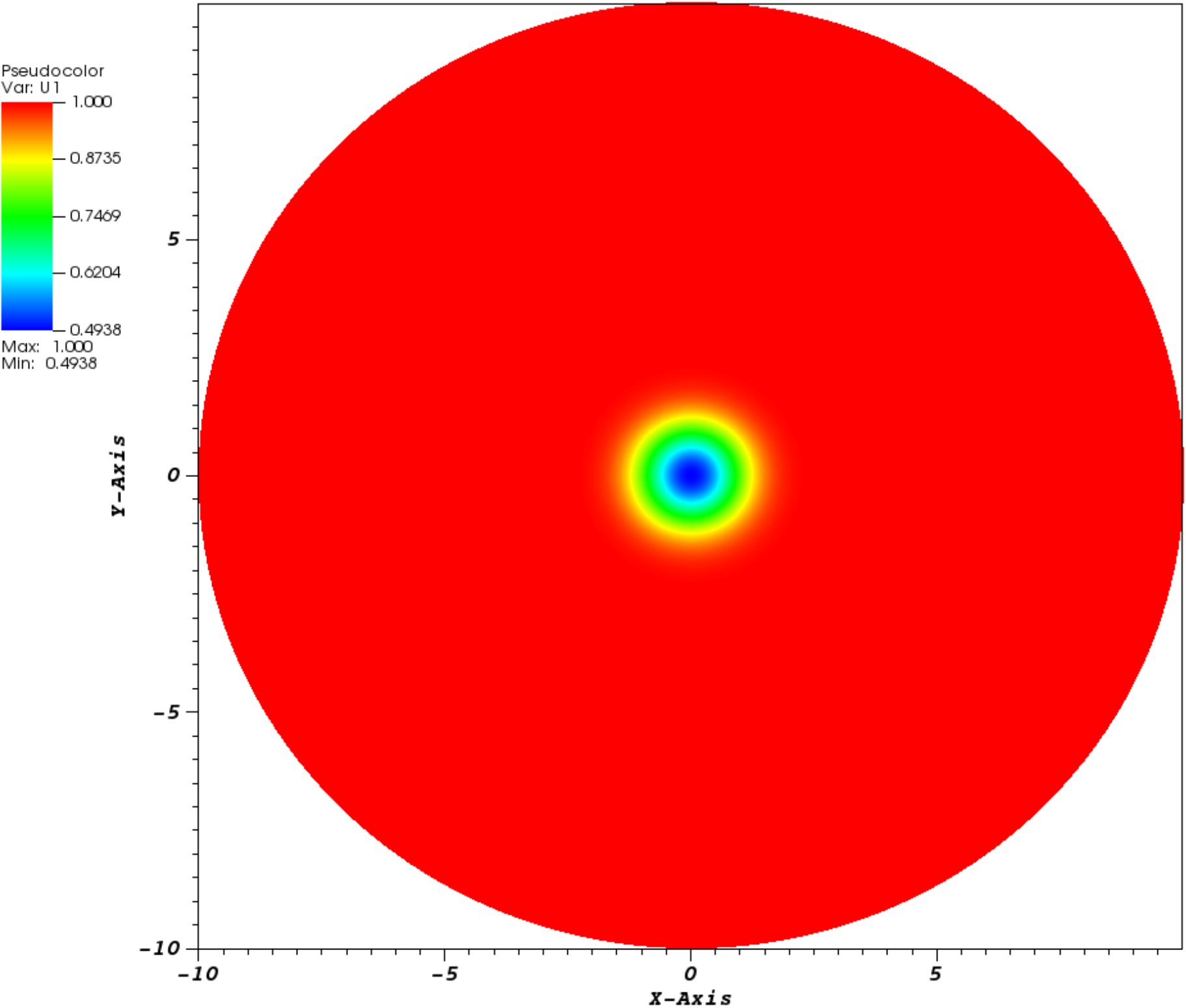}} 
\caption{2D stationary vortex for the density on a mesh given by $N_0=934$ elements at $T=50$.}\label{vortex_rho_2D}
\end{center}
\end{figure}

\begin{figure}[H]
\begin{center}
\includegraphics[width=0.4\textwidth]{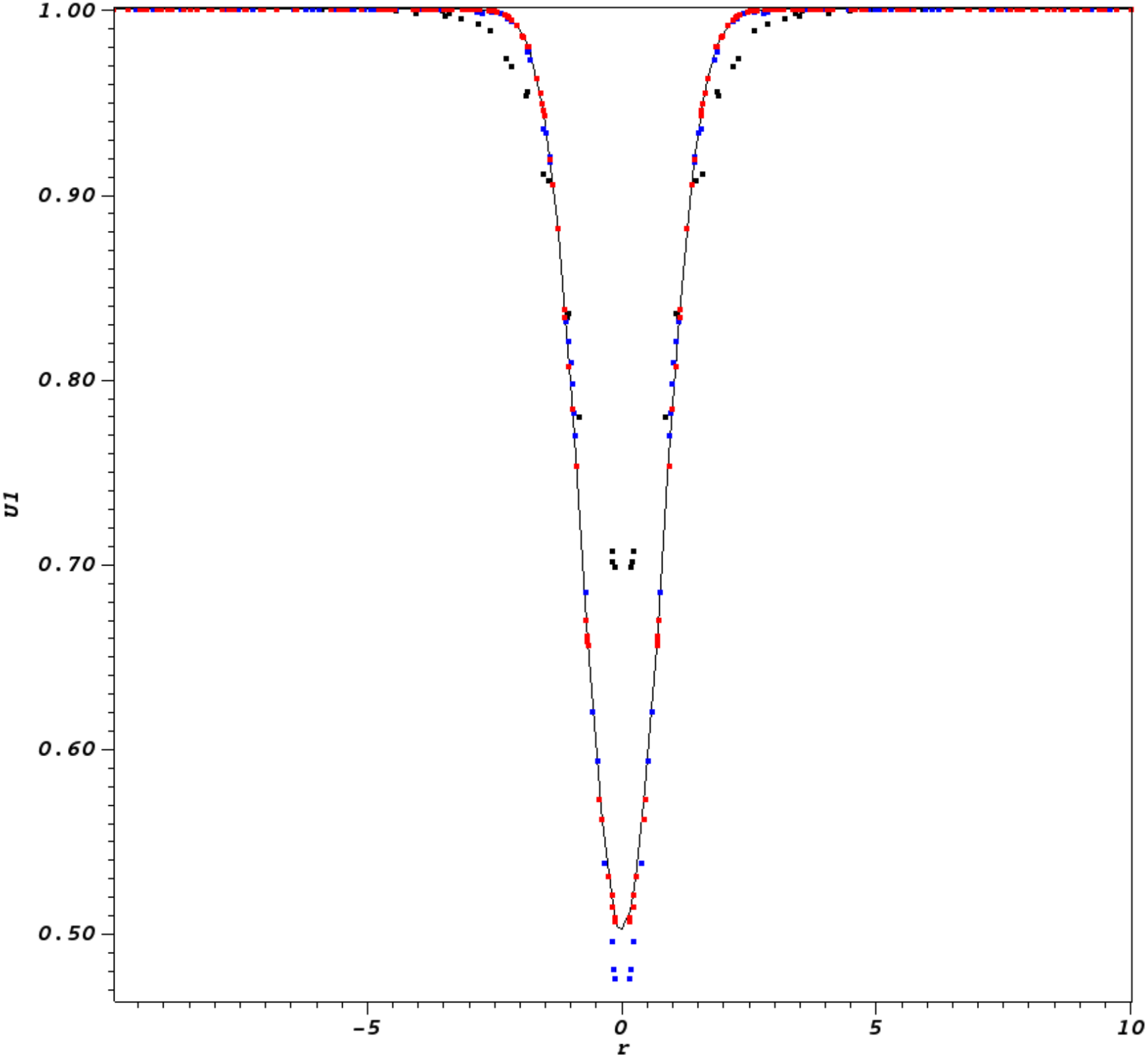}
\caption{2D stationary vortex scatter plot for the density on a mesh given by $N_0=934$ elements. Results for $B1$ (dashed black), $B2$ (dashed blue) and $B3$ (dashed red) and the exact (continuous black) at $T=50$.}\label{vortex_scatter_2D}
\end{center}
\end{figure}

The proposed RD method has been tested with $B1$, $B2$ and $B3$ basis polynomials on four different meshes with the number of triangular elements equal to $N_0=608$, $N_1=934$, $N_2=14176$ and $N_3=56192$, respectively\footnote{corresponding to typical element sizes $h=0.7$, $0.57$, $0.14$ and $0.075$, respectively}.
\rev{In order to stabilize the approximation, we set w.r.t. \eqref{phi_burman} for B1 $\theta_1=0.1$ and $\theta_2=0$, for B2 $\theta_1=0.01$ and $\theta_2=0$ and for B3 $\theta_1=0.001$ and $\theta_2=0.$}
The obtained convergence curves for the three schemes considered are displayed in Figure \ref{vortex_convergence_2D} and show an excellent correspondence to the theoretically predicted order. From Figures \ref{vortex_rho_2D} and \ref{vortex_scatter_2D} it can be seen that a very good approximation is already achieved on the mesh with $934$ cells. As expected, the $B1$ approximation is the most dissipative, while $B2$ becomes more accurate despite having a small density undershoot at the center of the vortex, and $B3$ shows an excellent approximation which is practically indistinguishable  from the exact solution.

\clearpage
\clearpage

In our numerical experiments this test case represents the only exception revealing the difference between the two approximating strategies of the fluxes, as described in Section \ref{Section_Extend_sys}.
From a formal point of view, it is equivalent to interpolate the flux or to compute the flux for the solution interpolated at a quadrature point. However, from a computational point of view, the first strategy leads to a quadrature-free algorithm, while the second one requires a quadrature formula. For this test case, the quadrature-free strategy proved to be unstable while the second one led to the results presented above. We have never encountered this issue in the 1D simulations, which have all been computed with the quadrature-free approximation. In order to be consistent with the presented results, the 2D test cases have accordingly been computed using the second approximation strategy. \rev{We note that instability of the quadrature-free version can be attributed to a violation of stability
conditions derived in the recent paper \cite{Barrenechea}.}

\subsubsection{2D Sod problem}\label{Section:2D_Sod}
Further, we have tested our high order RD scheme on a well-known 2D Sod benchmark problem. The initial conditions are given by
\begin{equation*}
(\rho_0,u_0,v_0,p_0) = 
\begin{cases}
[1, 0, 0, 1], \; &0 \leq r \leq 0.5, \\
[0.125, 0, 0, 0.1], \; &0.5 < r \leq 1,
\end{cases}
\end{equation*}
where $r = \sqrt{x^2+y^2}$ is the distance of the point $(x,y)$ from the origin.

\begin{figure}[H]
\begin{center}
\subfigure[$B1$]{\includegraphics[width=0.3\textwidth]{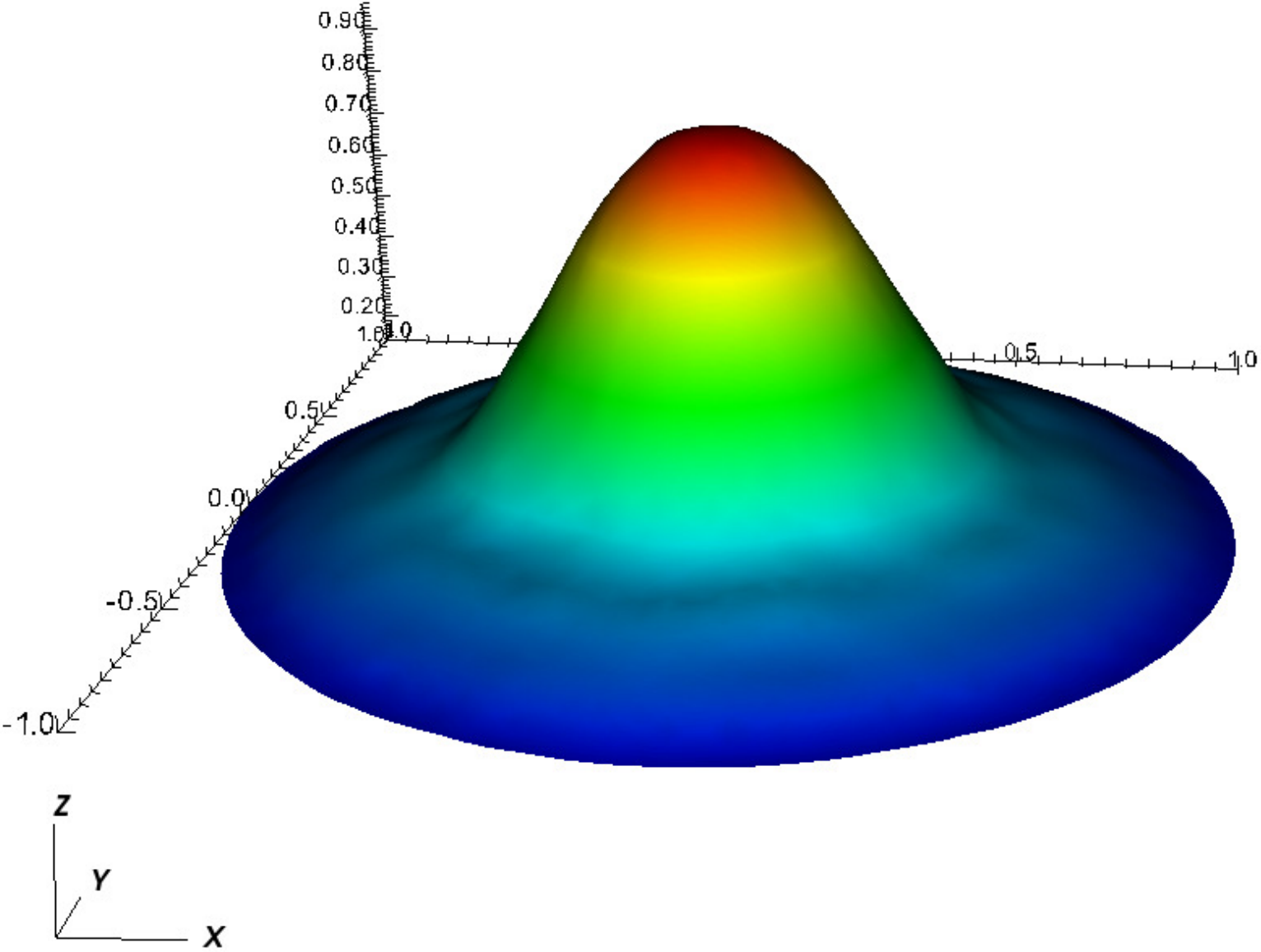}}
\subfigure[$B2$]{\includegraphics[width=0.3\textwidth]{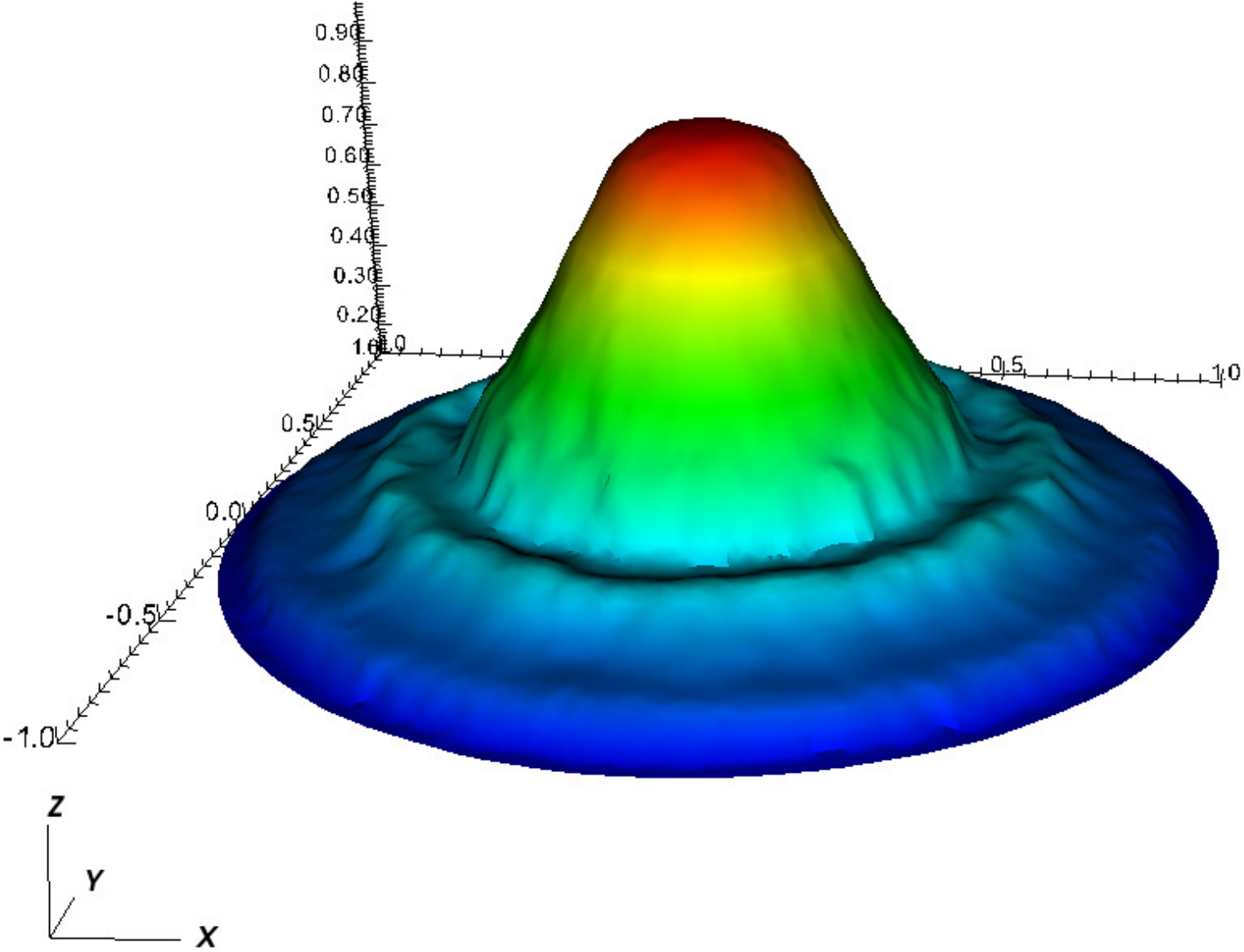}} 
\subfigure[$B3$]{\includegraphics[width=0.3\textwidth]{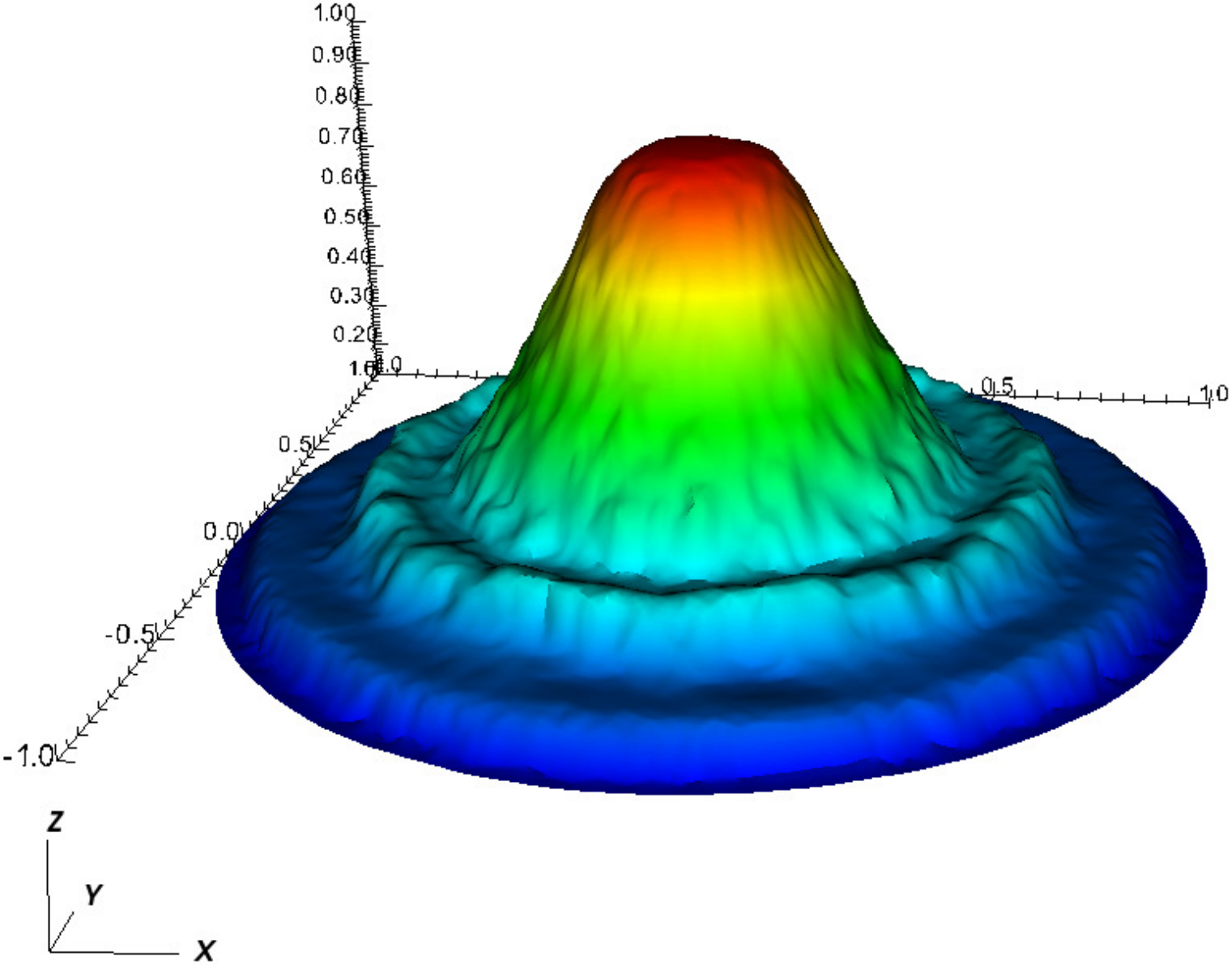}} 
\caption{2D Sod. Results at $T=0.25$ on a grid with $N=3576$ elements.}\label{Sod_2D_coarse}
\end{center}
\end{figure}

\begin{figure}[H]
\begin{center}
\subfigure[Density]{\includegraphics[width=0.3\textwidth]{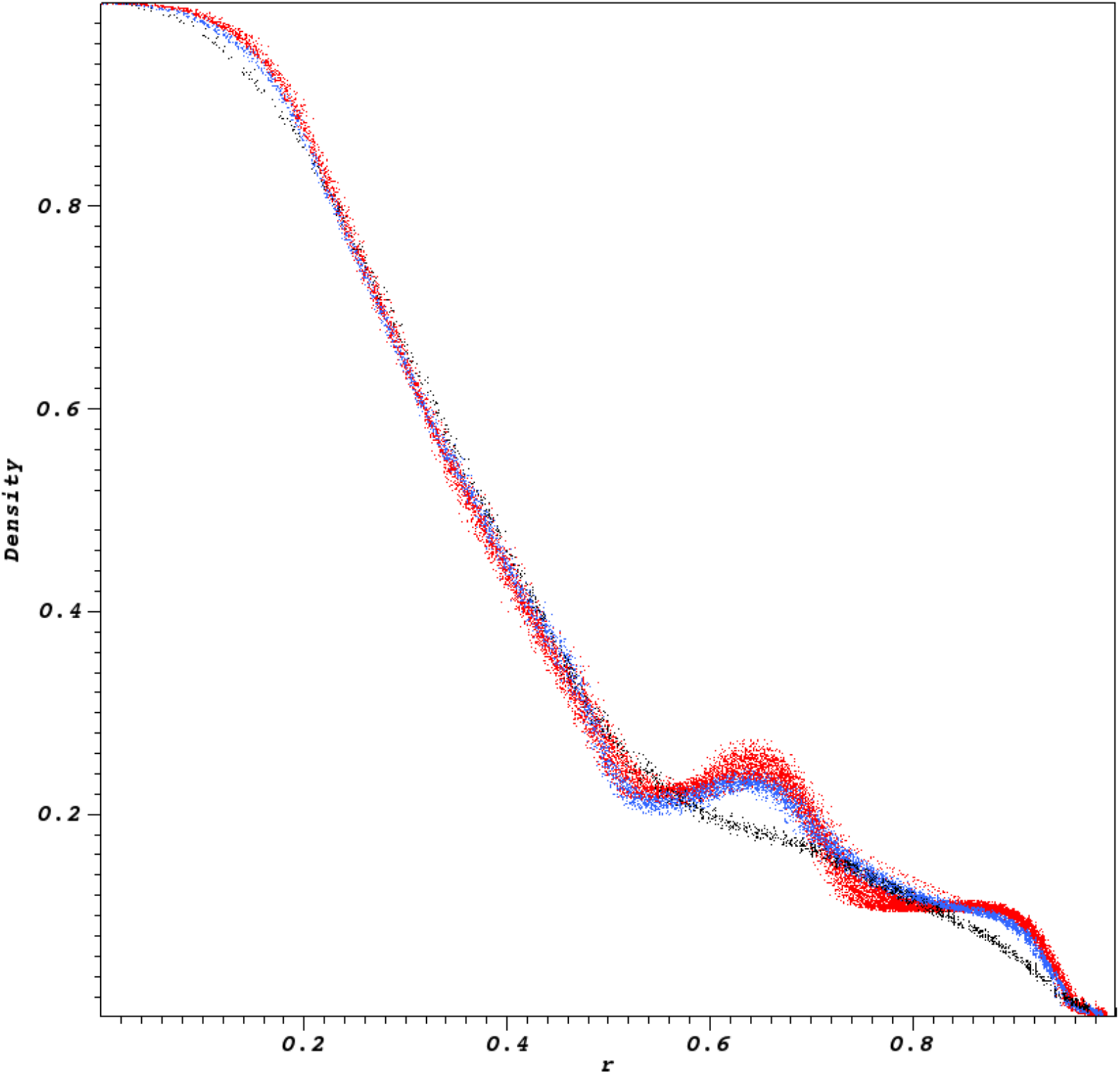}}
\subfigure[Velocity]{\includegraphics[width=0.3\textwidth]{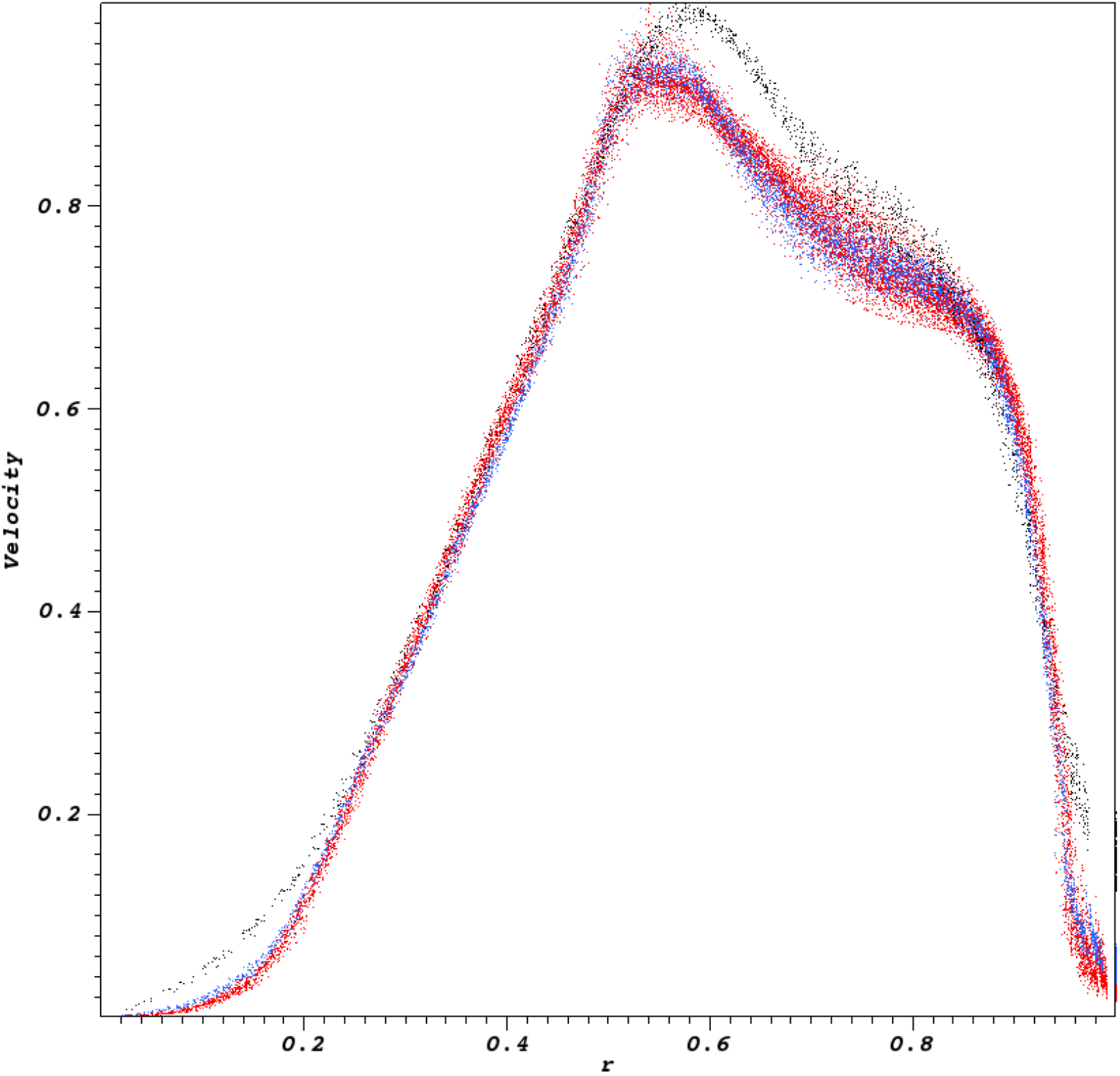}} 
\subfigure[Pressure]{\includegraphics[width=0.3\textwidth]{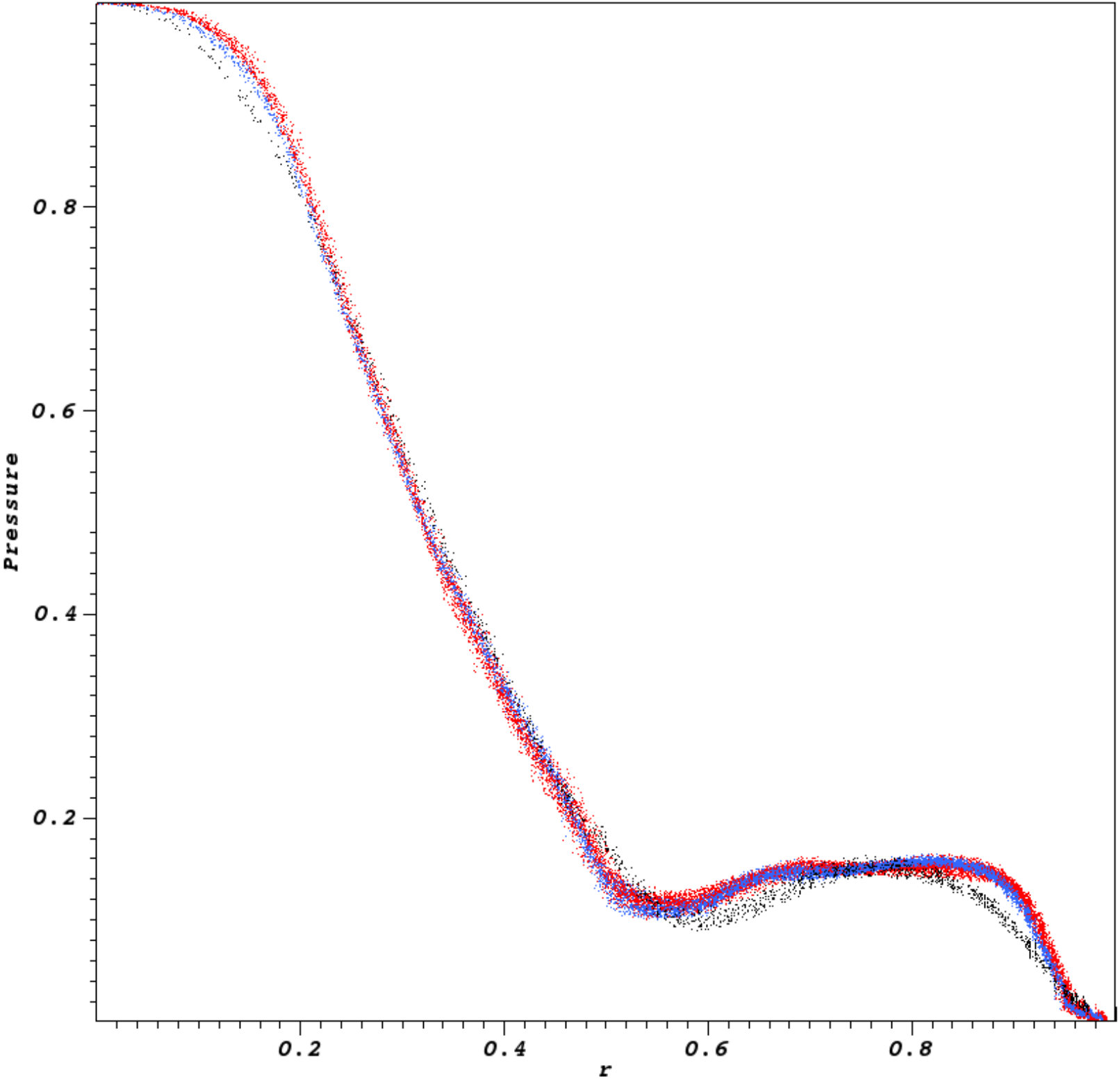}} 
\caption{2D Sod. Results for $B1$ (black), $B2$ (blue) and $B3$ (red) elements at $T=0.25$ on a grid with $N=3576$ elements. }\label{Sod_2D_scatter_coarse}
\end{center}
\end{figure}

\begin{figure}[H]
\begin{center}
\subfigure[$B1$]{\includegraphics[width=0.3\textwidth]{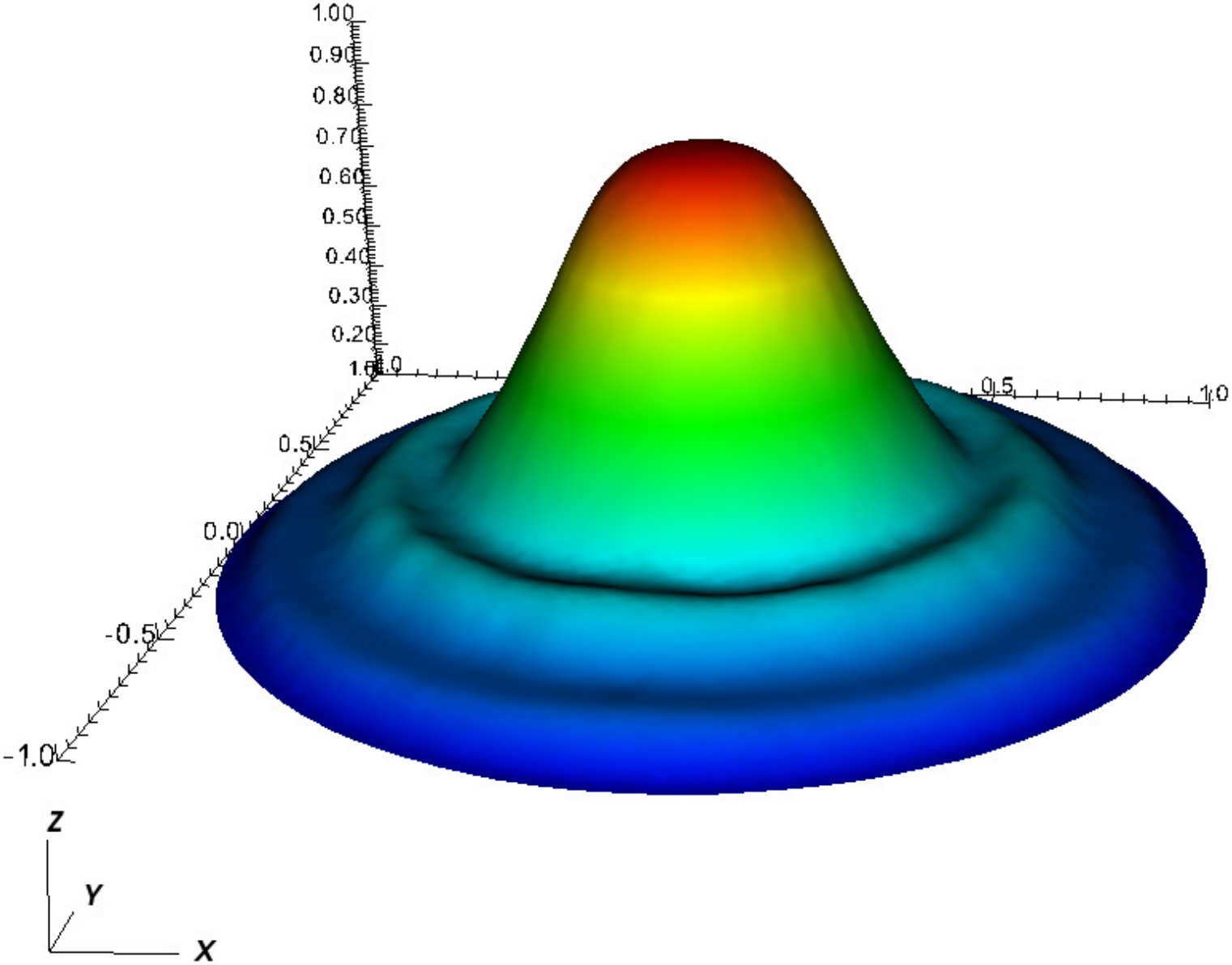}}
\subfigure[$B2$]{\includegraphics[width=0.3\textwidth]{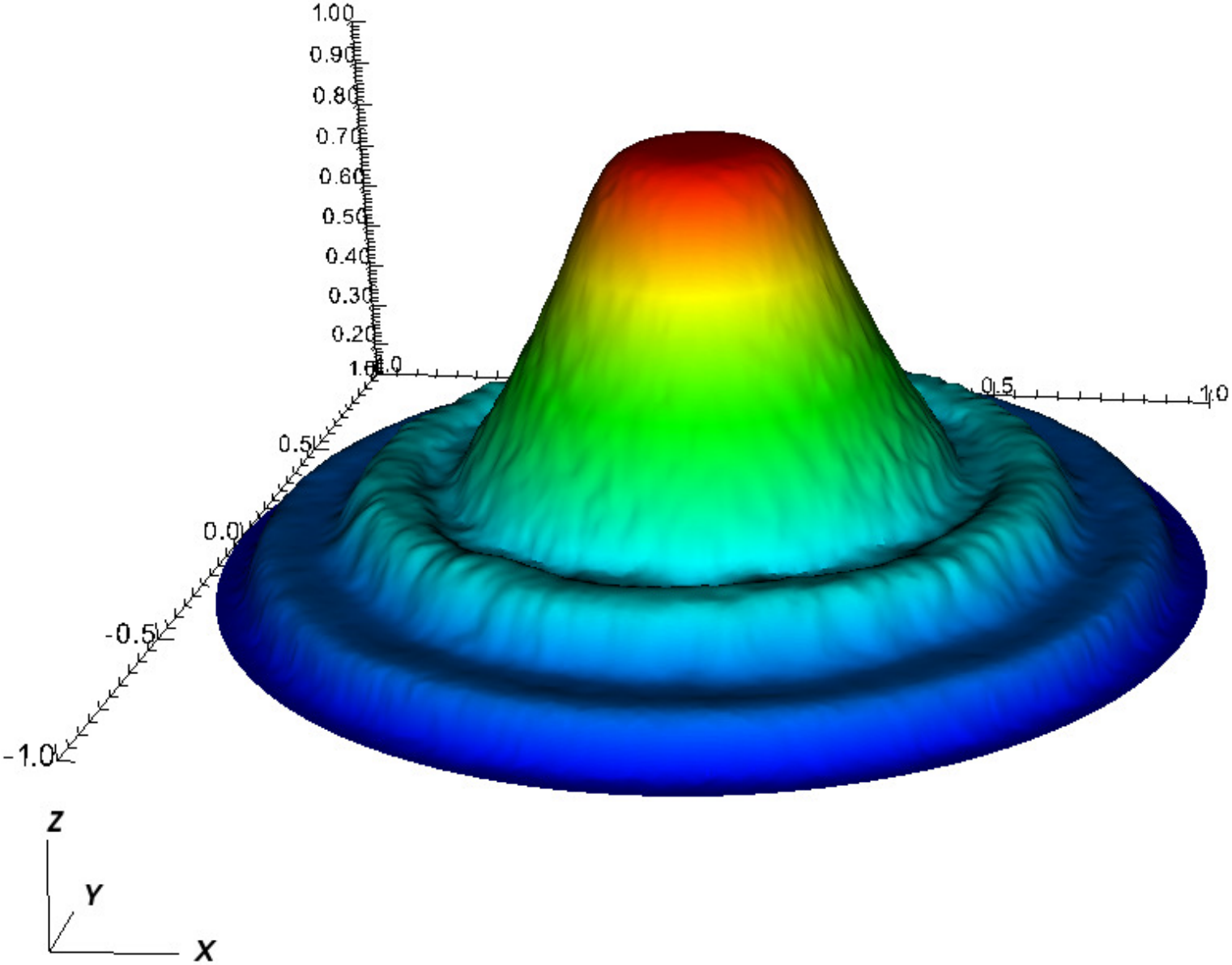}} 
\subfigure[$B3$]{\includegraphics[width=0.3\textwidth]{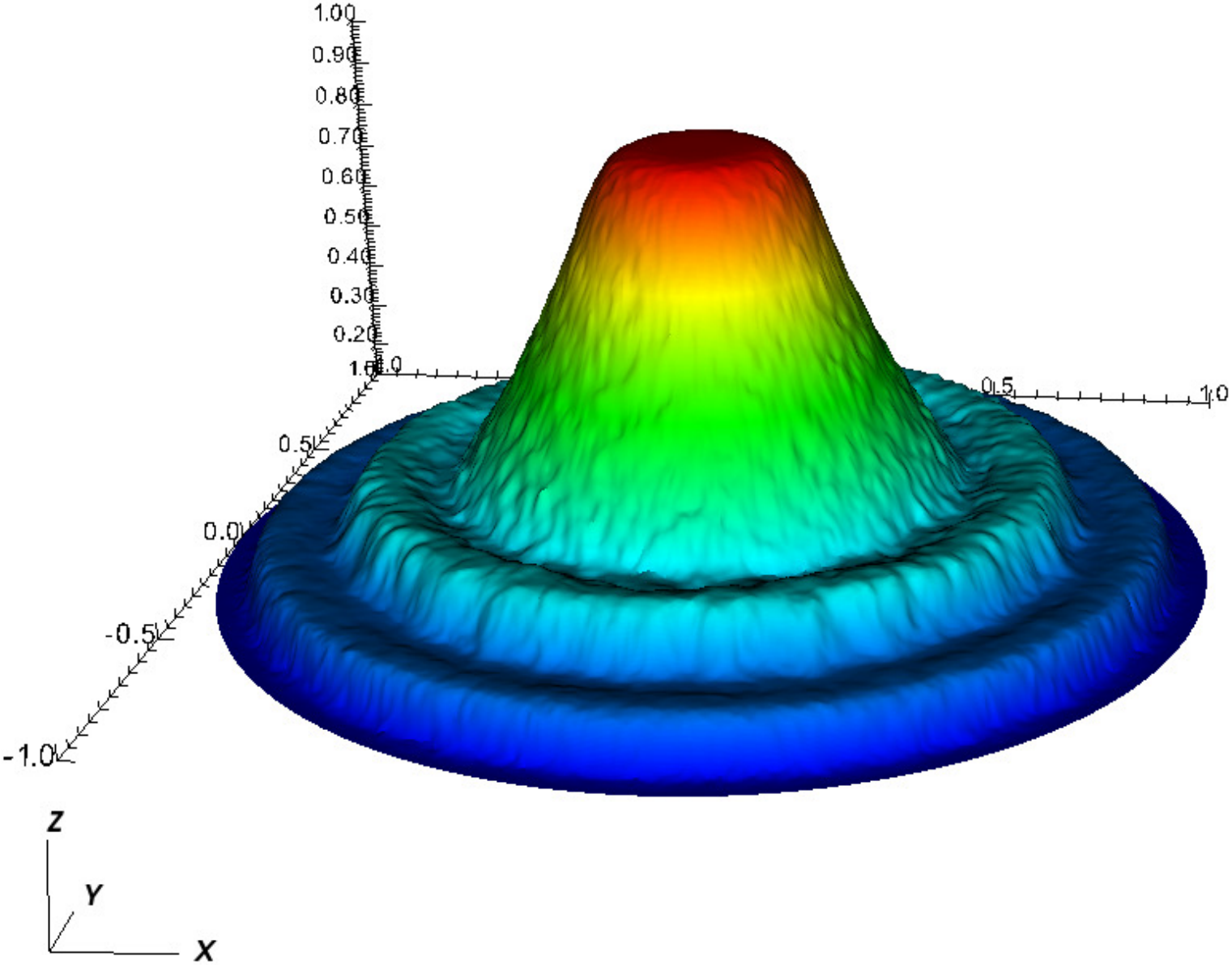}} 
\caption{2D Sod. Results at $T=0.25$ on a grid with $N=13548$ elements.}\label{Sod_2D_fine}
\end{center}
\end{figure}

\begin{figure}[H]
\begin{center}
\subfigure[Density]{\includegraphics[width=0.3\textwidth]{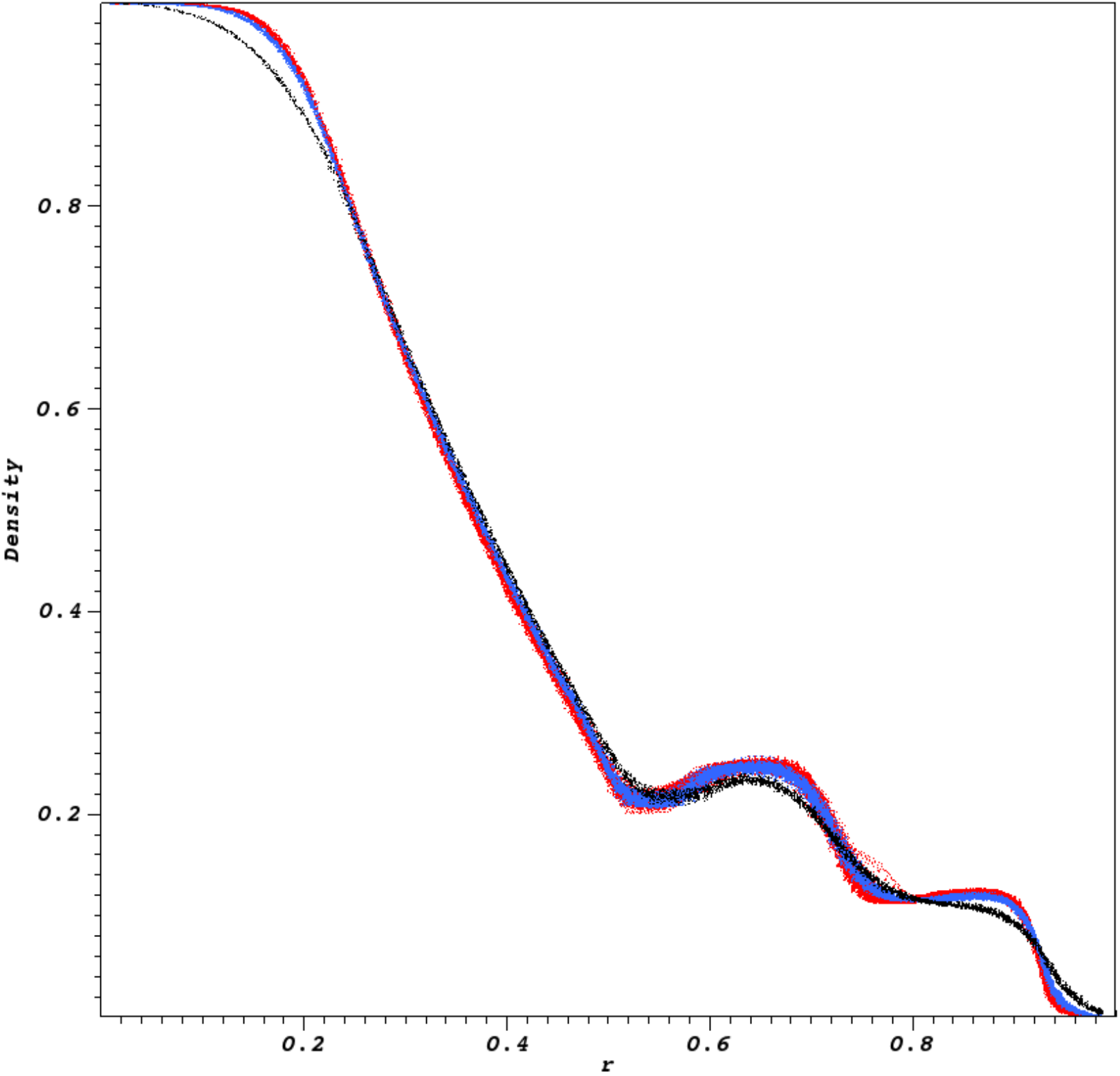}}
\subfigure[Velocity]{\includegraphics[width=0.3\textwidth]{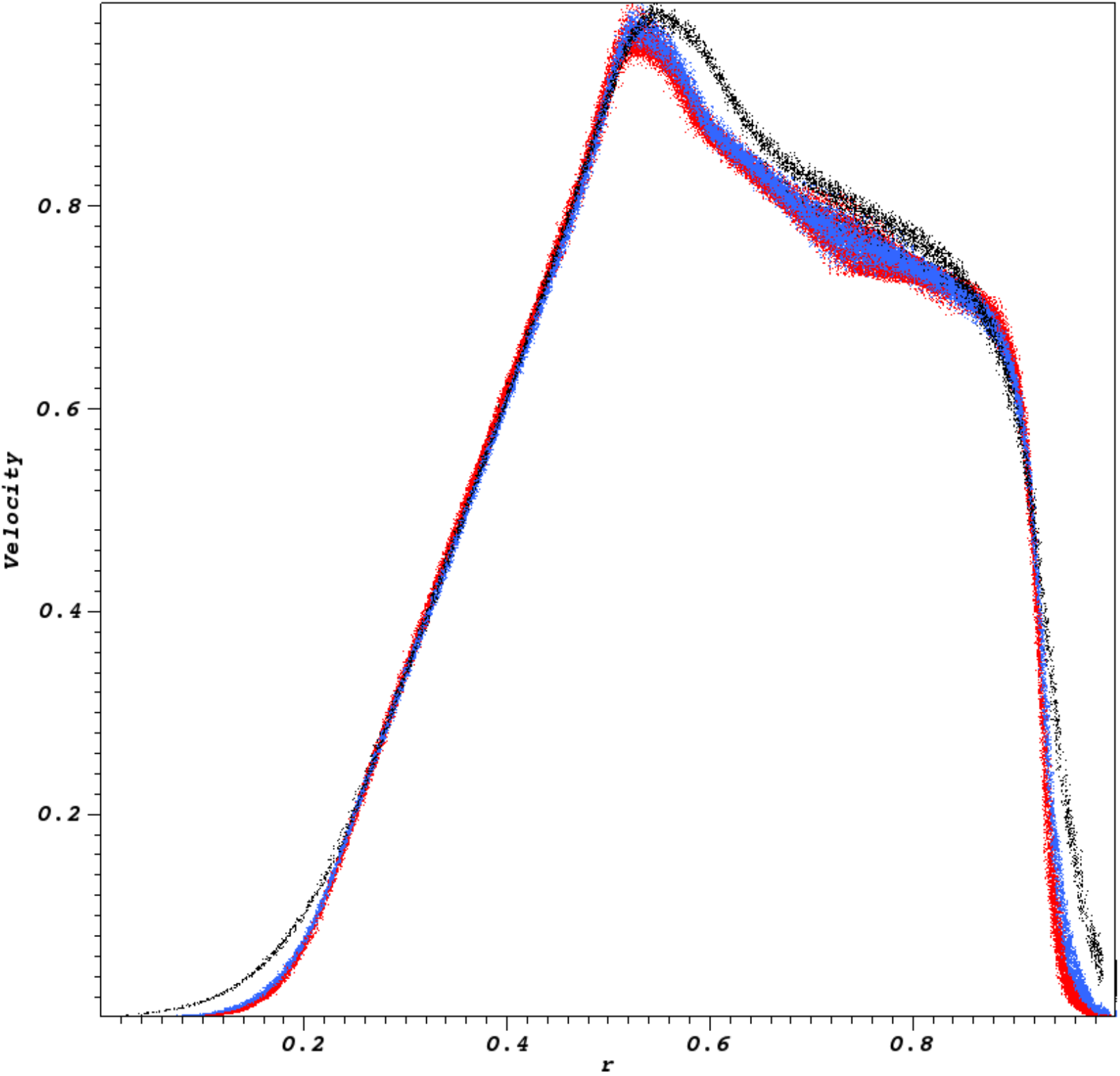}} 
\subfigure[Pressure]{\includegraphics[width=0.3\textwidth]{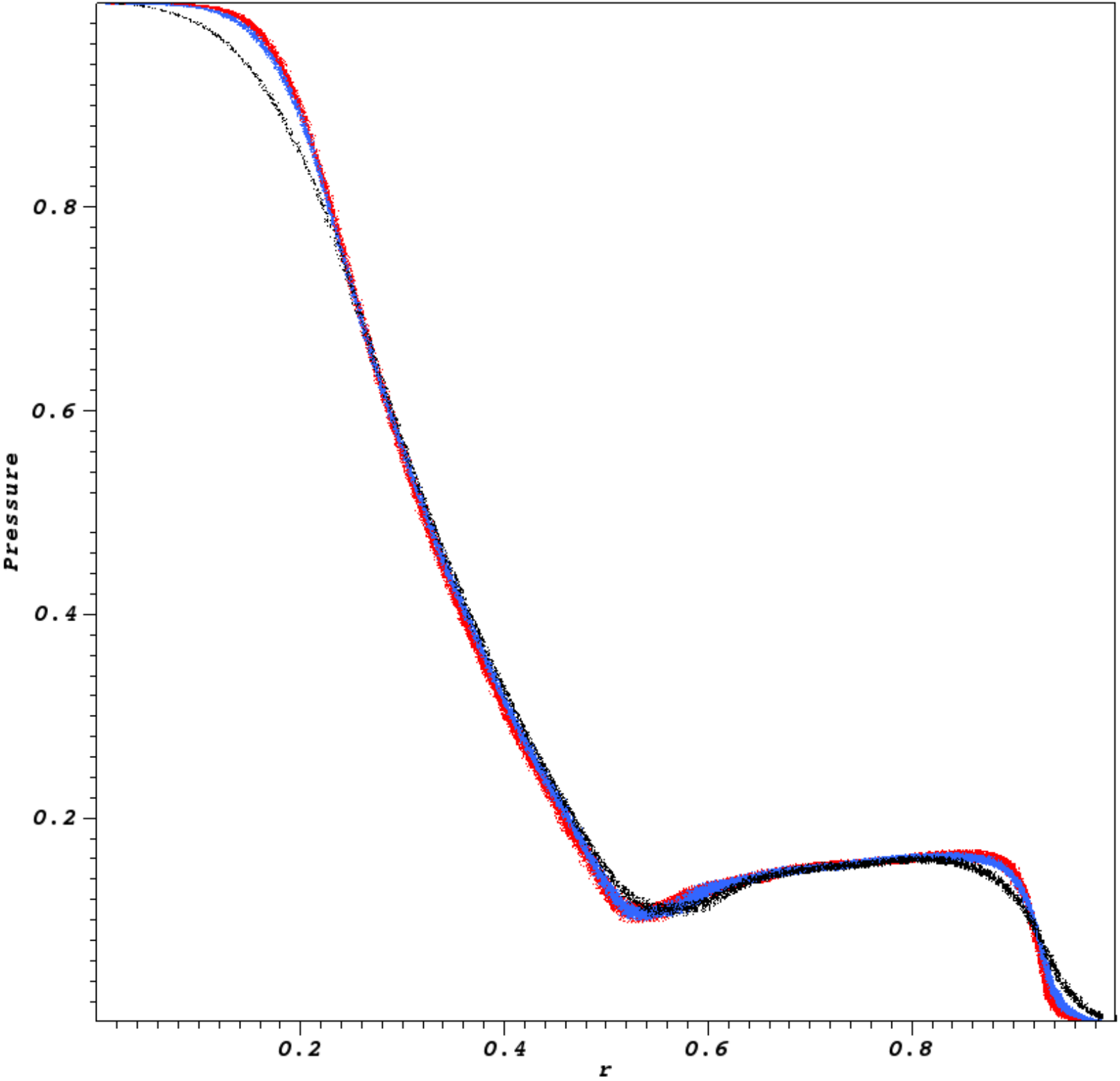}} 
\caption{2D Sod. Results for $B1$ (black), $B2$ (blue) and $B3$ (red) elements at $T=0.25$ on a grid with $N=13548$ elements.}\label{Sod_2D_scatter_fine}
\end{center}
\end{figure}

The computations have been performed a triangular mesh consisting of approximately \rev{$3500$ and $13500$}  elements\footnote{this last corresponds roughly to $115\times 115$ grid points on a Cartesian grid}. \rev{The stabilizing parameters of \eqref{phi_burman} have been set for B1, B2 and B3 as in Section \ref{Section:Vortex}. The 2D plot of the solution in Figure \ref{Sod_2D_coarse} and \ref{Sod_2D_fine}, as well as the behavior of the scatter plot of the density, velocity and pressure in Figure \ref{Sod_2D_scatter_coarse} and  \ref{Sod_2D_scatter_fine}} are in agreement with the previous test cases: increasing the order of the basis polynomials leads to the improvement of the solution quality and, thus, an enhanced representation of the shock waves.
\rev{Further, one may note also how increasing the mesh refinement, we can observe an increase of the symmetry of the solution from \ref{Sod_2D_scatter_coarse} to \ref{Sod_2D_scatter_fine}.}

\clearpage
\clearpage
\clearpage
\clearpage
\subsubsection{Mach 3 channel with forward-facing step}
\begin{figure}[H]
\begin{center}
\hspace{-0.1cm}\subfigure[Mesh with $N_0=2848$ elements]{\includegraphics[width=0.42\textwidth]{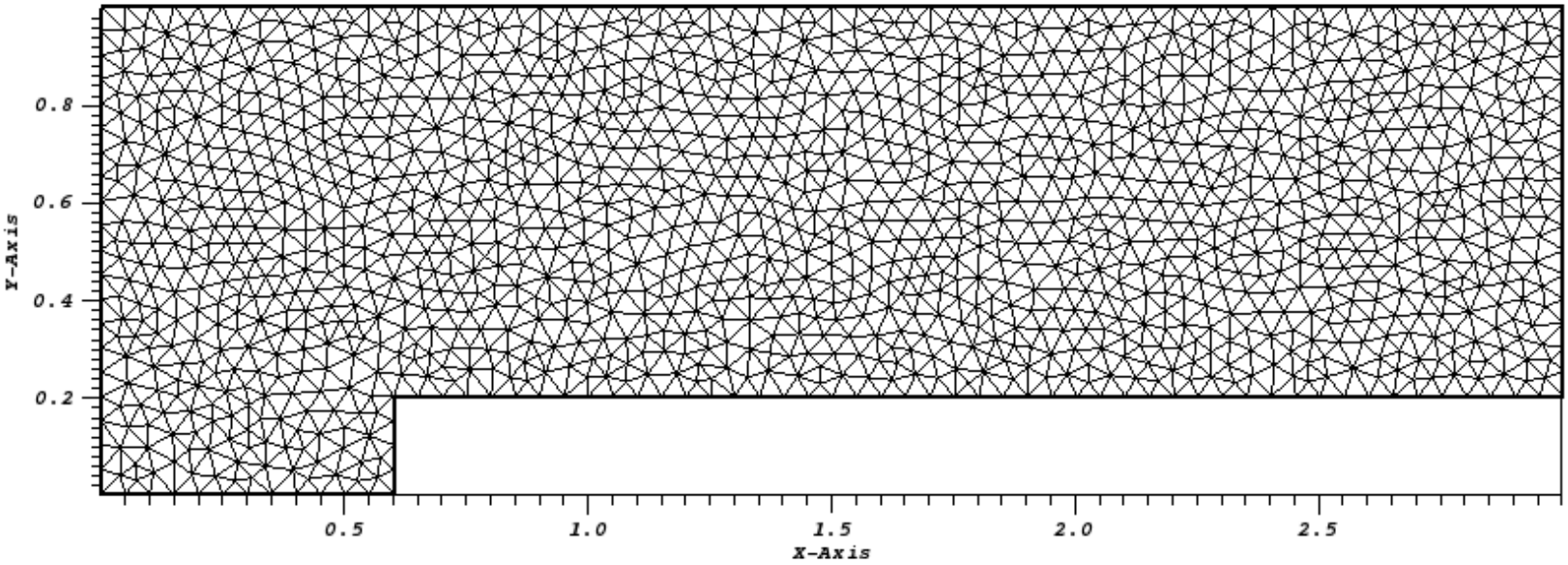}}
\hspace{0.75cm}\subfigure[Mesh with $N_1=11072$ elements]{\includegraphics[width=0.42\textwidth]{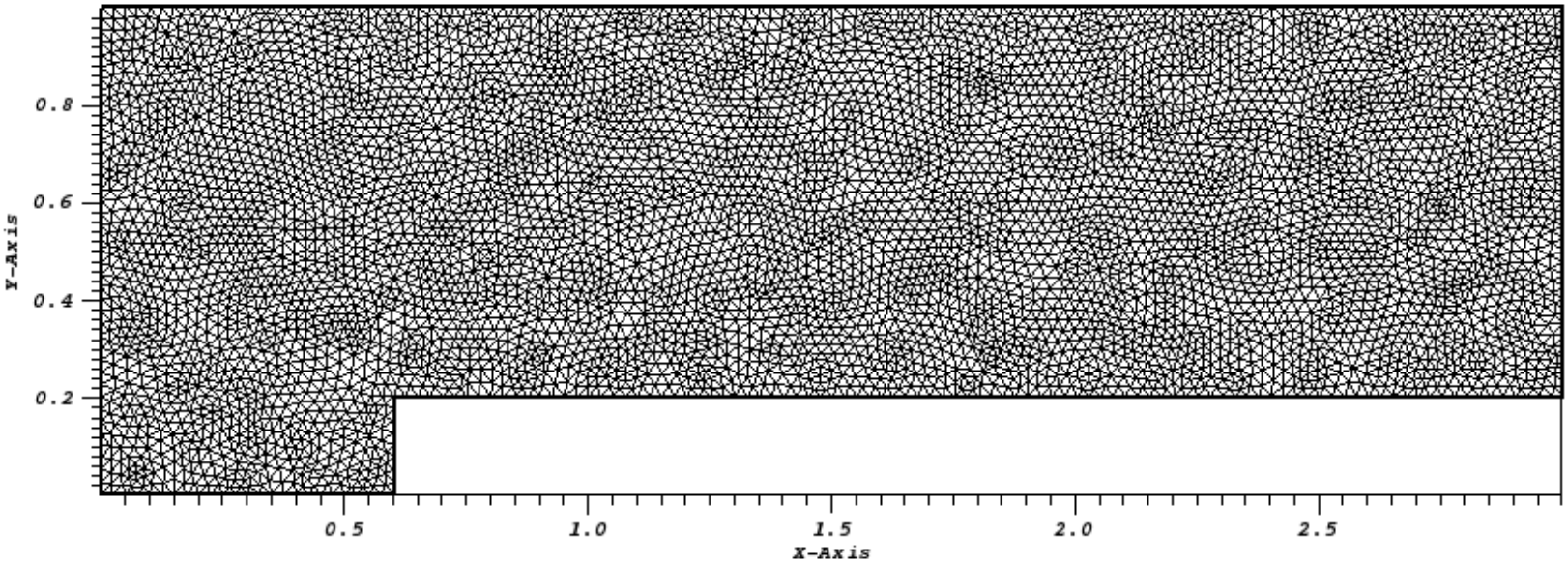}}\\
\hspace{-0.9cm} \subfigure[$B1$ with $N_0$]{\includegraphics[width=0.47\textwidth]{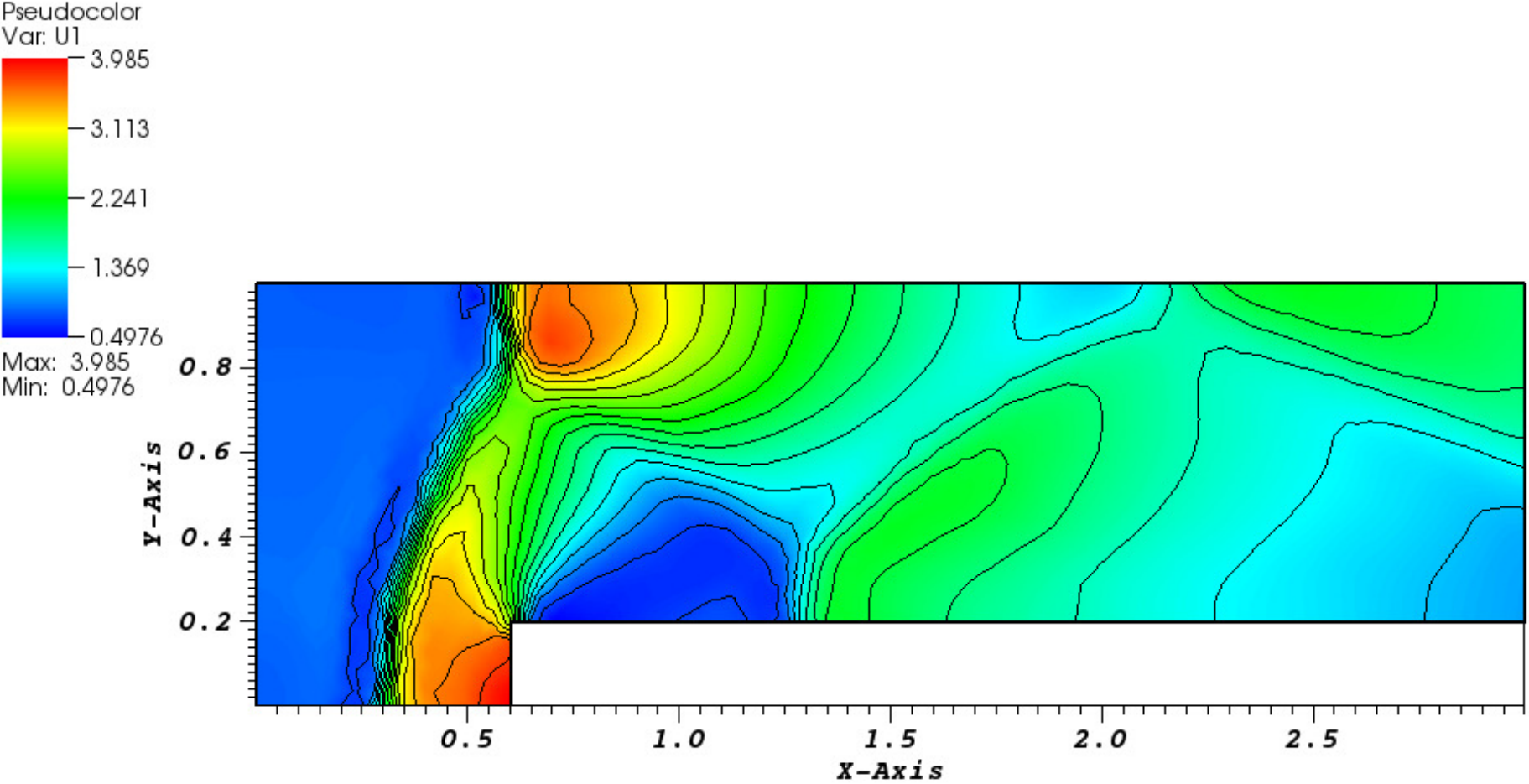}\label{step_2D_B1N0}}
\subfigure[$B1$ with $N_1$]{\includegraphics[width=0.475\textwidth]{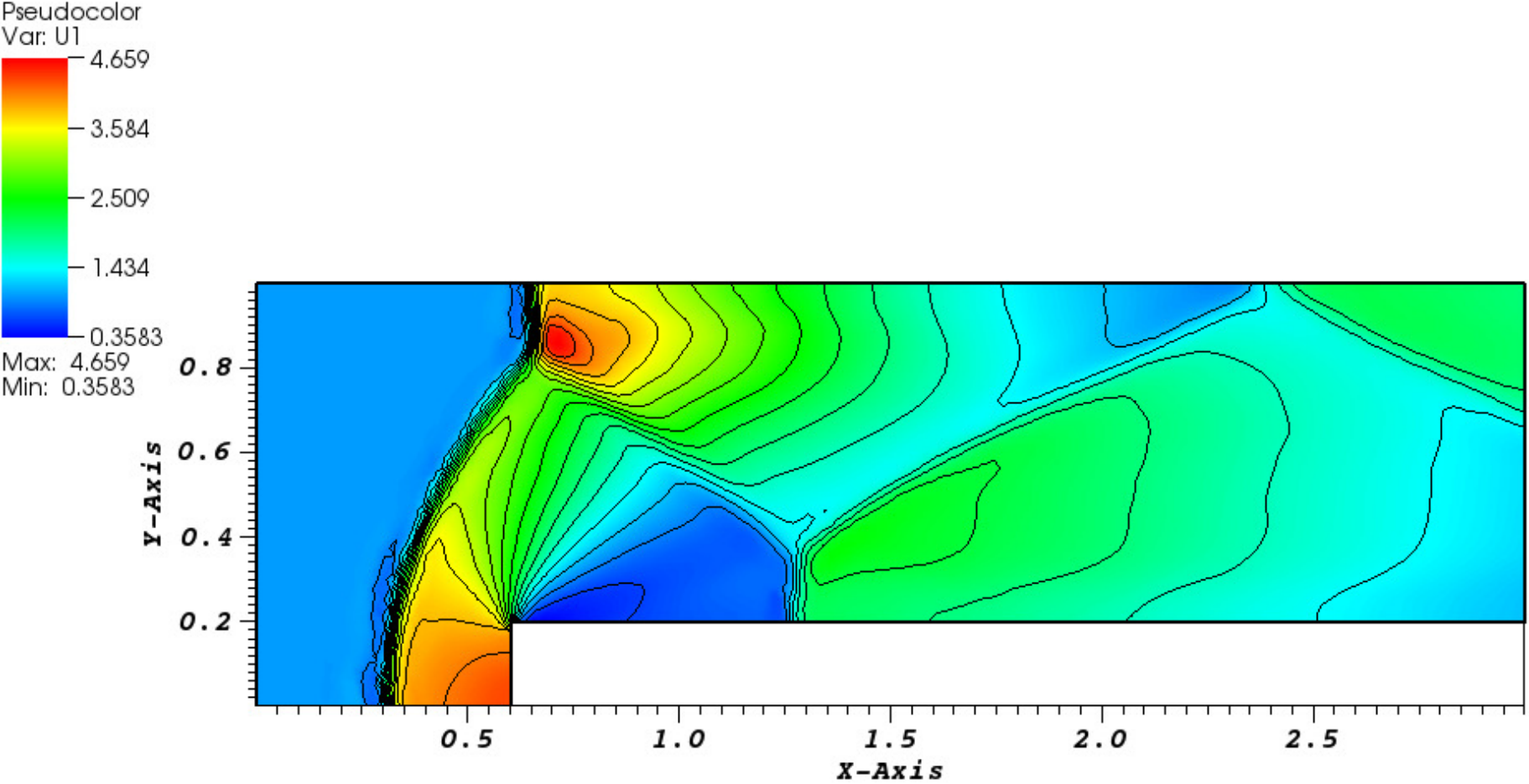}\label{step_2D_B1N1}}\\
\hspace{-0.9cm} \subfigure[$B2$ with $N_0$]{\includegraphics[width=0.47\textwidth]{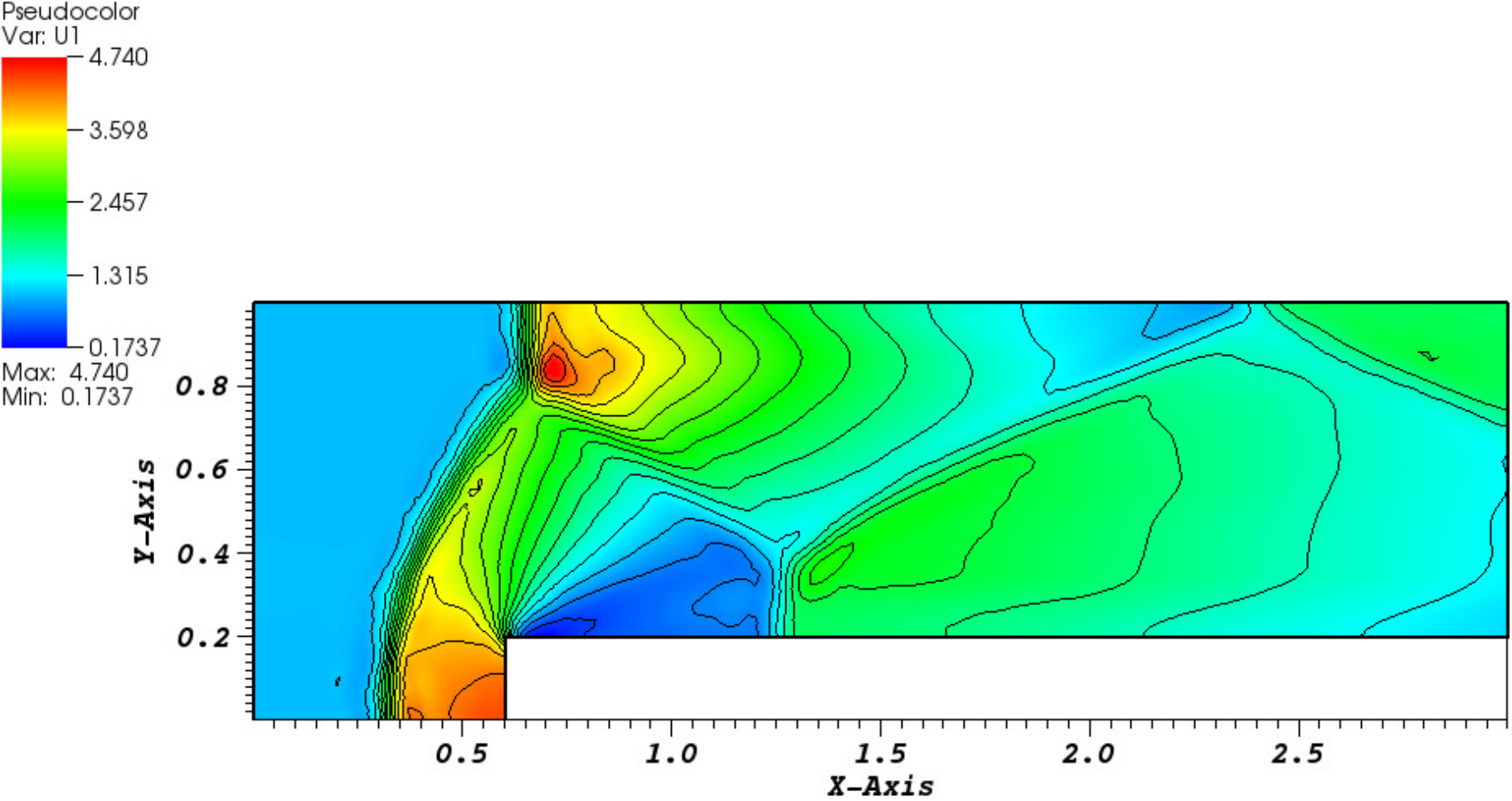}\label{step_2D_B2N0}} \subfigure[$B2$ with $N_1$]{\includegraphics[width=0.475\textwidth]{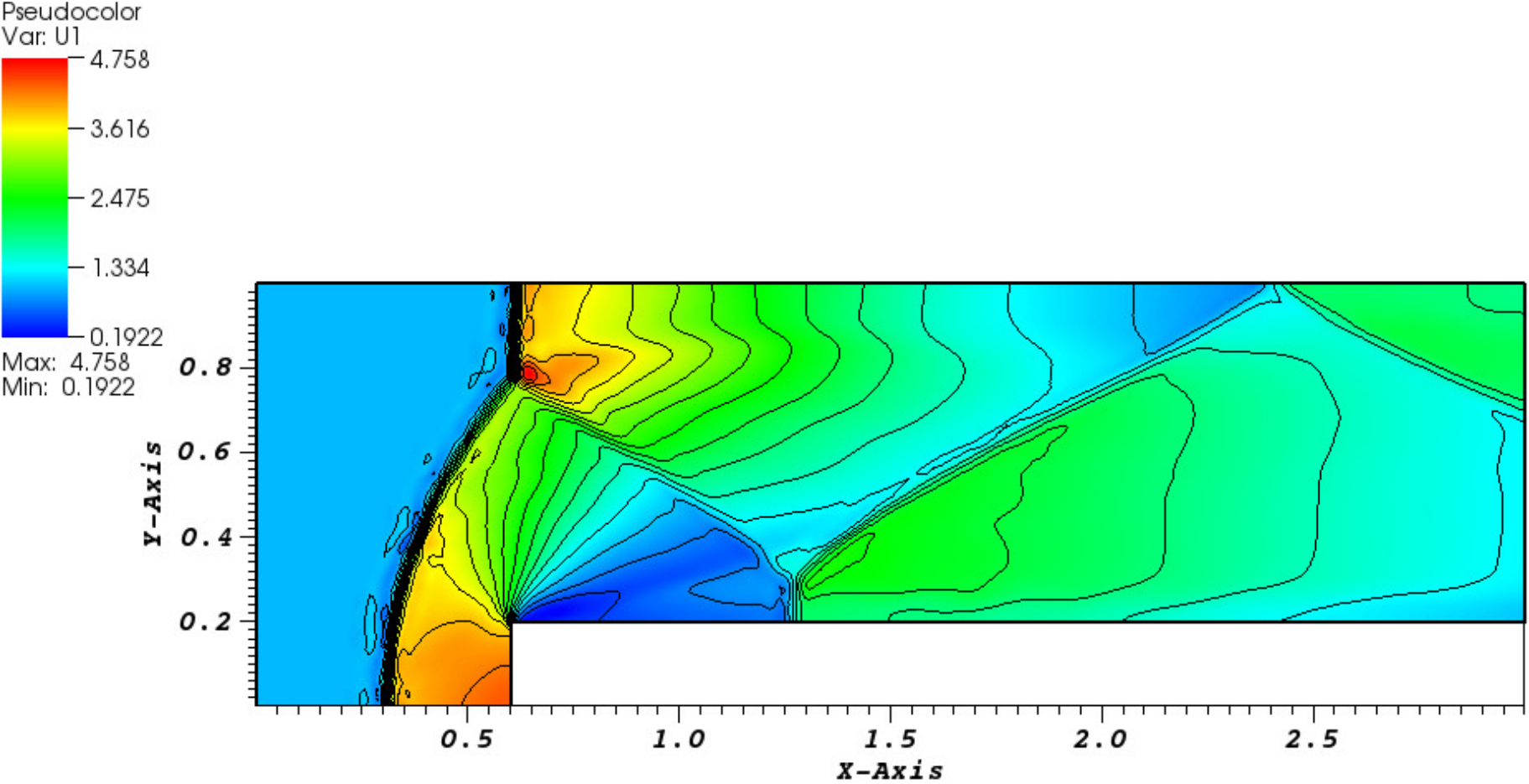}\label{step_2D_B2N1}}\\
\hspace{-0.9cm} \subfigure[$B3$ with $N_0$]{\includegraphics[width=0.47\textwidth]{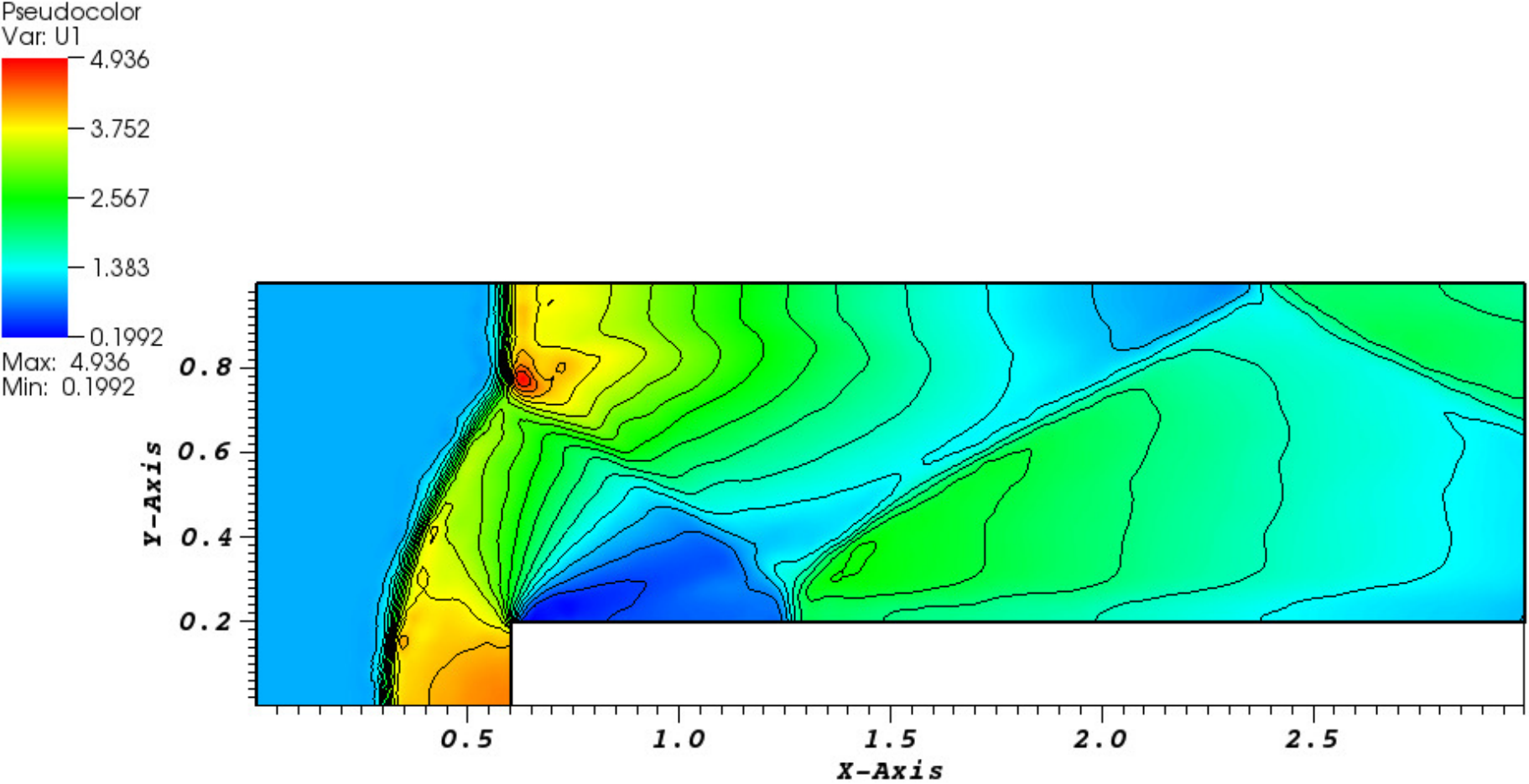}\label{step_2D_B3N0}} \subfigure[$B3$ with $N_1$]{\includegraphics[width=0.475\textwidth]{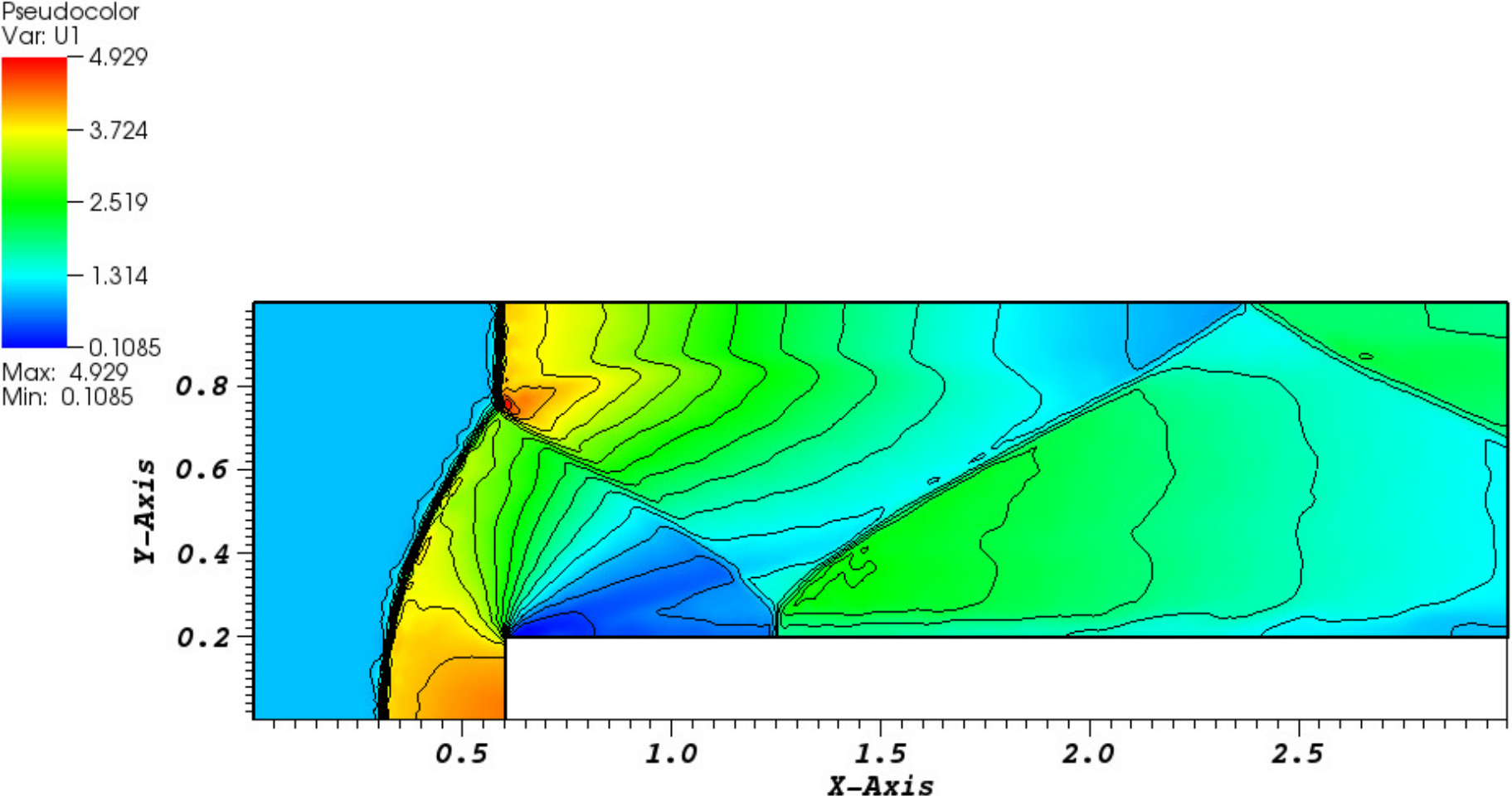}\label{step_2D_B3N1}}
\caption{Mach 3 channel with step. Results at $T=4.0$.}\label{step_2D}
\end{center}
\end{figure}

To assess the robustness of the proposed scheme in multidimensional problems involving strong shock waves, the Mach $3$ channel with a forward-facing step \cite{Woodward1984} test case has been used with $B1$, $B2$ and $B3$ elements on both coarse mesh having $N=2848$ cells\footnote{corresponds roughly to $30\times 100$ grid points} and finer one having $N=11072$ cells\footnote{corresponds roughly to $60\times 200$ grid points}) (see Fig.~\ref{step_2D}). \rev{The stabilizing parameters have been set for B1 to $\theta_1=0.1$ and $\theta_2=0$, for B2 to $\theta_1=0.3$ and $\theta_2=0$ and for B3 $\theta_1=0.05$ and $\theta_2=0$.} As expected, the quality of the solution increases when going from the second to fourth order scheme even on coarse meshes. Indeed, while in $B1$ case in Fig.~\ref{step_2D_B1N0} it is not possible to recognize the structure forming at the triple point, in $B2$ and $B3$ cases (Figs.~\ref{step_2D_B2N0} and \ref{step_2D_B3N0}, respectively) this structure  is already very well represented. On finer mesh, see Figs.~\ref{step_2D_B1N1},\ref{step_2D_B2N1} and \ref{step_2D_B3N1}, it is also possible to observe the gain in the quality of the approximation of shock waves when using a higher order RD method. 

\clearpage

\rev{
\subsubsection{Double Mach Reflection problem}
Finally, we present a widely used benchmark problem of a double Mack reflection problem as described in \cite{Woodward1984}. In this case, $B1$, $B2$ and $B3$ elements have been computed on a coarse mesh having $N=4908$ cells\footnote{corresponds roughly to $30\times 100$ grid points} and a finer one having $N=19248$ cells\footnote{corresponds roughly to $60\times 200$ grid points}) (see Fig.~\ref{Fig:DMR}). \rev{The stabilizing parameters have been set as in Section~\ref{Section:Vortex}.} Also here, as expected, the quality of the solution increases when going from the second to fourth order scheme on coarse meshes and more details are outlined on the finer mesh.
\begin{figure}[H]
\begin{center}
\hspace{0.3cm} \subfigure[Mesh with $N_0 = 4908$ elements]{\includegraphics[width=0.41\textwidth]{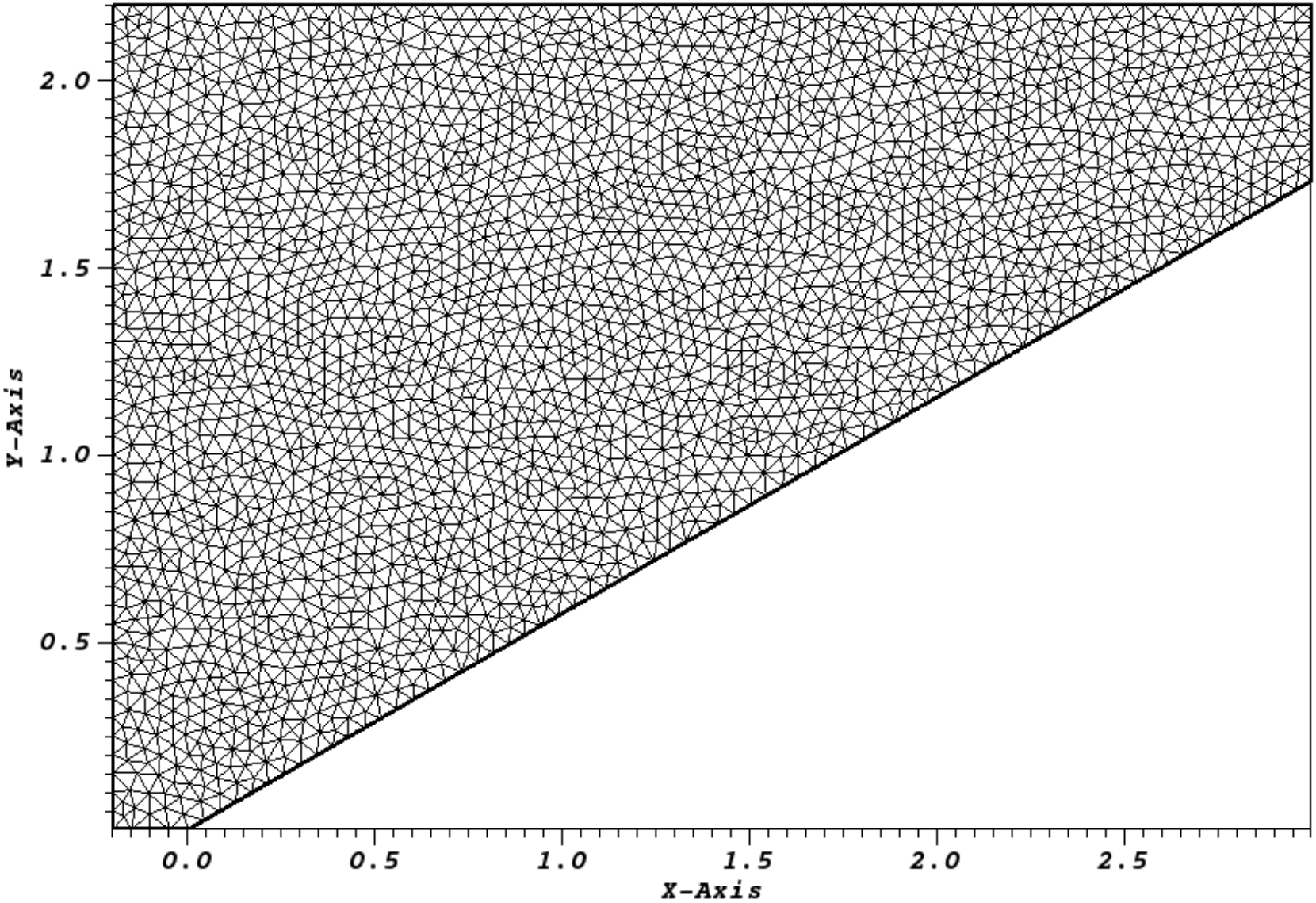}}
\hspace{0.4cm} \subfigure[Mesh with $N_1 = 19248$ elements]{\includegraphics[width=0.41\textwidth]{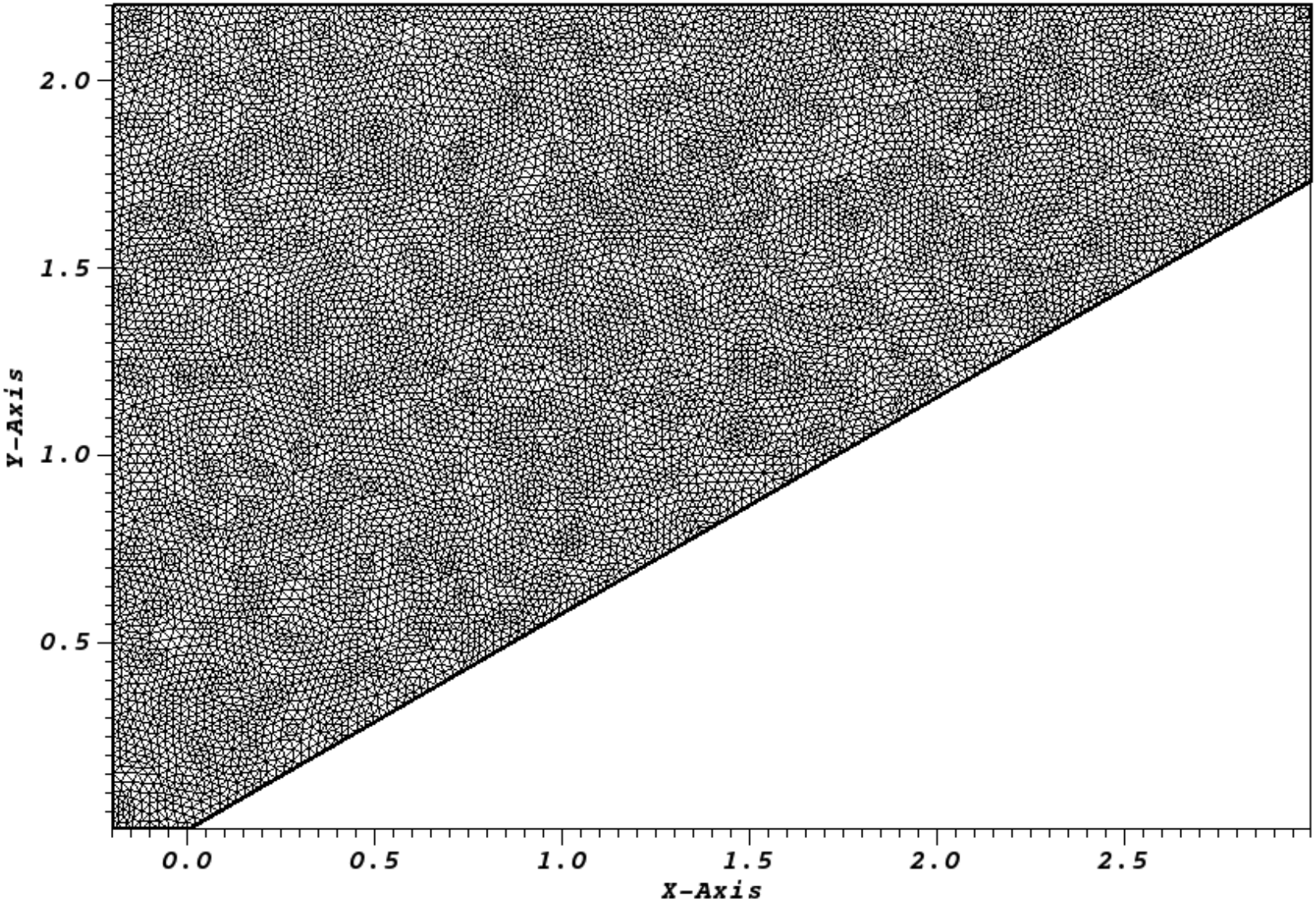}} \\
\subfigure[$B1$ with $N_0$]{\includegraphics[width=0.41\textwidth]{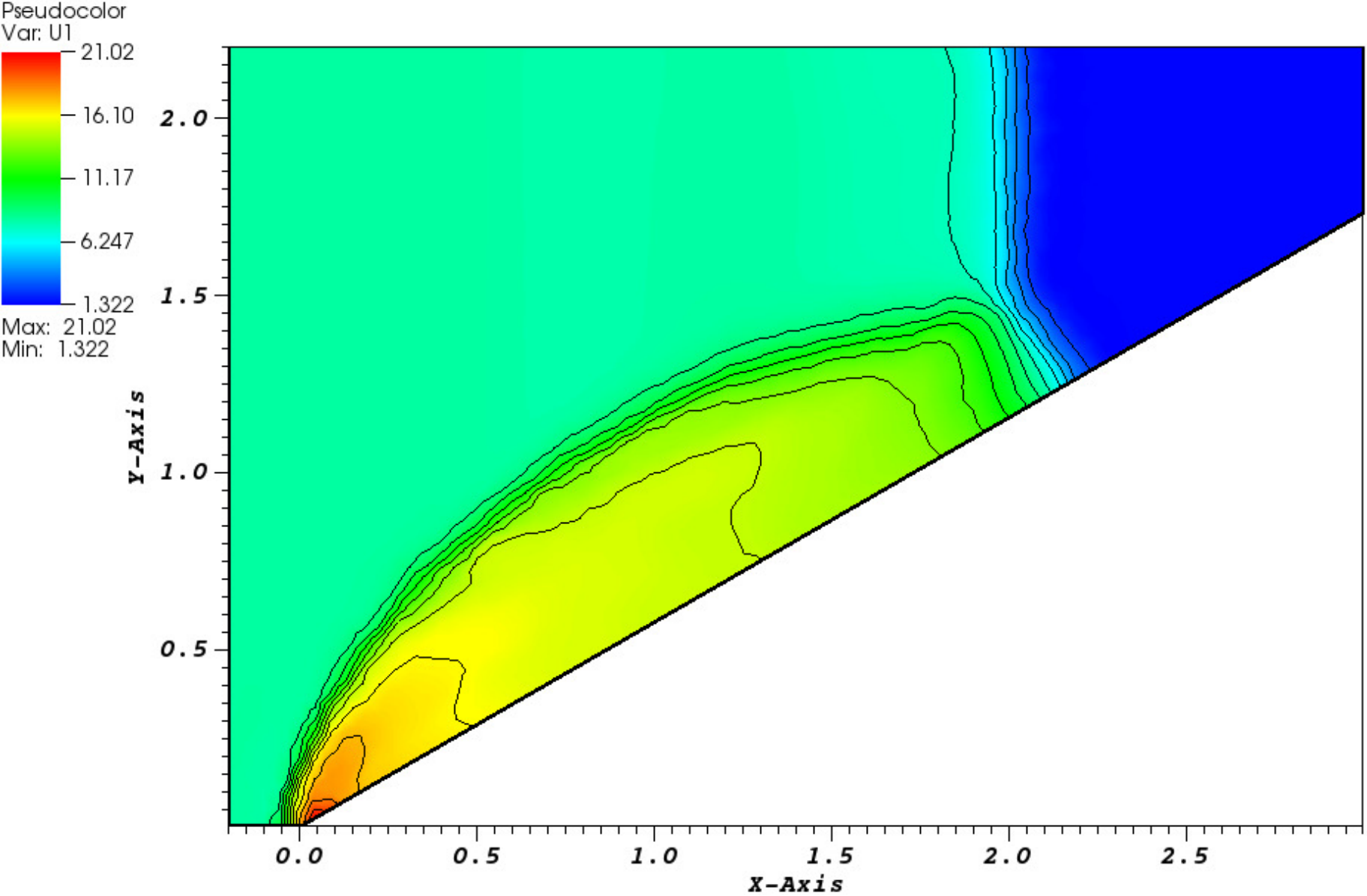}}
\hspace{0.6cm}\subfigure[$B1$ with $N_1$]{\includegraphics[width=0.41\textwidth]{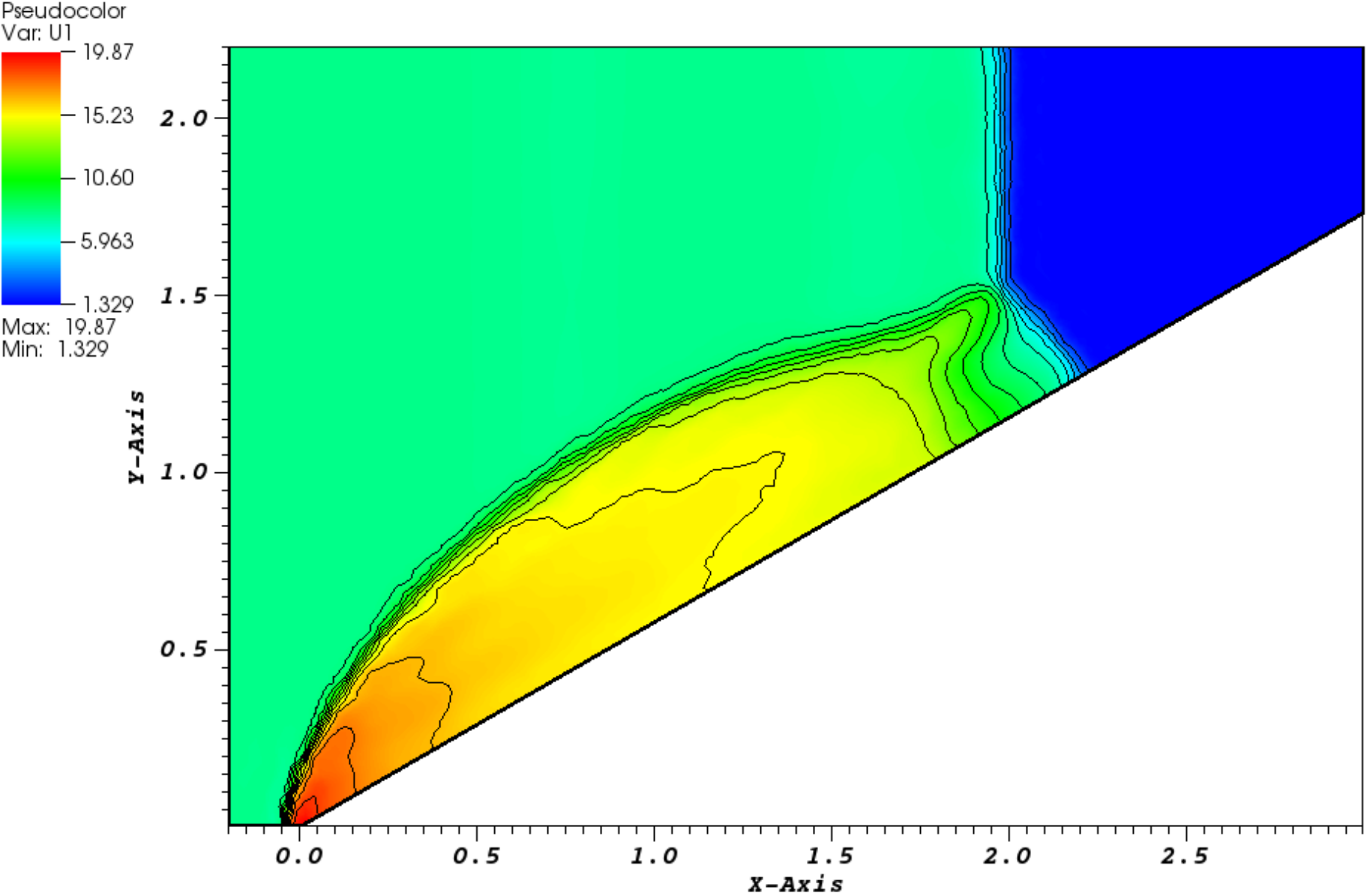}} \\
\subfigure[$B2$ with $N_0$]{\includegraphics[width=0.41\textwidth]{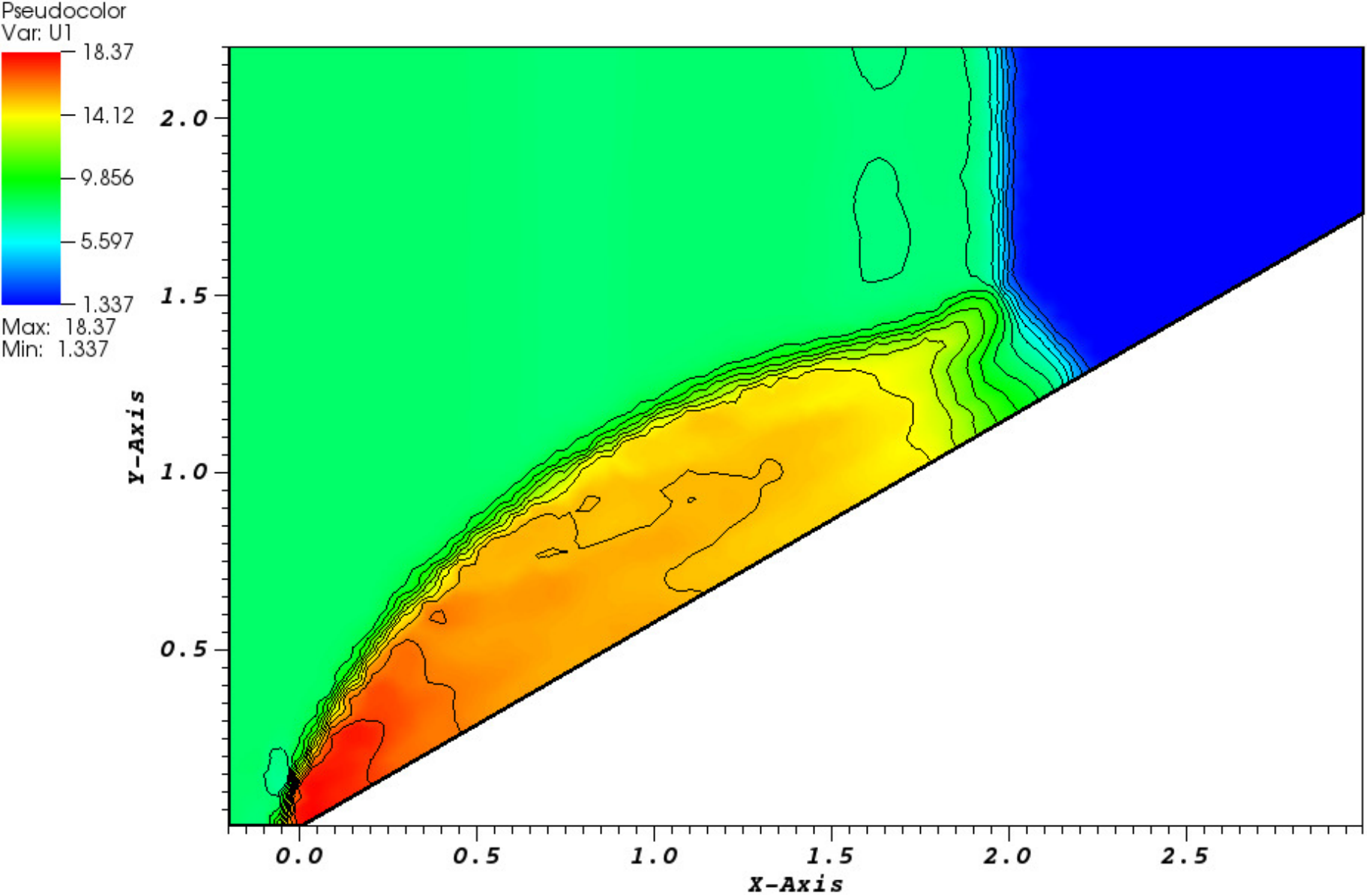}}
\hspace{0.6cm}\subfigure[$B2$ with $N_1$]{\includegraphics[width=0.41\textwidth]{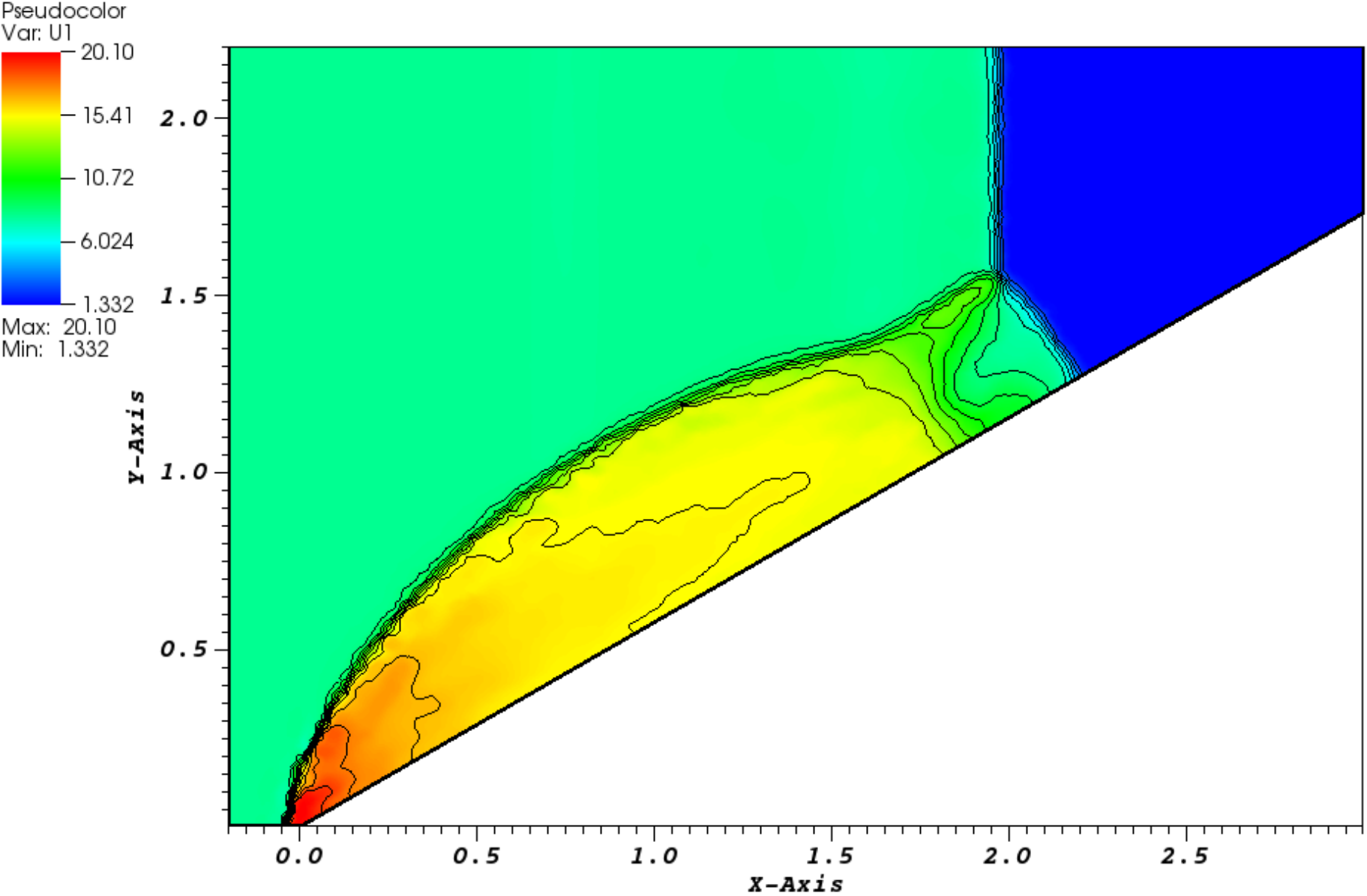}}\\
\subfigure[$B3$ with $N_0$]{\includegraphics[width=0.41\textwidth]{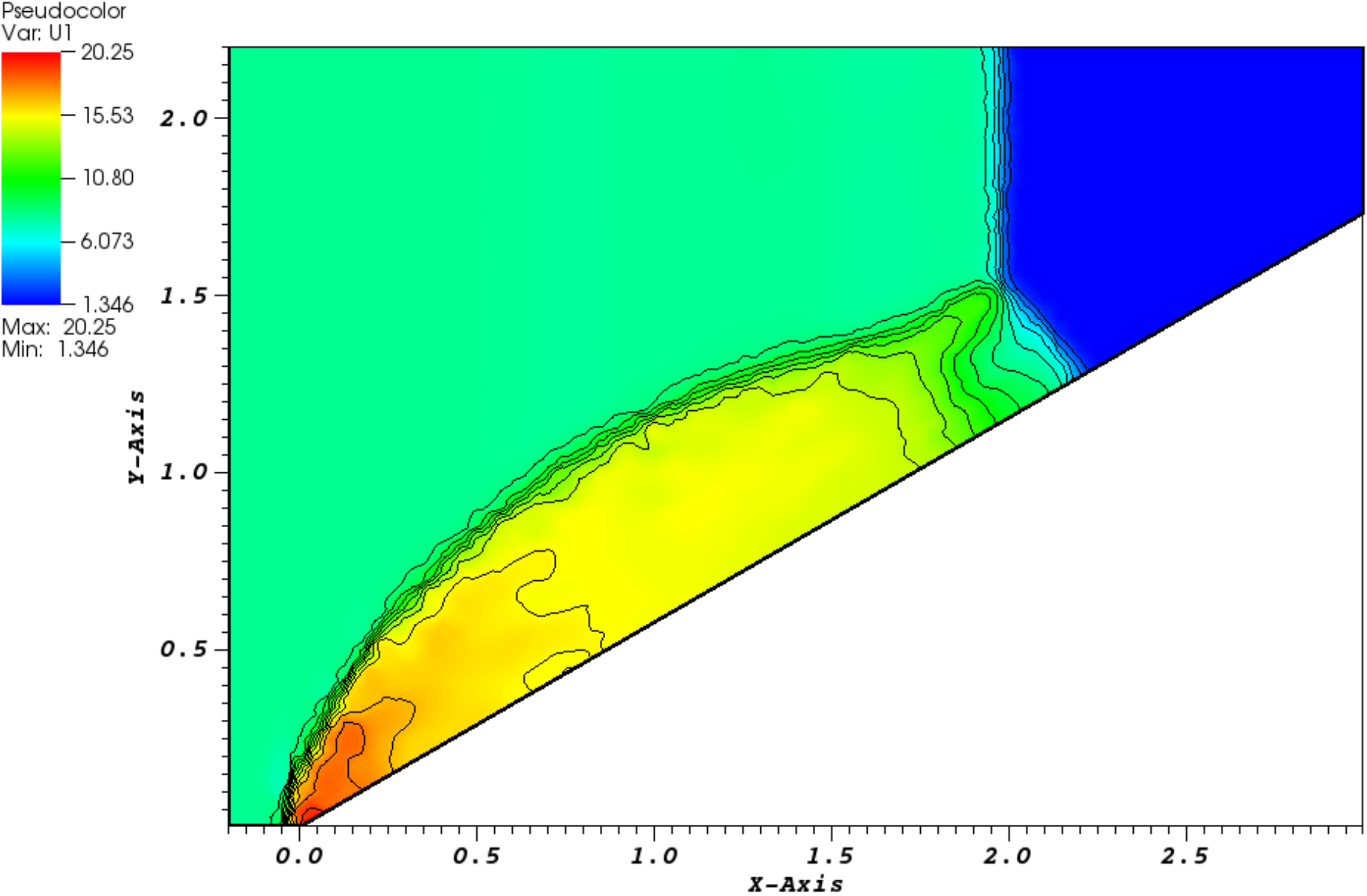}}
\hspace{0.6cm}\subfigure[$B3$ with $N_1$]{\includegraphics[width=0.41\textwidth]{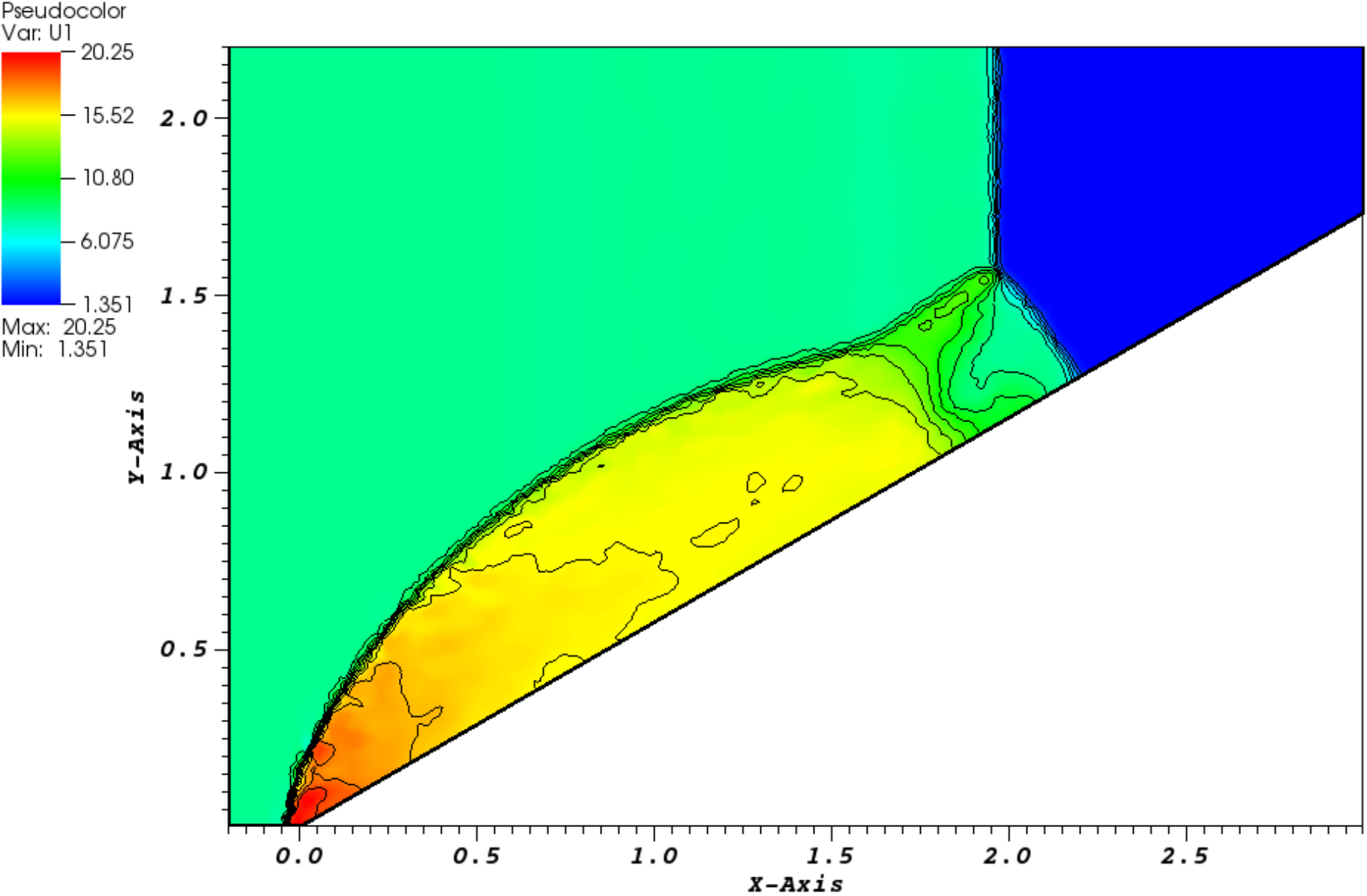}} 
\caption{\rev{Double Mach Reflection problem. Results at $T=0.2$}}
\label{Fig:DMR}
\end{center}
\end{figure}
}
\section{Conclusions}
In this paper, extending the ideas of \cite{mario} and \cite{Abgrall2017}, we have developed an explicit high order residual distribution scheme for multidimensional hyperbolic systems.
The main advantage of this approach consists in a simple way to detain high order of accuracy both in time and space, while having an explicit scheme in time which does not have the typical matrix "inversion" pertinent to classical finite element methodologies. 
We have demonstrated the mesh convergence to exact solutions with theoretically predicted high orders of accuracy by this new class of RD schemes.
Further, several benchmark problems have assessed the robustness of the scheme when dealing with strong discontinuities.
Extensions to other models, such as multiphase flows or Lagrangian hydrodynamics, and further investigations of high order residual distribution schemes will be considered in  forthcoming papers. \rev{Finally, we are currently extending the proposed approach to viscous problems by combining it with the discretisation technique as explained in \cite{DeSantis2015}. In this case, the time step results, of course, very small, and thus an implicit approach is needed with the challenge to have, nevertheless, a diagonal 'mass matrix'.}

\section*{Acknowledgements}
ST and PB have been funded by SNSF project 200021\_153604 "High fidelity simulation for compressible materials".  R.A. has been funded in part by the same project. The authors would like to thank M. Han Veiga and D. Torlo for their useful feedbacks on this paper. \rev{A special acknowledgement goes to D. Kuzmin who made very useful suggestions to lower the numerical dissipation of the Lax-Friedrichs scheme.}

\appendix
\section{\rev{A remark on the total residual as a weighted sum on triangular elements}}\label{appendix A}
\rev{
For a better understanding of \eqref{divkinto1} and following \cite{AbgrallViville2017}, we rewrite explicitly how the total residual is obtained in 2D on triangular elements in terms of a weighted sum of first order residuals.
This is carried out for the case of a B2 approximation.\\
Let us define first the nomenclature of the considered element to be as the left triangle of Figure \ref{Fig:sub_elements} with $k$ in \eqref{divkinto1} to be $2$. Take $\sigma$ as the generic degree of freedom on the considered element $K$, while $\sigma'$ the degree of freedom of the sub-triangles $K_i$ with $i=1,..,4$.
Further, let $\varphi$ be the Bernstein shape function of order 2,  $x_j$ the barycentric coordinates corresponding to the vertices of the triangular element with $j=1,2,3$, and $(y_i^j)_{i=1,2,3}$ the  local barycentric coordinates of the element $K_j$.
We define the first order approximation of the Flux as in \eqref{divkinto1} to read 
$$\bF_{K_i}^{(1)}=\sum\limits_{\sigma'_j \text{ vertex of } K_i} \bF_{\sigma'}y_i^j.$$
We set for the B2 approximation
$$\bF^{(2)}=\sum_{i=1}^{6} \bF_i\varphi_i.$$
The gradient of the Bernstein basis functions, which have been defined at the end of Section \ref{Sec_L1operator} are defined as 
\begin{equation}
\begin{split}
&\nabla \varphi_j=2x_j\nabla x_j,\,\, \text{for}\,\, j=1,2,3\\
&\nabla \varphi_4=2(x_1\nabla x_2+x_2\nabla x_1),\\
& \nabla \varphi_5=2(x_2\nabla x_3+x_3\nabla x_2),\\
& \nabla \varphi_5=2(x_1\nabla x_3+x_3\nabla x_1).
\end{split}
\end{equation}
Rewriting the total residual explicitly will result in the following
\begin{equation}
\begin{split}
\int_K \varphi_\sigma \text{ div }\bF^{(2)}\; d\mathbf{x} =\,&\sum_{K_i \in K, \sigma'
\in K_i} \big(\int_K \varphi_\sigma \text{ div } \varphi_{\sigma'} \big) \bF_{\sigma'} \\=\, &
\quad \, 2\int_K \varphi_{\sigma} x_1\nabla x_1 \bF_1+2\big( \int_K  \varphi_{\sigma} x_2 \nabla x_1 + \int_K  \varphi_{\sigma}  x_1 \nabla x_2 \big) \bF_4\\&+
2\int_K \varphi_{\sigma} x_2\nabla x_2 \bF_2+2\big( \int_K  \varphi_{\sigma} x_2 \nabla x_3 + \int_K  \varphi_{\sigma}  x_3 \nabla x_2 \big) \bF_5\\&+
2\int_K \varphi_{\sigma} x_3\nabla x_3 \bF_3+2\big( \int_K  \varphi_{\sigma} x_1 \nabla x_3 + \int_K  \varphi_{\sigma}  x_3 \nabla x_1 \big) \bF_6,
\end{split}
\end{equation}
which in terms of the sum over the sub-triangles is equivalent to
\begin{equation}
\begin{split}
\int_K \varphi_\sigma \text{ div }\bF^{(2)}\; d\mathbf{x} =\,&
\quad \, 2\int_K \varphi_{\sigma} x_1\nabla x_1 \bF_1+2\int_K \varphi_{\sigma}  x_1 \nabla x_2 \bF_4 + 2\int_K \varphi_{\sigma}  x_1 \nabla x_3 \bF_6\quad  \quad \quad (K_1)\\&+
2\int_K \varphi_{\sigma} x_2\nabla x_2 \bF_2+2\int_K \varphi_{\sigma}  x_2 \nabla x_1 \bF_4 + 2\int_K \varphi_{\sigma}  x_2 \nabla x_3 \bF_5\quad  \quad \quad (K_2)\\&+
2\int_K \varphi_{\sigma} x_3\nabla x_3 \bF_3+2\int_K \varphi_{\sigma}  x_3 \nabla x_1 \bF_6 + 2\int_K \varphi_{\sigma}  x_3 \nabla x_2 \bF_5\quad  \quad \quad (K_3)\\&-
2\int_K \varphi_{\sigma} x_3\nabla x_3 \bF_4-2\int_K \varphi_{\sigma} x_1\nabla x_1 \bF_5-2\int_K \varphi_{\sigma} x_2\nabla x_2 \bF_6 \quad  \quad \quad (K_4)\\=\,&
\quad \,  \sum_{i=1}^4 2 \int_K \varphi_{\sigma} \text{ div }\bF^{(1)}\; d\mathbf{x}
\end{split}
\label{explicit_weightedsum}
\end{equation}
From \eqref{explicit_weightedsum} one can easily see the relation with \eqref{divkinto1}, where in particular $$\omega_{K_i}=2 \int_K \varphi_{\sigma} d\mathbf{x}.$$
For Bernstein shape functions of order 3 the same idea is applied.
}

\clearpage

\section*{References}
\bibliography{biblio}

\begin{thebibliography}{10}
\expandafter\ifx\csname url\endcsname\relax
  \def\url#1{\texttt{#1}}\fi
\expandafter\ifx\csname urlprefix\endcsname\relax\def\urlprefix{URL }\fi
\expandafter\ifx\csname href\endcsname\relax
  \def\href#1#2{#2} \def\path#1{#1}\fi

\bibitem{mario}
M.~Ricchiuto, R.~Abgrall, Explicit {Runge-Kutta} residual distribution schemes
  for time dependent problems: Second order case, Journal of Computational
  Physics 229~(16) (2010) 5653--5691.

\bibitem{abg}
R.~Abgrall, Residual distribution schemes: Current status and future trends,
  Computers and Fluids 35~(7) (2006) 641--669.

\bibitem{shu-dec}
Y.~Liu, C.-W. Shu, M.~Zhang, Strong stability preserving property of the
  deferred correction time discretisation, Journal of Computational Mathematics
  26~(5) (2008) 633--656.

\bibitem{Minion2}
M.~Minion, Semi-implicit spectral deferred correction methods for ordinary
  differential equaions, Communication in Mathematical Physics 1~(3) (2003)
  471--500.

\bibitem{enumath}
R.~{Abgrall}, P.~{Bacigaluppi}, S.~{Tokareva}, {How to avoid mass matrix for
  linear hyperbolic problems.}, in: {Numerical mathematics and advanced
  applications -- ENUMATH 2015. Selected papers based on the presentations at
  the European conference, Ankara, Turkey, September 14--18, 2015}, Cham:
  Springer, 2016, pp. 75--86.
\newblock \href {http://dx.doi.org/10.1007/978-3-319-39929-4_8}
  {\path{doi:10.1007/978-3-319-39929-4_8}}.

\bibitem{Abgrall2017}
R.~Abgrall, High order schemes for hyperbolic problems using globally
  continuous approximation and avoiding mass matrices, Journal of Scientific
  Computing~(2) (2017) 461--494.

\bibitem{Abgrall2010}
R.~Abgrall, J.~Trefil{\'i}k, {An Example of High Order Residual Distribution
  Scheme Using non-Lagrange Elements}, Journal of Scientific Computing 45~(1)
  (2010) 3--25.

\bibitem{Lohmann2017}
C.~Lohmann, D.~Kuzmin, J.~N. Shadid, S.~Mabuza, {Flux-corrected transport
  algorithms for continuous Galerkin methods based on high order Bernstein
  finite elements}, {Journal of Computational Physics} 344 (2017) 151--186.

\bibitem{enordhigh}
R.~Abgrall, A.~Larat, M.~Ricchiuto, Construction of very high order residual
  distribution schemes for steady inviscid flow problems on hybrid unstructured
  meshes, Journal of Computational Physics 230~(11) (2011) 4103--4136.

\bibitem{Cangiani2013}
A.~Cangiani, J.~Chapman, E.~H. Georgoulis, M.~Jensen, {Implementation of the
  Continuous-Discontinuous Galerkin Finite Element Method}, in: A.~Cangiani,
  R.~L. Davidchack, E.~Georgoulis, A.~N. Gorban, J.~Levesley, M.~V. Tretyakov
  (Eds.), Numerical Mathematics and Advanced Applications 2011, Springer Berlin
  Heidelberg, Berlin, Heidelberg, 2013, pp. 315--322.

\bibitem{AbgrallViville2017}
R.~Abgrall, Q.~Viville, H.~Beaugendre, C.~Dobrzynski, {Construction of a
  p-Adaptive Continuous Residual Distribution Scheme}, J. Sci. Comput. 72~(3)
  (2017) 1232--1268.

\bibitem{SWjcp}
M.~Ricchiuto, R.~Abgrall, H.~Deconinck, Application of conservative residual
  distribution schemes to the solution of the shallow water equations on
  unstructured meshes, Journal of Computational Physics 222 (2007) 287--331.

\bibitem{Ricchiuto2007}
H.~Deconinck, M.~Ricchiuto, Encyclopedia of Computational Mechanics, John Wiley
  \& Sons, 2007, Ch. Residual Distribution Schemes: Foundations and Analysis.

\bibitem{abgrall2017jcp}
R.~Abgrall, {Some remarks about conservation for residual distribution
  schemes}, Journal of Computational Physics, submitted.
\newblock \href {http://dx.doi.org/arXiv:1708.03108}
  {\path{doi:arXiv:1708.03108}}.

\bibitem{Abgrall99}
R.~Abgrall, Toward the ultimate conservative scheme: Following the quest., J.
  Comput. Phys. 167~(2) (2001) 277--315.

\bibitem{hughes}
T.~J. {Hughes}, M.~{Mallet}, A new finite element formulation for computational
  fluid dynamics: {III}. the generalized streamline operator for
  multidimensional advective-diffusive systems, Computer Methods in Applied
  Mechanics and Engineering 58~(3) (1986) 305--328.

\bibitem{burman}
E.~Burman, P.~Hansbo., Edge stabilization for {Galerkin} approximations of
  convection-diffusion-reaction problems, Computational Methods in Applied
  Mechanical Engineering 193 (2004) 1437--1453.

\bibitem{DeSantis2015}
R.~{Abgrall}, D.~{De Santis}, {Linear and non-linear high order accurate
  residual distribution schemes for the discretization of the steady
  compressible Navier-Stokes equations.}, {J. Comput. Phys.} 283 (2015)
  329--359.
\newblock \href {http://dx.doi.org/10.1016/j.jcp.2014.11.031}
  {\path{doi:10.1016/j.jcp.2014.11.031}}.

\bibitem{icm}
R.~Abgrall, On a class of high order schemes for hyperbolic problems, in: S.~Y.
  Jang, Y.~Kim, D.~Lee, I.~Yie (Eds.), Proceedings of the International
  Conference of Mathematicians, Vol.~IV, Kyung Moon SA Co. Ltd., 2014, pp.
  699--725.

\bibitem{CanadaCFD}
R.~Abgrall, {About Non Linear Stabilization for Scalar Hyperbolic Problems},
  in: J.~B. Roderick~Melnik, Roman~Makarov (Ed.), Recent Progress and Modern
  Challenges in Applied Mathematics, Modeling and Computational Science, Fields
  Institute Communications, Springer, 2017, pp. 89--116.

\bibitem{ENORD}
R.~Abgrall, Essentially non oscillatory residual distribution schemes for
  hyperbolic problems, Journal of Computational Physics 214~(2) (2006)
  773--808.

\bibitem{Morton}
R.~D. Richtmyer, K.~W. Morton, Difference Methods for Initial-Value Problems,
  Inter-science, New-York, 1967.

\bibitem{abg2001d}
R.~Abgrall, P.~L. Roe, High order fluctuation schemes on triangular meshes,
  Journal of Scientific Computing 19~(1) (2003) 3--36.

\bibitem{Jameson1995}
A.~Jameson, Positive schemes and shock modelling for compressible flows,
  International Journal for Numerical Methods in Fluids 20 (1995) 743--776.

\bibitem{ChengShu2014}
J.~Cheng, C.-W. Shu, Positivity-preserving {L}agrangian scheme for
  multi-material compressible flow, Journal of Computational Physics 257 (2014)
  143--168.

\bibitem{shuOsher1989}
C.~Shu, S.~Osher, Efficient implementation of essentially non-oscillatory
  shock-capturing schemes, {II}, Journal of Computational Physics 83 (1989)
  32--78.

\bibitem{Yee1999}
H.~Yee, N.~Sandham, M.~Djomehri, Low-dissipative high-order shock-capturing
  methods using characteristic-based filters, Journal of Computational Physics
  150~(1) (1999) 199 -- 238.

\bibitem{Barrenechea}
G.~Barrenechea, P.~Knobloch, Analysis of a group finite element formulation,
  Applied Numerical Mathematics 118 (2017) 238--248.

\bibitem{Woodward1984}
P.~{Woodward}, P.~{Colella}, The numerical simulation of two-dimensional fluid
  flow with strong shocks, Journal of Computational Physics 54 (1984) 115--173.

\end{thebibliography}
\bibliographystyle{elsarticle-num}

\end{document}